%% file: main.tex
\documentclass[a4paper, 11pt]{article}
\usepackage[utf8]{inputenc}

\RequirePackage{etex} 
\usepackage{authblk}
\usepackage{blindtext}
\usepackage{comment}
\usepackage[utf8]{inputenc}
\usepackage{amsmath}
\usepackage{multirow}
\usepackage{amsthm}
\usepackage{amsfonts}
\usepackage{amssymb}
\usepackage{empheq}
\usepackage{mathtools}
\usepackage{amsbsy}
\usepackage{braket}
\usepackage[font=footnotesize]{caption}
\usepackage[export]{adjustbox}
\usepackage[shortlabels]{enumitem}
\usepackage{xcolor}
\usepackage{multicol}
\usepackage{graphicx}
\usepackage[english]{babel}
\usepackage[margin=0.8in]{geometry}

\usepackage{float}
\usepackage{longtable}
\usepackage{subcaption}
\usepackage{listings}
\usepackage{hyperref}
\usepackage{cleveref}
\usepackage{autonum}
\usepackage{multirow}
\usepackage{siunitx}
\usepackage{tikz}
\usepackage{pgfplots}
\pgfplotsset{compat=1.18} 
\usepackage{enumitem}
\usepackage[makeroom]{cancel}
\usepackage[square, numbers, sort&compress]{natbib} 
\theoremstyle{definition}

\numberwithin{definition}{section}

\newtheorem{proposition}{Proposition}
\numberwithin{proposition}{section}
\newtheorem{assumption}{Assumption}
\numberwithin{assumption}{section}

\newtheorem{theorem}{Theorem}
\numberwithin{theorem}{section}

\theoremstyle{plain}
\newtheorem{remark}
{Remark}
\theoremstyle{plain}

\newcommand{\jump}[1]{\left[\mkern-1.5mu \left[#1\right] \mkern-1.5mu\right]}
\newcommand{\avg}[1]{\left\{ \mkern-5mu \left\{#1 \right\} \mkern-5mu \right\}}
\newcommand{\wavg}[1]{\{ \mkern-5mu \{#1 \} \mkern-5mu \}_{\omega}}
\newcommand{\wavgmu}[1]{\{ \mkern-5mu \{#1 \} \mkern-5mu \}_{\omega_{\mu}}}
\newcommand{\wavgd}[1]{\{ \mkern-5mu \{#1 \} \mkern-5mu \}_{\omega_{\delta_1}}}
\newcommand{\wavgD}[1]{\{ \mkern-5mu \{#1 \} \mkern-5mu \}_{\omega_{\mathbf{D}}}}

\renewcommand{\div}[1]{\nabla \mkern-2.5mu \cdot \mkern-2.5mu {#1}}
\newcommand{\divh}[1]{\nabla_h \mkern-2.5mu \cdot \mkern-2.5mu {#1}}

\newcommand{\tn}{|\!\!\:|\!\!\:|}
\newcommand\T{\rule{0pt}{2.6ex}}
\newcommand\B{\rule[-1.2ex]{0pt}{0pt}}
\definecolor{myred}{rgb}{0.9, 0.0, 0.0}
\definecolor{myblue}{rgb}{0.0, 0.28, 0.67}
\definecolor{mygreen}{rgb}{0.0, 0.7, 0.0}
\definecolor{myyellow}{rgb}{1.0, 0.55, 0.0}
\definecolor{mypurple}{rgb}{0.5, 0.1, 0.5}

\usepackage{soul}

\title{\textbf{Unified discontinuous Galerkin analysis of a thermo/poro-viscoelasticity model}}
\date{}

\author[1]{Stefano Bonetti}
\author[1]{Mattia Corti}
\affil[1]{\small{\textit{MOX-Dipartimento di Matematica, Politecnico di Milano, Piazza Leonardo da Vinci 32, Milan, 20133, Italy}}}

\begin{document}

\maketitle

\begin{center}
\begin{minipage}[c]{0.9\textwidth}
\input{Sections/Abstract}
\end{minipage}
\end{center}

\section{Introduction}
\label{sec:introduction}
\input{Sections/Introduction.tex}

\section{Model problem}
\label{sec:model_problem}
\input{Sections/ModelProblem.tex}

\section{Stability analysis}
\label{sec:stability_analysis}
\input{Sections/StabilityAnalysys.tex}

\section{Polytopal discontinuous Galerkin discretization}
\label{sec:discretization}
\input{Sections/Discretization}

\section{Numerical results}
\label{sec:numerical_results}
\input{Sections/NumericalResults}

\section{Conclusions}
\label{sec:conclusions}
\input{Sections/Conclusion}

\section*{Acknowledgements}
We acknowledge Paola F. Antonietti and Michele Botti (MOX - Dipartimento di Matematica, Politecnico di Milano, Italy) for the helpful discussions. MC has been funded from the European Union (ERC SyG, NEMESIS, project number 101115663). Views and opinions expressed are however those of the authors only and do not necessarily reflect those of the European Union or the European Research Council Executive Agency. Neither the European Union nor the granting authority can be held responsible for them. 
SB has been funded by the research grant PRIN2020 n. 20204LN5N5 funded by the Italian Ministry of Universities and Research (MUR). SB and MC are members of INdAM-GNCS. The present research is part of the activities of the project Dipartimento di Eccellenza 2023-2027, Dipartimento di Matematica, Politecnico di Milano.

\bibliography{bibliography.bib}

\appendix
\section{Derivation of Newmark-\texorpdfstring{$\beta$-$\theta$}{} method}
\label{sec:appendix}
\input{Sections/Appendix}

\section{Proof of convergence estimate (Theorem \ref{thm:conv})}
\label{sec:appendixconv}
\input{Sections/AppendixConvergence}

\end{document}

%% file: Sections/Abstract.tex
We present and analyze a discontinuous Galerkin method for the numerical modeling of a Kelvin-Voigt thermo/poro-viscoelastic problem. We present the derivation of the model and we develop a stability analysis in the continuous setting that holds both for the full inertial and quasi-static problems and that is robust with respect to most of the physical parameters of the problem. For spatial discretization, we propose an arbitrary-order weighted symmetric interior penalty scheme that supports general polytopal grids and is robust with respect to strong heterogeneities in the model coefficients. For the semi-discrete problem, we prove the extension of the stability result demonstrated in the continuous setting and we provide an a-priori error estimate. A wide set of numerical simulations is presented to assess the convergence and robustness properties of the proposed method. Moreover, we test the scheme with literature and physically sound test cases for \textit{proof-of-concept} applications in the geophysical context.

%% file: Sections/Introduction.tex
In recent years, there has been an increasing interest in studying the poroelasticity equations \cite{Kreuzer2024,Rodrigo2016,Rohan2013,VanDuijn2023,Both2017,Bociu2022}. The equations of linear poroelasticity are commonly known as Biot’s equations, and they find origin in the works of Biot \cite{Biot1941} and Terzaghi \cite{Terzaghi1943}. This model aims to study and describe the interaction between the fluid flow and elastic deformation within a porous medium. This problem was initially associated with geophysical applications, where the subsoil is modeled as a fully saturated poroelastic material (examples of application can be found, e.g., in \cite{Detournay1993,Russell1983,Castelletto2019}). Through the years, the classical poroelastic equations have been enhanced to couple them with other physical phenomena by including quantities that may influence (and may be influenced by) the fluid flow and the elastic deformation. Some examples can be found in poroelasto-acoustic coupling \cite{Antonietti2021}, which finds application in the context of earthquakes, and in the thermo-poroelasticity theory \cite{Brun2018,Brun2019}, which is used for modeling geothermal energy production procedures and greenhouse gas sequestration.
\par
A growing interest in the poroelasticity field has been motivated by its application to biological tissues. Indeed, organs, bones, and engineered tissue scaffolds can be modeled starting from the poroelasticity theory \cite{Bociu2021,Bociu2019,Lee2019}. Interesting results in the biological modeling, for example, in the brain context, can be obtained using the linear poroelasticity model \cite{Causemann2022}. More sophisticated multiple-network poroelastic theory \cite{Corti2022} considers multiple fluid compartments and, consequently, coupling terms coming from the interaction between different fluids.
\par
However, biological tissues typically exhibit both elastic and viscoelastic behavior due to the combined action of elastin and collagen \cite{Mow1980}. Linear viscoelastic effects can be incorporated into traditional linear Biot dynamics by considering the viscoelastic strain and possibly adjusting the formula for the local fluid content accordingly (depending on the specific scenario considered, either focusing on incompressible or compressible constituents). Some early works on poro-viscoelastic solids and their modeling are \cite{Coussy2003,Sanchez1980}. This work considers the simplest (linear) visco-elasticity inclusion: the Kelvin-Voigt type \cite{Bociu2016,Bociu2023}. We remark that the viscoelastic effects may be of interest also for geophysical applications. Indeed, traditionally, they were considered in the displacement dynamics by invoking the so-called secondary consolidation \cite{Murad1996,Gaspar2008}, for instance, for studies on clays.  
\par
This paper uses the full inertial Biot's system for describing poro-viscoelasticity, which is formally equivalent to thermoelasticity and thermo-viscoelasticity problems \cite{Showalter2000}. For instance, fully dynamic thermoelasticity can be used in the context of earthquake modeling. Indeed, a mechanical source of elastic waves induces a temperature field whose heat flow equalizes the temperature difference with the surroundings, giving rise to energy dissipation. On the other hand, a heat source generates elastic or viscoelastic waves \cite{Carcione2018}. Some \textit{proof-of-concept} numerical results in the direction of thermoelasticity are presented in this paper. Including viscous effects in the thermoelasticity theory may be useful for describing the behavior of the elastomers under large strains \cite{Boukamel2001}.
\par
This paper aims to provide a general and unified framework for studying thermo/poro-\hspace{-1pt}viscoelasticity problems in both their quasi-static and dynamic forms (i.e., considering inertial terms). Specifically, the treatment of the model problem highlights which terms are more or less significant depending on the values of the physical parameters and, consequently, on the reference application. Second, the stability analysis of the problem is developed to have only weak constraints on the model parameters; in particular, it is designed to be robust to the presence of the inertial terms within the two equations and the parameters describing the viscoelasticity characteristics in the model problem. We notice that the stability estimate obtained in this paper is not robust with respect to the conductivity coefficient, the second Lamé parameter and the first viscoelasticity retardation time. However, in the numerical results section, we observe that the proposed scheme is (at least) numerically robust with respect to low values of these model parameters.
\par
For the spatial discretization of the problem, we adopt a discontinuous Galerkin finite element method on polytopal grids (PolyDG) \cite{Antonietti2013,Cangiani2014}. The PolyDG schemes are appealing because of their geometrical flexibility, facilitating local mesh refinement and coarsening. Moreover, they allow us to efficiently handle highly heterogeneous media by better representing inner discontinuities. Geometrical flexibility is a desirable feature in both the geophysical and biological contexts, where discretization should capture peculiar features of the domain without dramatically affecting the number of elements.
Another advantage of polytopal discretization techniques is the possibility of exploiting arbitrary- and, in particular, high-order approximations. In second-order hyperbolic problems, high-order discretizations minimize the dispersion and diffusive phenomena. However, we remark that for subsurface applications, we may not have the regularity needed to justify high-order discretizations due to the large variabilities of the model's coefficients \cite{Flemisch2018,Stefansson2021}. The proposed method's arbitrary-order accuracy, combined with its geometric flexibility, allows using \textit{hp}-refinement techniques, taking full advantage of agglomeration techniques and reducing the overall computational cost.
\par
Examples of PolyDG schemes can be found in  \cite{Antonietti2021} for poroelasticity, in \cite{Antonietti2023,Bonetti2023,Bonetti2024} for thermo-hydro-mechanics (both in the quasi-static and fully-dynamic regimes), and in \cite{Corti2022} for multiple-network poroelasticity. Moreover, in \cite{Antonietti2021,Antonietti2022,DeLaPuente2008}, PolyDG methods for wave propagation problems in porous media are analyzed. Last, in \cite{Bajpai2023}, a discontinuous Galerkin method for the Kelvin-Voigt viscoelastic fluid model is proposed. In the literature, other discretization strategies for the poroelasticity problems include, e.g., finite Volume methods \cite{Berge2020,Schneider2017}, Hybrid Finite Volume method \cite{Asadi2021}, Hybrid High-Order \cite{BottiDiPietro2020}, Hybridizable Discontinuous Galerkin \cite{Fu2019}, lowest order Raviart-Thomas coupled with conforming finite element methods \cite{Both2017}, and eXtended finite element method \cite{Jafari2022}.
\par
The major highlights of this paper are: \textit{(a)} a detailed discussion of the model problem, in which the importance of the inertial terms as a function of the physical parameters of the problem is highlighted; \textit{(b)} an analysis covering both the static problem case and the fully inertial problem that investigates the stability of the problem under mild assumptions on the problem coefficients (for both the continuous and the semi-discrete problem) and an a-priori error estimate (for the semi-discrete problem); and \textit{(c)} a wide set of numerical simulations that are intended to prove the convergence and robustness properties of the method, test the method against literature benchmarks, and explore the applicability of the method in \textit{physically-sound} test cases. 
\par
The rest of the paper is organized as follows: the model problem, its derivation, and its weak formulation are reported in Section~\ref{sec:model_problem}. In Section~\ref{sec:stability_analysis}, we prove the stability of the continuous problem. In Section~\ref{sec:discretization}, we design the discretization of the problem. In particular, in Section~\ref{sec:DG_formulation}, we detail the PolyDG space discretization, and in Section~\ref{sec:time_discretization}, we show the time-discretization techniques we exploit for the hyperbolic-hyperbolic case and the hyperbolic-parabolic case. Moreover, we extend the stability result obtained in Section~\ref{sec:stability_analysis} for the discrete setting and we derive an a-priori error estimate. Last, we report numerical results in Section~\ref{sec:numerical_results}. Namely, we start by assessing the method's performance in terms of convergence properties and robustness, and then we address benchmark and physically sound test cases.

%% file: Sections/ModelProblem.tex
Let $\Omega\subset\mathbb{R}^d$, $d = \{2;3\}$, be an open, polygonal/polyhedral domain with Lipschitz boundary $\partial\Omega$.
Given a final time $T_f > 0$, the problem reads: \textit{find $(\mathbf{u}, \varphi)$ such that:}
\begin{subequations}
    \label{eq:system}
    \begin{empheq}[left=\empheqlbrace]{align}
    & \rho \ddot{\mathbf{u}} - \div{\boldsymbol{\sigma}}(\mathbf{u}, \dot{\mathbf{u}}, \varphi) = \mathbf{f} & \text{in }\Omega\times(0,T_f], \label{eq:momentum_cons_1} \\
    & d_0 \left(\dot{\varphi} + \tau_1 \ddot{\varphi} \right) + \gamma \left( \div{\dot{\mathbf{u}}} + \tau_2 \div{\ddot{\mathbf{u}}} \right) - \div{\left( \mathbf{D} \nabla \varphi \right)} = g
    & \text{in }\Omega\times(0,T_f], 
    \label{eq:mass_cons_1}   \end{empheq}
\end{subequations}
where the unknowns $(\mathbf{u},\varphi)$ stand for the solid displacement and a generalized pressure variable. In problem \eqref{eq:system}, equation \eqref{eq:momentum_cons_1} represents the momentum conservation law. At the same time, equation \eqref{eq:mass_cons_1} can represent both a mass conservation or an energy conservation law based on the physical interpretation given to the variable $\varphi$, which we discuss later in this section. Thus, the terms $\mathbf{f}$ and $g$ are source terms that represent a body force and a mass or energy source, respectively. Finally, problem \eqref{eq:system} is supplemented by imposing suitable boundary and initial conditions.
\par
In \eqref{eq:system} and in the rest of the article, we use the short-hand notation $\dot{\psi}$ and $\ddot{\psi}$ for denoting the first and second partial derivatives with respect to time of a function $\psi:\Omega\times(0,T_f]\to\mathbb{R}$, respectively.  
\par
The constitutive law for the total stress tensor $\boldsymbol{\sigma}$ is obtained as in \cite{Coussy2003} under the small deformations assumption, taking into account the isotropic effect of the generalized-pressure field on the porous matrix:
\begin{equation}
    \label{eq:const_law_sigma}
    \boldsymbol{\sigma}(\mathbf{u},\dot{\mathbf{u}}, \varphi) = 2 \mu \boldsymbol{\epsilon}(\mathbf{u}) + 2 \mu \delta_1 \boldsymbol{\epsilon}(\dot{\mathbf{u}}) +  \left( \lambda \nabla \cdot \mathbf{u} + \lambda \delta_2 \nabla \cdot \dot{\mathbf{u}} - \gamma \varphi \right) \mathbf{I}, 
\end{equation}
where $\mathbf{I}$ is the identity tensor and $\boldsymbol{\epsilon}(\mathbf{u})= \frac{1}{2}(\nabla \mathbf{u} + \nabla \mathbf{u}^T)$ is the strain tensor. Moreover, the dependency on the velocity of the deformation $\mathbf{u}$ models the secondary consolidation phenomena in poroelasticity for coil applications \cite{barden_primary_1968}. 
\par
The generalized pressure variable $\varphi$ can be considered a pressure in a poro-(visco)elasticity framework or a temperature in a thermo-(visco)elasticity one. In the first case, the second equation rises from the linearized equation of mass conservation \cite{Coussy2003}:
\begin{equation}
d_0\dot{\varphi}+\gamma(\nabla\cdot\dot{\mathbf{u}}) + \nabla\cdot\mathbf{w} = \tilde{g},
\end{equation}
where $\mathbf{w}$ is the filtration velocity, $d_0$ is the specific storage coefficient, and $\gamma$ is the Biot-Willis parameter \cite{Coussy2003}. Then, equation \eqref{eq:mass_cons_1} is obtained by adopting a Darcy law that considers the acceleration of the fluid \cite{Nield2017}:
\begin{equation}
\label{eq:dynamics_darcy}
    \dfrac{\rho_f}{\phi} \mathbf{K}\dot{\mathbf{w}}=- \mathbf{K}\nabla \varphi - \mathbf{w},
\end{equation}
where $\rho_f$ is the fluid density, $\phi$ is the porosity (the ratio between the void space in a porous medium and its whole volume such that $0 < \phi_0 \leq \phi \leq \phi_1 < 1$ \cite{Souzanchi2013}), $\mathbf{K}$ is the permeability divided by the dynamic fluid viscosity and $\tilde{g}$ is the forcing term. In particular, to obtain the formulation in equation   \eqref{eq:mass_cons_1}, we need to define $\tau_1 = \rho_f \mathbf{K} \phi^{-1}$,  $\tau_2 = \rho_f \mathbf{K} \phi^{-1}$, $\mathbf{D}=\mathbf{K}$, $g=\tau_1\dot{\tilde{g}}+\tilde{g}$. In the dynamical Darcy law in Equation \eqref{eq:dynamics_darcy}, we neglect the convection term ($(\mathbf{w}\cdot\nabla)\mathbf{w}$), because it is generally small and negligible in most applications \cite{Nield2017}.
\par
On the contrary, in the thermoelastic case, considering $\varphi$ as a temperature, equation \eqref{eq:mass_cons_1} can be derived from the linearized equation of energy conservation \cite{Coussy2003}:
\begin{equation}
d_0\dot{\varphi}+\gamma(\nabla\cdot\dot{\mathbf{u}}) + \nabla\cdot\mathbf{q} = \tilde{g},
\end{equation}
where $\mathbf{q}$ is the heat flux, $d_0$ is the thermal dilatation coefficient, and $\gamma$ is the thermal stress parameter \cite{Coussy2003}. Then, equation \eqref{eq:mass_cons_1} is obtained by adopting a Maxwell-Vernotte-Cattaneo law that considers a relaxation time in the heat conduction equation \cite{Nield2017}:
\begin{equation}
    \tau\dot{\mathbf{q}}= -\boldsymbol{\Theta} \nabla \varphi - \mathbf{q},
\end{equation}
where $\tau$ is the relaxation time, $\boldsymbol{\Theta}$ is the effective thermal conductivity and $\tilde{g}$ is the forcing term. In particular, to obtain the formulation in equation   \eqref{eq:mass_cons_1}, we need to define $\tau_1 = \tau$,  $\tau_2 = \tau$, $\mathbf{D}=\boldsymbol{\Theta}$, $g=\tau_1\dot{\tilde{g}}+\tilde{g}$. In both cases, the choice of considering two different values $\tau_1$ and $\tau_2$ in the final model is for generalization purposes. 
\par
The coefficients appearing in \eqref{eq:system}, along with their unit of measure and physical meaning in poroelasticity or thermoelasticity frameworks, are reported in Table~\ref{table:Biot_coeff}:
\begin{table}[H]
    \centering 
    \begin{tabular}{ c | c | l }
    \textbf{Parameter} & \textbf{Unit} & \textbf{Quantity} \\ \hline
    $\rho_f$ & \si[per-mode = symbol]
    {\kilogram \per \meter\cubed} & Saturating fluid density  \T\B\\ \hline
    $\rho_s$ & \si[per-mode = symbol]
    {\kilogram \per \meter\cubed} & Solid matrix density  \T\B\\ \hline
    $\phi$ & -  & \textbf{P:} Porosity \T\B\\ \hline
     \multirow{2}{*}{$\rho$} & \si{\kilogram \per \meter\cubed} & \textbf{P:} Density  $\rho = \phi \rho_f + (1 - \phi) \rho_s$ \T\B\\
    & \si{\kilogram \per \meter\cubed} & \textbf{T:} Density $\rho = \rho_s$ \T\B\\ \hline
    $\mu, \lambda$ & \si{\pascal} & Lamé parameters   \T\B\\ \hline
    $\delta_1, \delta_2$ & \si{\second} & Viscoelasticity retardation times \T\B\\ \hline
    \multirow{2}{*}{$\gamma$} & - & \textbf{P:} Biot--Willis constant \T\B\\
    & \si{\pascal \per\kelvin} & \textbf{T:} Thermal stress coefficient \T\B\\ \hline
    \multirow{2}{*}{$d_0$} & \si{\per\pascal} & \textbf{P:} Specific storage coefficient \T\B\\
    & \si{\pascal \per\kelvin\squared} & \textbf{T:} Thermal capacity \T\B\\ \hline
    \multirow{2}{*}{$\mathbf{D}$} & \si{\metre\squared \per \pascal \per \second} & \textbf{P:} Permeability divided by dynamic fluid viscosity \T\B\\
    & \si{\metre\squared \per \kelvin\squared \per \second} & \textbf{T:} Effective thermal conductivity \T\B\\ \hline    \multirow{2}{*}{$\tau_1$,$\tau_2$} & \si{\second} & \textbf{P:}  Relaxation times \T\B\\
    & \si{\second} & \textbf{T:} Maxwell-Vernotte-Cattaneo relaxation times \T\B\\ \hline
    \end{tabular}
    \\[10pt]
    \caption{Model coefficients appearing in problem~\eqref{eq:system} with explicit indication of the associated framework of the description: poroelasticity (P) or thermoelasticity (T). Where no indication is provided, the unit and the description are valid for both frameworks.}
    \label{table:Biot_coeff}
\end{table}
\par
As already mentioned, from the general formulation of equation \eqref{eq:system} we can recover some specific models.
\subsubsection*{Poroelasticity model}
The classical poroelastic Biot problem \cite{Biot1941}, can be recovered by taking the pressure field $p=\varphi$ and considering $\tau_1=\tau_2=\delta_1=\delta_2=0$.
\begin{equation}
    \label{eq:PE_system}
    \left\{
    \begin{aligned}
    & \rho \ddot{\mathbf{u}} - 2\div{(\mu \boldsymbol{\epsilon}(\mathbf{u}))} - \nabla {(\lambda \nabla \cdot \mathbf{u})} + \nabla{(\gamma p)} = \mathbf{f} \qquad && \text{in }\Omega\times(0,T_f], \\
    & d_0 \dot{p}+ \gamma \left( \div{\dot{\mathbf{u}}} \right) - \div{\left( \mathbf{D} \nabla p \right)} = g
    \qquad && \text{in }\Omega\times(0,T_f].
   \end{aligned}
    \right.
\end{equation}
\subsubsection*{Thermoelasticity model}
The thermoelasticity problem with the Maxwell-Vernotte-Cattaneo relaxation law can be recovered by taking the temperature field $T=\varphi$ and considering $\tau_1=\tau_2=\tau$ and $\delta_1=\delta_2=0$.
\begin{equation}
    \label{eq:TE_system}
    \left\{
    \begin{aligned}
    & \rho \ddot{\mathbf{u}} - 2\div{(\mu \boldsymbol{\epsilon}(\mathbf{u}))} - \nabla {(\lambda \nabla \cdot \mathbf{u})} + \nabla{(\gamma T)} = \mathbf{f} \qquad && \text{in }\Omega\times(0,T_f], \\
    & d_0 \left(\dot{T} + \tau \ddot{T} \right) + \gamma \left( \div{\dot{\mathbf{u}}} + \tau \div{\ddot{\mathbf{u}}} \right) - \div{\left( \mathbf{D} \nabla T \right)} = g
    \qquad && \text{in }\Omega\times(0,T_f].
    \end{aligned}
    \right.
\end{equation}
By taking also $\tau=0$, we recover the formulation with the Fourier law for the temperature.
\subsubsection*{Poro-viscoelasticity model}
The classical poro-viscoelasticity problem \cite{Bociu2023} can be recovered by taking the pressure field $p=\varphi$ and considering $\tau_1=\tau_2=0$.
\begin{equation}
    \label{eq:PVE_system}
    \left\{
    \begin{aligned}
    & \rho \ddot{\mathbf{u}} - \div{( 2 \mu (\boldsymbol{\epsilon}(\mathbf{u}) + \delta_1 \boldsymbol{\epsilon}(\dot{\mathbf{u}})))} -  \nabla{\left( \lambda (\nabla \cdot \mathbf{u} + \delta_2 \nabla \cdot \dot{\mathbf{u}}\right))} + \nabla{(\gamma p)} = \mathbf{f} \quad && \text{in }\Omega\times(0,T_f], \\
    & d_0 \dot{p} + \gamma \div{\dot{\mathbf{u}}} - \div{\left( \mathbf{D} \nabla p \right)} = g
    \quad && \text{in }\Omega\times(0,T_f].
    \end{aligned}
    \right.
\end{equation}
By taking $\delta_1=0$, we recover the secondary consolidation poroelastic model.
\subsection{Weak formulation}
\label{subsection:weak_formulation}
Before presenting the variational formulation of problem~\eqref{eq:system}, we introduce the required notation. 
For $X\subseteq\Omega$, we denote by $L^p(X)$ the standard Lebesgue space of index $p\in [1, \infty]$ and by $H^q(X)$ the Sobolev space of index $q \geq 0$ of real-valued functions defined on $X$. 
The notation $\mathbf{L}^p(X)$ and $\mathbf{H}^q(X)$ is adopted in place of $\left[ L^p(X) \right]^d$ and $\left[ H^q(X) \right]^d$, respectively. These spaces are equipped with natural inner products and norms denoted by $(\cdot, \cdot)_X = (\cdot, \cdot)_{L^2(X)}$ and $||\cdot||_X = ||\cdot||_{L^2(X)}$, with the convention that the subscript can be omitted in the case $X=\Omega$.
\par
We denote by $\langle \cdot,\cdot \rangle$ the duality pairing between the space $Y^*$ and $Y$; the former being the dual space of the latter. Moreover, we denote by $\| \cdot \|_*$ the dual norm in the space $Y^*$.
\par
For the sake of brevity, in what follows, we make use of the symbol $x \lesssim y$ to denote $x \le C y$, where $C$ is a positive constant independent of the discretization parameters but possibly dependent on physical coefficients and final time $T_f$.
\par
To derive the weak formulation of problem \eqref{eq:system}, we start by providing the definition of the functional spaces that take into account the essential boundary conditions, namely
\begin{equation}
\label{eq:func_spaces}
\begin{aligned}
V & = H_0^1(\Omega) = \left\{ \varphi \in H^1(\Omega) \ \text{s.t.} \  \varphi_{\vert \partial \Omega} = 0 \right\}, \quad \mathbf{V} = \left[ V \right]^d.
\end{aligned}
\end{equation}
Next, we multiply \eqref{eq:system} times suitable test functions, integrate in space, and we get: 
\par
\textit{For any time $t \in (0, T_f]$, find $(\mathbf{u}, \varphi)(t) \in \mathbf{V} \times V$ such that:}
\begin{equation}
    \label{eq:Biot_weak_form}
    \begin{aligned}
     \mathcal{M}_u(\ddot{\mathbf{u}}, \mathbf{v}) & + (2\mu\delta_1 \, \boldsymbol{\epsilon}(\dot{\mathbf{u}}),\boldsymbol{\epsilon}(\mathbf{v})) + (\lambda \delta_2 \, \div{\dot{\mathbf{u}}},\div{\mathbf{v}}) \\
        & + (2\mu \boldsymbol{\epsilon}(\mathbf{u}),\boldsymbol{\epsilon}(\mathbf{v})) + (\lambda \div{\mathbf{u}},\div{\mathbf{v}}) - (\gamma \div{\mathbf{v}}, \varphi) = (\mathbf{f}, \mathbf{v}) && \forall \mathbf{v}\in\mathbf{V}, 
    \\[3pt]
    \mathcal{M}_{\varphi, \tau_1}(\ddot{\varphi},\psi) & + (\gamma \tau_2 \, \div{\ddot{\mathbf{u}}}, \psi) + \mathcal{M}_{\varphi}(\dot{\varphi},\psi) + (\gamma \div{\dot{\mathbf{u}}}, \psi) + (\mathbf{D} \nabla \varphi, \nabla \psi) = (g, \psi) && \forall \psi \in V, 
    \end{aligned}
\end{equation}
where for every $(\mathbf{u}, \varphi), \, (\mathbf{v},\psi) \in \mathbf{V} \times V$ we have set: $\mathcal{M}_u(\mathbf{u}, \mathbf{v}) = (\rho \, \mathbf{u}, \mathbf{v})$, $\mathcal{M}_{\varphi, \tau_1}(\varphi, \psi) = (d_0 \tau_1 \, \varphi, \psi)$, and $\mathcal{M}_{\varphi}(\varphi, \psi) = (d_0 \, \varphi, \psi)$.

%% file: Sections/StabilityAnalysys.tex
This section is devoted to proving the stability bounds of the continuous weak solution of problem \eqref{eq:Biot_weak_form}. Before stating the theorem, we introduce the following auxiliary norms of the displacement $\mathbf{u}$ in $\mathbf{V}$:
\begin{align}
    \|\mathbf{u}\|_\mathbf{V}^2 = & \,2\|\sqrt{\mu}\epsilon(\mathbf{u})\|^2+\|\sqrt{\lambda}\nabla\cdot\mathbf{u}\|^2, \\
    \|\mathbf{u}\|_{\mathbf{V}\delta}^2 = & \,2\|\sqrt{\mu\delta_1}\epsilon(\mathbf{u})\|^2+\|\sqrt{\lambda\delta_2}\nabla\cdot\mathbf{u}\|^2.
\end{align}
Moreover, we introduce an assumption on the regularity of model parameters:
\begin{assumption}[Data regularity]
\label{ass:coeff_reg}
We assume the following regularities for the coefficients, the forcing terms, and the initial conditions:
\begin{enumerate}
    \item[$(i)$] $\mu,\, \lambda,\, \rho,\, d_0,\, \gamma,\, \tau_1,\, \tau_2,\, \delta_1,\, \delta_2 \in L^\infty(\Omega)$, moreover we require $\mu>0$;
    \item[$(ii)$] $\boldsymbol{\mathrm{D}}\in \mathbf{L}^\infty(\Omega)$ and $\exists d>0\;\forall\boldsymbol{\xi}\in \mathbb{R}^d:\; d|\boldsymbol{\xi}|^2 \leq \boldsymbol{\xi}^\top\mathbf{D}\boldsymbol{\xi} \quad \forall \boldsymbol{\xi}\in\mathbb{R}^d$;
    \item[$(iii)$] $\mathbf{f}\in C^0(0,T_f; \mathbf{H}^{-1}(\Omega))$;
    \item[$(iv)$] $g\in C^0(0,T_f; H^{-1}(\Omega))$;
    \item[$(v)$] $\mathbf{u}_0 \in \mathbf{H}^1_0(\Omega)$, $\dot{\mathbf{u}}_0 \in \mathbf{L}^2(\Omega)$, $\varphi_0 \in L^2(\Omega)$, and $\dot{\varphi}_0 \in L^2(\Omega)$.
\end{enumerate} 
\end{assumption}
\begin{theorem}
\label{thm:stab}
Let us consider -- for any time $t\in(0,T_f]$ -- $(\mathbf{u},\varphi)(t)\in\mathbf{V}\times V$ to be the solution of the weak problem \eqref{eq:Biot_weak_form} with homogeneous Dirichlet boundary conditions. Under Assumption \ref{ass:coeff_reg}, then the following stability estimate holds:
\begin{equation}
\begin{split}
& \|\sqrt{\tau_2\rho}\mathbf{u}\|^2 + \int_0^t \left(\|\tau_2\sqrt{\rho}\dot{\mathbf{u}}\|^2 + \|\sqrt{\rho}\mathbf{u}\|^2 + \left\|\dfrac{\tau_2}{2} \mathbf{u}\right\|^2_\mathbf{V} + \|\sqrt{\tau_2}\mathbf{u}\|^2_{\mathbf{V}\delta} + \|\sqrt{\tau_2 \tau_1 d_0} \varphi\|^2 \right) \lesssim \\ & \|\sqrt{\tau_2\rho}\mathbf{u}_0\|^2 + t I_0 + \int_0^t \left(\|\mathbf{F}\|^2_{*} 
+ \dfrac{\left\|\sqrt{\tau_2} G\right\|^2_{*}}{{\|\sqrt{\mathbf{D}}\|^2}} + t \left\|\dfrac{\tau_2\mathbf{f}}{\sqrt{\mu\delta_1}}\right\|^2_{*} + t\dfrac{\|\tau_2 g\|^2_{*}}{\|\sqrt{\mathbf{D}}\|^2} + t\left\|\dfrac{\mathbf{F}}{\sqrt{\mu\delta_1}}\right\|^2_{*} + t\dfrac{\|G\|^2_{*}}{\|\sqrt{\mathbf{D}}\|^2}\right),
\end{split}
\end{equation}
where:
\begin{equation}
\label{eq:int_forcing_cont}
\begin{aligned}
     \mathbf{F} = & \,\int_0^t \mathbf{f} + \rho \dot{\mathbf{u}}_0 - 2 \nabla\cdot\mu\delta_1\boldsymbol{\epsilon}(\mathbf{u}_0) - \nabla(\lambda\delta_2\nabla\cdot\mathbf{u}_0), \\
     G = & \, \int_0^t g + d_0 (\varphi_0 + \tau_1 \dot{\varphi}_0) + \gamma(\nabla\cdot\mathbf{u}_0+\tau_2\nabla\cdot\dot{\mathbf{u}}_0), \\
     I_0 = & \, \|\tau_2\sqrt{\rho}\dot{\mathbf{u}}_0\|^2 + \|\sqrt{\rho}\mathbf{u}_0\|^2 + \|\tau_2 \mathbf{u}_0\|^2_\mathbf{V} + \|\sqrt{\tau_2 \tau_1 d_0} \varphi_0\|^2.
\end{aligned}
\end{equation}
In this theorem, the (hidden) stability constant is independent of the physical parameters.
\end{theorem}
\begin{proof}
To start the proof of stability, we introduce two auxiliary variables $\Psi(t) = \int_0^t \varphi(s)\mathrm{d} s$ and $\mathbf{W}(t) = \int_0^t \mathbf{u}(s)\mathrm{d} s$. Then, we define two additional problems. The first is obtained by integrating Equation \eqref{eq:mass_cons_1} in time:
\begin{equation}
    \label{eq:system_1}
    \begin{dcases}
    \rho \ddot{\mathbf{u}} - 2 \div{\mu \boldsymbol{\epsilon}(\mathbf{u})} - 2 \div{\mu \delta_1 \boldsymbol{\epsilon}(\dot{\mathbf{u}})} -  \nabla(\lambda \nabla \cdot \mathbf{u}) - \nabla(\lambda \delta_2 \nabla \cdot \dot{\mathbf{u}}) + \nabla(\gamma \varphi) = \mathbf{f}, \\
    d_0 \left({\varphi} + \tau_1 \dot{\varphi} \right) + \gamma \left( \div{{\mathbf{u}}} + \tau_2 \div{\dot{\mathbf{u}}} \right) - \div{\left( \mathbf{D} \nabla \Psi \right)} = G, 
    \end{dcases}
\end{equation}
The second formulation is obtained by integrating both Equations \eqref{eq:momentum_cons_1} and \eqref{eq:mass_cons_1} in time:
\begin{equation}
    \label{eq:system_2}
    \begin{dcases}
   \rho \dot{\mathbf{u}} - 2 \div{\mu \boldsymbol{\epsilon}(\mathbf{W})} - 2 \div{\mu \delta_1 \boldsymbol{\epsilon}(\mathbf{u})} -  \nabla(\lambda \nabla \cdot \mathbf{W}) - \nabla(\lambda \delta_2 \nabla \cdot \mathbf{u}) + \nabla(\gamma \Psi) = \mathbf{F}, \\
    d_0 \left({\varphi} + \tau_1 \dot{\varphi} \right) + \gamma \left( \div{{\mathbf{u}}} + \tau_2 \div{\dot{\mathbf{u}}} \right) - \div{\left( \mathbf{D} \nabla \Psi \right)} = G. 
    \end{dcases}
\end{equation}
The two additional forcing terms $\mathbf{F}$ and $G$ are defined in Equation \eqref{eq:int_forcing_cont}. The first step of the proof consists in obtaining the weak formulations of \eqref{eq:system_1}, \eqref{eq:system_2}. To do so, we multiply the first equations of \eqref{eq:system_1}, \eqref{eq:system_2} times $\mathbf{v} \in \mathbf{V}$, and the second equations of \eqref{eq:system_1}, \eqref{eq:system_2} times $\psi \in V$. Then, we integrate by parts and sum the contributions of the two systems of equations, respectively. The total weak formulation of \eqref{eq:system_1} reads:
\begin{equation}
    \label{eq:system_1_wf}
    \begin{aligned}
    & (\rho \ddot{\mathbf{u}}, \mathbf{v}) + (d_0 \varphi, \psi) + (d_0 \tau_1 \dot{\varphi}, \psi) + (2\mu \boldsymbol{\epsilon}(\mathbf{u}),\boldsymbol{\epsilon}(\mathbf{v})) + (2\mu\delta_1 \boldsymbol{\epsilon}(\dot{\mathbf{u}}),\boldsymbol{\epsilon}(\mathbf{v})) + (\lambda \div{\mathbf{u}}, \div{\mathbf{v}}) \\
    & + (\lambda\delta_2 \div{\dot{\mathbf{u}}}, \div{\mathbf{v}}) + (\mathbf{D} \nabla \Psi, \nabla \psi) - (\varphi, \gamma\div{\mathbf{v}}) + (\gamma \div{{\mathbf{u}}}, \psi) + (\gamma\tau_2 \div{\dot{\mathbf{u}}}, \psi) = (\mathbf{f}, \mathbf{v}) + (G, \psi)
    \end{aligned}
\end{equation}
and the total weak formulation of \eqref{eq:system_2} reads:
total weak formulation of \eqref{eq:system_1} reads:
\begin{equation}
    \label{eq:system_2_wf}
    \begin{aligned}
    & (\rho \dot{\mathbf{u}}, \mathbf{v}) + (d_0 \varphi, \psi) + (d_0 \tau_1 \dot{\varphi}, \psi) + (2\mu \boldsymbol{\epsilon}(\mathbf{W}),\boldsymbol{\epsilon}(\mathbf{v})) + (2\mu\delta_1 \boldsymbol{\epsilon}(\mathbf{u}),\boldsymbol{\epsilon}(\mathbf{v})) + (\lambda \div{\mathbf{W}}, \div{\mathbf{v}}) \\
    & + (\lambda\delta_2 \div{\mathbf{u}}, \div{\mathbf{v}}) + (\mathbf{D} \nabla \Psi, \nabla \psi) - (\Psi, \gamma\div{\mathbf{v}}) + (\gamma \div{{\mathbf{u}}}, \psi) + (\gamma\tau_2 \div{\dot{\mathbf{u}}}, \psi) = (\mathbf{F}, \mathbf{v}) + (G, \psi).
    \end{aligned}
\end{equation}
Now, in \eqref{eq:system_1_wf} we take $(\mathbf{v}, \psi) = (\tau_2 \dot{\mathbf{u}}, \varphi)$ and multiply everything by $\tau_2$:
\begin{equation}
\begin{split}
 \dfrac{1}{2}\dfrac{\mathrm{d}}{\mathrm{d}t} \Big(\|\tau_2\sqrt{\rho}\dot{\mathbf{u}}\|^2 + & \|\tau_2 \mathbf{u}\|^2_\mathbf{V} + \|\sqrt{\tau_2 \tau_1 d_0} \varphi\|^2+ \|\sqrt{\tau_2\mathbf{D}} \nabla \Psi\|^2\Big) + \|\tau_2 \dot{\mathbf{u}}\|^2_{\mathbf{V}\delta} + \|\sqrt{\tau_2 d_0} \varphi\|^2 \\
= &  - (\gamma\tau_2 \nabla\cdot\mathbf{u},\varphi)  + (\tau_2^2\mathbf{f},\dot{\mathbf{u}}) + (\tau_2 G, \varphi)
\end{split}
\end{equation}
and in \eqref{eq:system_2_wf} we take $(\mathbf{v}, \psi) =(\mathbf{u},\Psi)$:
\begin{equation}
\begin{split}
\dfrac{1}{2}\dfrac{\mathrm{d}}{\mathrm{d}t} \Big(\|\sqrt{\rho}\mathbf{u}\|^2 + & \|\mathbf{W}\|^2_\mathbf{V} + \|\sqrt{d_0} \Psi\|^2 \Big) + \|\mathbf{u}\|^2_{\mathbf{V}\delta} + \|\sqrt{\mathbf{D}}\nabla \Psi\|^2 \\
= & (\gamma\tau_2 \div{\mathbf{u}},\varphi ) - \dfrac{\mathrm{d}}{\mathrm{d}t}(\gamma\tau_2 \nabla\cdot \mathbf{u},\Psi ) - (d_0 \tau_1 \dot{\varphi},\Psi) + (\mathbf{F}, \mathbf{u}) + (G,\Psi)
\end{split}
\end{equation}

Summing up the two equations and integrating in time, we obtain:
\begin{equation}
\begin{split}
 \|\tau_2&\sqrt{\rho}\dot{\mathbf{u}}\|^2 +  \|\sqrt{\rho}\mathbf{u}\|^2 + \|\tau_2 \mathbf{u}\|^2_\mathbf{V} + \|\mathbf{W}\|^2_\mathbf{V} + \|\sqrt{\tau_2 \tau_1 d_0} \varphi\|^2+ \|\sqrt{\tau_2\mathbf{D}} \nabla \Psi\|^2 +  \|\sqrt{d_0} \Psi\|^2 \\
+ & 2\int_0^t \|\mathbf{u}\|^2_{\mathbf{V}\delta} + 2\int_0^t \|\tau_2 \dot{\mathbf{u}} \|^2_{\mathbf{V}\delta} +2 \int_0^t \|\sqrt{\tau_2 d_0} \varphi\|^2 + 2\int_0^t\|\sqrt{\mathbf{D}}\nabla \Psi\|^2 + 2(\gamma\tau_2 \nabla\cdot \mathbf{u},\Psi ) \\ + & 2\int_0^t (d_0 \tau_1 \dot{\varphi},\Psi) = I_0 + 2\int_0^t(\tau_2^2\mathbf{f},\dot{\mathbf{u}}) + 2\int_0^t(\tau_2 G, \varphi) + 2\int_0^t (\mathbf{F}, \mathbf{u}) + 2\int_0^t  (G,\Psi),
\end{split}
\end{equation}
where $I_0 = \|\tau_2\sqrt{\rho}\dot{\mathbf{u}}_0\|^2 + \|\sqrt{\rho}\mathbf{u}_0\|^2 + \|\tau_2 \mathbf{u}_0\|^2_\mathbf{V} + \|\sqrt{\tau_2 \tau_1 d_0} \varphi_0\|^2$. First, we treat the forcing terms using integration by parts in time and Poincaré and Korn inequalities \cite{Duvaut1976}:
\begin{align}
2\int_0^t(\tau_2^2\mathbf{f},\dot{\mathbf{u}}) & \lesssim \int_0^t \left\|\dfrac{\tau_2\mathbf{f}}{\sqrt{\mu\delta_1}}\right\|^2_{*} + \int_0^t \|\tau_2\dot{\mathbf{u}}\|^2_{\mathbf{V}\delta}\\
2\int_0^t(\tau_2 G, \varphi) & = - 2\int_0^t(\tau_2 g, \Psi) + 2(\tau_2 G, \Psi) \\
&\lesssim \int_0^t \dfrac{\|\tau_2 g\|^2_{*}}{\|\sqrt{\mathbf{D}}\|^2} + \int_0^t \|\sqrt{\mathbf{D}}\nabla\Psi\|^2 + \dfrac{\left\|\sqrt{\tau_2} G\right\|^2_{*}}{{\|\sqrt{\mathbf{D}}\|^2}}+ \|\sqrt{\tau_2\mathbf{D}}\nabla\Psi\|^2 \\
 2\int_0^t (\mathbf{F},\mathbf{u}) & \lesssim \int_0^t \left\|\dfrac{\mathbf{F}}{\sqrt{\mu\delta_1}}\right\|^2_{*} + \int_0^t \|\mathbf{u}\|^2_{\mathbf{V}\delta}\\
 2\int_0^t  (G,\Psi) & \lesssim \int_0^t \dfrac{\|G\|^2_{*}}{\|\sqrt{\mathbf{D}}\|^2} + \int_0^t \|\sqrt{\mathbf{D}}\nabla\Psi\|^2.
\end{align}
Then, we obtain:
\begin{equation}
\label{eq:stab_preliminary2}
\begin{split}
 \|\tau_2&\sqrt{\rho}\dot{\mathbf{u}}\|^2 + \|\sqrt{\rho}\mathbf{u}\|^2 + \|\tau_2 \mathbf{u}\|^2_\mathbf{V} + \|\mathbf{W}\|^2_\mathbf{V} + \|\sqrt{\tau_2 \tau_1 d_0} \varphi\|^2 +  \|\sqrt{d_0} \Psi\|^2 + \int_0^t \|\mathbf{u}\|^2_{\mathbf{V}\delta} \\
+ & \int_0^t \|\tau_2 \dot{\mathbf{u}} \|^2_{\mathbf{V}\delta} 
 + 2 \int_0^t \|\sqrt{\tau_2 d_0} \varphi\|^2 + 2(\gamma\tau_2 \nabla\cdot \mathbf{u},\Psi ) + 2\int_0^t (d_0 \tau_1 \dot{\varphi},\Psi) \lesssim I_0 \\ + &  \dfrac{\left\|\sqrt{\tau_2} G\right\|^2_{*}}{{\|\sqrt{\mathbf{D}}\|^2}} + \int_0^t \left\|\dfrac{\tau_2\mathbf{f}}{\sqrt{\mu\delta_1}}\right\|^2_{*} + \int_0^t \dfrac{\|\tau_2 g\|^2_{*}}{\|\sqrt{\mathbf{D}}\|^2} + \int_0^t \left\|\dfrac{\mathbf{F}}{\sqrt{\mu\delta_1}}\right\|^2_{*} + \int_0^t \dfrac{\|G\|^2_{*}}{\|\sqrt{\mathbf{D}}\|^2}.
\end{split}
\end{equation}
We are left to control the tenth and eleventh terms on the left-hand side of \eqref{eq:stab_preliminary2}. Concerning the first of the two, under the assumption $\tau_2\geq\tau_1$:
\begin{align}
 2\int_0^t (d_0 \tau_1 \dot{\varphi},\Psi) = \dfrac{\mathrm{d}}{\mathrm{d}t}\|\sqrt{\tau_1 d_0}\Psi\|^2 - 2\int_0^t \| \sqrt{\tau_1 d_0} \varphi\|^2
\end{align}
For the second one, we test problem \eqref{eq:system_2} by $(\tau_2\mathbf{u},0)$ and we get:
\begin{equation}
\dfrac{1}{2}\dfrac{\mathrm{d}}{\mathrm{d}t} \left(\|\sqrt{\tau_2\rho}\mathbf{u}\|^2 + \|\sqrt{\tau_2}\mathbf{W}\|^2_\mathbf{V} \right) + \|\sqrt{\tau_2}\mathbf{u}\|^2_{\mathbf{V}\delta} = (\gamma\tau_2 \nabla\cdot \mathbf{u},\Psi ) + (\tau_2\mathbf{F}, \mathbf{u})
\end{equation}
Finally, by substituting the two quantities into \eqref{eq:stab_preliminary2}, we obtain that:
\begin{equation}
\begin{split}
\dfrac{\mathrm{d}}{\mathrm{d}t}(\|&\sqrt{\tau_2\rho}\mathbf{u}\|^2 + \|\sqrt{\tau_2}\mathbf{W}\|^2_\mathbf{V}+ \|\sqrt{\tau_1 d_0}\Psi\|^2) + \|\tau_2\sqrt{\rho}\dot{\mathbf{u}}\|^2 + \|\sqrt{\rho}\mathbf{u}\|^2 + \left\|\dfrac{\tau_2}{2} \mathbf{u}\right\|^2_\mathbf{V} + \|\sqrt{\tau_2}\mathbf{u}\|^2_{\mathbf{V}\delta}   \\
+ & \|\mathbf{W}\|^2_\mathbf{V} + \|\sqrt{\tau_2 \tau_1 d_0} \varphi\|^2 + \|\sqrt{d_0} \Psi\|^2 + \int_0^t \|\mathbf{u} \|^2_{\mathbf{V}\delta} + \int_0^t \|\tau_2 \dot{\mathbf{u}} \|^2_{\mathbf{V}\delta} + 2 \int_0^t \|\sqrt{(\tau_2-\tau_1) d_0} \varphi\|^2 \\ \lesssim & \|\mathbf{F}\|^2_{*} 
 + \dfrac{\left\|\sqrt{\tau_2} G\right\|^2_{*}}{{\|\sqrt{\mathbf{D}}\|^2}} + \int_0^t \left\|\dfrac{\tau_2\mathbf{f}}{\sqrt{\mu\delta_1}}\right\|^2_{*} + \int_0^t \dfrac{\|\tau_2 g\|^2_{*}}{\|\sqrt{\mathbf{D}}\|^2} + \int_0^t \left\|\dfrac{\mathbf{F}}{\sqrt{\mu\delta_1}}\right\|^2_{*} + \int_0^t \dfrac{\|G\|^2_{*}}{\|\sqrt{\mathbf{D}}\|^2}.
\end{split}
\end{equation}
After a final integration in time and neglecting the positive integrals on the introduced auxiliary variables, we obtain the thesis.
\end{proof}
\begin{remark}
Theorem \ref{thm:stab} provides two stability controls in $L^\infty((0,T_f],\mathbf{L}^2(\Omega))$ and $L^2((0,T_f],\mathbf{H}^1(\Omega))$ on the displacement solution $\mathbf{u}$. Concerning $\varphi$, we proved a control in $L^2((0,T_f],L^2(\Omega))$. Moreover, a final control for $\dot{\mathbf{u}}$ is provided in $L^2((0,T_f],\mathbf{L}^2(\Omega))$.
\end{remark}
\begin{remark}
Theorem \ref{thm:stab} provides a robust estimate with respect to the quasi-static case ($\rho=0$). Moreover, the elastic case $\delta_1=\delta_2=0$ is simple to be obtained by changing the control on the forcing term using Gr\"{o}nwall lemma \cite{Quarteroni2017}. We observe that the estimate is not robust with respect to the conductivity $\mathbf{D}$, the second Lamé parameter $\mu$, and the viscoelasticity retardation time $\delta_1$. However, in Section~\ref{sec:robtest} we show that the proposed scheme is numerically robust with respect to the cases $\mathbf{D}, \mu, \delta_1 \ll 1$.
\end{remark}

%% file: Sections/Discretization.tex
This section aims to derive the fully-discrete scheme of problem~\eqref{eq:system}. After introducing some preliminary concepts and assumptions about PolyDG methods, we detail the spatial discretization obtained via the PolyDG approximation, cf. Section~\ref{sec:DG_formulation}. Then, for the time-integration of this problem, we consider two different cases, depending on the value of the time-relaxation parameter $\tau_1$. When $\tau_1 > 0$, both the equations in \eqref{eq:system} are second-order hyperbolic, so we integrate in time with a Newmark-$\beta$ scheme. When $\tau_1 = 0$, equation \eqref{eq:momentum_cons_1} is second-order hyperbolic, while equation \eqref{eq:mass_cons_1} is parabolic, then we couple an implicit Newmark-$\beta$ scheme for \eqref{eq:momentum_cons_1}, with a $\theta$-method for \eqref{eq:mass_cons_1}, cf. Section~\ref{sec:time_discretization}.

\subsection{Preliminaries}

First, we present the mesh assumptions, the discrete spaces, and some instrumental results for designing and analyzing PolyDG schemes. 
We introduce a subdivision $\mathcal{T}_h$ of the computational domain $\Omega$, whose elements are polygons/polyhedra in dimension $d = 2, 3$, respectively. Next, we define the interfaces (or internal faces) as subsets of the intersection of any two neighboring elements of $\mathcal{T}_h$. If $d=2$, an interface is a line segment, while if $d=3$, an interface is a planar polygon that we assume can be further decomposed into a set of triangles. The same holds for the boundary faces collected in the set $\mathcal{F}_B$, which yields a simplicial subdivision of $\partial\Omega$.  Accordingly, we define $\mathcal{F}_I$ to be the set of internal faces and the set of all the faces as $\mathcal{F}_h=\mathcal{F}_B\cup\mathcal{F}_I$.
\par
As a basis for constructing the PolyDG approximation, we define fully-discontinuous polynomial spaces on the mesh $\mathcal{T}_h$. 
Given an element-wise constant polynomial degree $\ell:\mathcal{T}_h\to\mathbb{N}_{>0}
$, which determines the order of the approximation, the discrete spaces are defined such as
\begin{equation}
    \label{eq:discrete_spaces}
    \begin{aligned}
    V_h^{\ell} &= \left\{ v_h \in L^2(\Omega) : v_h |_{\kappa} \in \mathbb{P}^{\ell_{\kappa}}(\kappa) \ \ \forall \kappa \in \mathcal{T}_h \right\}, \quad \mathbf{V}_h^{\ell} = \left[ V_h^{\ell} \right]^d.\\
    \end{aligned}
\end{equation}
where, for each $\kappa\in\mathcal{T}_h$, the space $\mathbb{P}^{\ell_{\kappa}}(\kappa)$ is spanned by polynomials of maximum degree  $\ell_{\kappa}=\ell_{|\kappa}$. To analyze the convergence of the spatial discretization, we consider a mesh sequence $\{\mathcal{T}_h\}_{h\to0}$ satisfying the following properties:

\begin{assumption}
\label{ass:mesh_Th1}
The mesh sequence $\{\mathcal{T}_h\}_h$ satisfies the following properties \cite{dipietro2020}:
\begin{enumerate}[start=1,label={\bfseries A.\arabic*}]
    \item \label{ass:A1} Shape Regularity: $\forall K\in\mathcal{T}_h\;it\;holds: c_1 h_K^d\lesssim|K|\lesssim  c_2h_K^d$.
    \item \label{ass:A2} Contact Regularity: $\forall F\in\mathcal{F}$ with $F\subseteq \overline{K}$ for some $K\in\mathcal{T}_h$, it holds $h_K^{d-1}\lesssim |F|$, where $|F|$ is the Hausdorff measure of the face $F$.
    \item \label{ass:A3} Submesh Condition: There exists a shape-regular, conforming, matching simplicial submesh $\widetilde{\mathcal{T}_h}$ such that:
    \begin{itemize}
        \item[$(i)$] $\forall \widetilde{K}\in\widetilde{\mathcal{T}_h}\;\quad\exists K\in\mathcal{T}_h$ such that $\widetilde{K}\subseteq K$.
        \item[$(ii)$] The family $\{\widetilde{\mathcal{T}_h}\}_h$ is shape and contact regular.
        \item[$(iii)$] $\forall \widetilde{K}\in\widetilde{\mathcal{T}_h}, K\in\mathcal{T}_h$ with $\widetilde{K} \subseteq K$, it holds $h_K \lesssim h_{\widetilde{K}}$.
    \end{itemize}
\end{enumerate}
\end{assumption}
We remark that under \textit{\ref{ass:A1}} the following inequality (called \textit{discrete trace-inverse inequality}) holds (cf. \cite{dipietro2020} for all the details):
\begin{equation}
    \label{eq:trace_inverse_ineq}
    ||v||_{L^2(\partial \kappa)} \lesssim \frac{\ell_{\kappa}}{h_{\kappa}^{1/2}} ||v||_{L^2(\kappa)} \quad \forall v \in \mathbb{P}^{\ell_{\kappa}}(\kappa),
\end{equation}
where the hidden constant is independent of $\ell_{\kappa}, h_{\kappa}$, and the number of faces per element. In the following discussion, we introduce the Weighted Symmetric Interior Penalty method (WSIP) \cite{Ern2009}. The key ingredient of this method is to exploit weighted averages instead of the arithmetic ones used in the standard Interior Penalty formulation, cf. \cite{Wheeler1978,Arnold1982}. The use of weighted averages was first introduced for elliptic problems in \cite{Stenberg1998} and then developed for dG methods for dealing with advection-diffusion problems with locally vanishing diffusion \cite{Ern2009}. In \cite{Bonetti2024}, the PolyDG-WSIP discretization of a thermo-hydromechanics problem is presented.
\par
For the definition of the PolyDG-WSIP method, we introduce the weight function $\omega^+:\mathcal{F}_I\to [0,1]$ \cite{Heinrich2002}. Given an interior face $F \in \mathcal{F}_I$, we denote the values taken by $\omega^+$ and $\omega^- = 1-\omega^+$ on the face $F$ as $\omega\rvert_F^{+}$ and $\omega\rvert_F^{-}$, respectively. Given the function $\omega$, we can introduce the notion of weighted averages and jump operators, denoted with $\wavg{\cdot}$ and $\jump{\cdot}$, and normal jump, denoted with $\jump{\cdot}_n$ \cite{Arnold2002,Ern2009}:
\begin{equation}
	\label{eq:avg_jump_operators}
	\begin{aligned}
		& \jump{a} = a^+ \mathbf{n^+} + a^- \mathbf{n^-}, \quad && \jump{\mathbf{a}} = \mathbf{a}^+ \odot \mathbf{n^+} + \mathbf{a}^- \odot \mathbf{n^-}, \quad &&\jump{\mathbf{a}}_n = \mathbf{a}^+ \cdot \mathbf{n^+} + \mathbf{a}^- \cdot \mathbf{n^-}, \\ 
		& \wavg{a} = \omega^+ a^+ + \omega^- a^, \quad && \wavg{\mathbf{a}} = \omega^+ \mathbf{a}^+ + \omega^- \mathbf{a}^-, \quad && \wavg{\mathbf{A}} = \omega^+ \mathbf{A}^+ + \omega^- \mathbf{A}^-,
	\end{aligned}
\end{equation}
where $\mathbf{a} \odot \mathbf{n} = \mathbf{a}\mathbf{n}^T$ and $a, \ \mathbf{a}, \ \mathbf{A}$ are (regular enough) scalar-, vector-, and tensor-valued functions, respectively. The notation $(\cdot)^{\pm}$ is used for the trace on $F$ taken within the interior of $\kappa^\pm$ and $\mathbf{n}^\pm$ is the outer unit normal vector to $\partial \kappa^\pm$. When the subscript $\omega$ is omitted, we consider $\omega^+ = \omega^- = 1/2$. Accordingly, on boundary faces $F\in\mathcal{F}_B$, we set
$$
 \jump{a} = a \mathbf{n},\ \
\wavg{a} = a,\ \
\jump{\mathbf{a}} = \mathbf{a} \odot \mathbf{n},\ \
\wavg{\mathbf{a}} = \mathbf{a},\ \
\jump{\mathbf{a}}_n = \mathbf{a} \cdot \mathbf{n},\ \
\wavg{\mathbf{A}} = \mathbf{A},
$$ 
for the averages, this corresponds to consider $\omega^\pm$ single-valued and equal to $1$.

From now on, for the sake of simplicity, we assume that the model parameters are element-wise constant. Moreover, for later use, we can introduce the quantities
\begin{equation}
    \mu_{\kappa} = \mu \rvert_{\kappa}, \quad \lambda_{\kappa} = \lambda \rvert_{\kappa}, \quad \text{and } \overline{\mathbf{D}}_{\kappa} = |\sqrt{\mathbf{\mathbf{D}}\rvert_{\kappa}}|_2^2,
\end{equation}
where $|\cdot|_2$ denotes the $\ell^2$-norm in $\mathbb{R}^{d \times d}$. 

\subsection{Discontinuous Galerkin semi-discrete problem}
\label{sec:DG_formulation}

The aim of this section is to introduce the PolyDG-WSIP approximation of problem~\eqref{eq:system} and to derive the stability estimate of the semi-discrete formulation. The PolyDG-WSIP semi-discretization of problem \eqref{eq:Biot_weak_form} reads: \textit{for any $t \in (0, T_f]$, find $(\mathbf{u}_h, \varphi_h)(t) \in \mathbf{V}^{\ell}_h \times V^{\ell}_h $ such that $\forall \, (\mathbf{v}_h, \psi_h) \in \mathbf{V}^{\ell}_h \times V^{\ell}_h$:}
\begin{equation}
\label{eq:Biot_semidiscr}
\begin{aligned}
\mathcal{M}_{u}(\ddot{\mathbf{u}}_h, \mathbf{v}_h) & + \mathcal{M}_{\varphi, \tau_1}(\ddot{\varphi}_h, \psi_h) + \mathcal{C}_{\tau_2,h}(\ddot{\mathbf{u}}_h, \psi_h) + \mathcal{A}_{e,\delta_1,h}(\dot{\mathbf{u}}_h, \mathbf{v}_h) + \mathcal{A}_{\text{div},\delta_2,h}(\dot{\mathbf{u}}_h, \mathbf{v}_h) \\
& + \mathcal{M}_{\varphi}(\dot{\varphi}_h, \psi_h) + \mathcal{C}_h(\dot{\mathbf{u}}_h, \psi_h)
+ \mathcal{A}_{e,h}(\mathbf{u}_h, \mathbf{v}_h) + \mathcal{A}_{\text{div},h}(\mathbf{u}_h, \mathbf{v}_h) \\ & + \mathcal{A}_{\varphi,h}(\varphi_h, \psi_h) - \mathcal{C}_h(\mathbf{v}_h, \varphi_h) = (\mathbf{f},\mathbf{v}_h)+(g, \psi_h) 
\end{aligned}
\end{equation}
\textit{supplemented by initial conditions $(\mathbf{u}_{h,0}, \, \varphi_{h,0}, \, \dot{\mathbf{u}}_{h,0}, \, \dot{\varphi}_{h,0})$ that are fitting approximations of the initial conditions of problem~\eqref{eq:system}}. 
The bilinear forms labelled with the subscript ${}_h$  appearing in \eqref{eq:Biot_semidiscr} read

\begin{align}
    \mathcal{A}_{e,h}(\mathbf{u}, \mathbf{v}) = & (2 \mu\boldsymbol{\epsilon}_h(\mathbf{u}),\boldsymbol{\epsilon}_h(\mathbf{v}))
    + \sum_{F \in \mathcal{F}_h} \int_F \sigma \jump{\mathbf{u}} \mkern-2.5mu : \mkern-2.5mu \jump{\mathbf{v}} \\-& \sum_{F \in \mathcal{F}_h} \int_F \big( \avg{2 \mu\boldsymbol{\epsilon}_h(\mathbf{u})}_{\omega_{\mu}} \mkern-2.5mu : \mkern-2.5mu \jump{\mathbf{v}} + \jump{\mathbf{u}} \mkern-2.5mu : \mkern-2.5mu \avg{2 \mu\boldsymbol{\epsilon}_h(\mathbf{v})}_{\omega_{\mu}}\big),
\end{align}
\begin{align}
\mathcal{A}_{e,\delta_1,h}(\mathbf{u}, \mathbf{v}) = & (2 \mu\delta_1 \, \boldsymbol{\epsilon}_h(\mathbf{u}),\boldsymbol{\epsilon}_h(\mathbf{v}))
        + \sum_{F \in \mathcal{F}_h} \int_F \sigma_{\delta_1} \jump{\mathbf{u}} \mkern-2.5mu : \mkern-2.5mu \jump{\mathbf{v}} \\ -
        & \sum_{F \in \mathcal{F}_h} \int_F \big( \avg{2 \mu\delta_1 \, \boldsymbol{\epsilon}_h(\mathbf{u})}_{\omega_{\delta_1}} \mkern-2.5mu : \mkern-2.5mu \jump{\mathbf{v}} + \jump{\mathbf{u}} \mkern-2.5mu : \mkern-2.5mu \avg{2 \mu\delta_1 \, \boldsymbol{\epsilon}_h(\mathbf{v})}_{\omega_{\delta_1}} \big),
\end{align}
\begin{align}
\mathcal{A}_{\text{div},h}(\mathbf{u}, \mathbf{v}) = & (\lambda \divh{\mathbf{u}},  \divh{\mathbf{v}}) + \sum_{F \in \mathcal{F}_h} \int_F \big(\xi \jump{\mathbf{u}}_\mathbf{n} \jump{\mathbf{v}}_\mathbf{n} \big) \\
- & \sum_{F \in \mathcal{F}_h}\int_F \big(\avg{\lambda \divh{\mathbf{u}}}_{\omega_{\lambda}} \jump{\mathbf{v}}_\mathbf{n} + \jump{\mathbf{u}}_\mathbf{n} \avg{ \lambda\divh{\mathbf{v}}}_{\omega_{\lambda}} \big),
\end{align}
\begin{align}
        \mathcal{A}_{\text{div},\delta_2,h}(\mathbf{u}, \mathbf{v}) = \ & (\lambda\delta_2 \, \divh{\mathbf{u}},  \divh{\mathbf{v}}) +\sum_{F \in \mathcal{F}_h} \int_F \xi_{\delta_2} \jump{\mathbf{u}}_\mathbf{n} \jump{\mathbf{v}}_\mathbf{n} \\
        - & \sum_{F \in \mathcal{F}_h} \int_F \big(\avg{\lambda \delta_2 \, \divh{\mathbf{u}}}_{\omega_{\delta_2}} \jump{\mathbf{v}}_\mathbf{n} + \jump{\mathbf{u}}_\mathbf{n} \avg{ \lambda\delta_2 \, \divh{\mathbf{v}}}_{\omega_{\delta_2}} \big), 
\end{align}
\begin{align}
    \mathcal{A}_{\varphi,h}(\varphi,\psi) = & \left(\mathbf{D} \, \nabla_h \varphi, \nabla_h \psi\right) + \sum_{F \in \mathcal{F}_h} \int_F \zeta \jump{\varphi} \mkern-2.5mu \cdot \mkern-2.5mu \jump{\psi} \\ - & \sum_{F \in \mathcal{F}_h} \int_F \big(\avg{\mathbf{D} \, \nabla_h \varphi}_{\omega_{\mathbf{D}}} \mkern-2.5mu \cdot \mkern-2.5mu \jump{\psi} + \jump{\varphi} \mkern-2.5mu \cdot \mkern-2.5mu \avg{\mathbf{D} \, \nabla_h \psi}_{\omega_{\mathbf{D}}} - \zeta \jump{\varphi} \mkern-2.5mu \cdot \mkern-2.5mu \jump{\psi} \big),
\end{align}
\begin{align}
    & \mathcal{C}_h(\mathbf{u}, \psi) = \left( \gamma \divh{\mathbf{u}}, \psi \right) - \sum_{F \in \mathcal{F}_h} \int_F \avg{\psi} \jump{\gamma \mathbf{u}}_n,
\end{align}
\begin{align}
    & \mathcal{C}_{\tau_2,h}(\mathbf{u}, \psi) = \left( \gamma \tau_2 \, \divh{\mathbf{u}}, \psi \right) - \sum_{F \in \mathcal{F}_h} \int_F \avg{\psi} \jump{\gamma \tau_2 \, \mathbf{u}}_n.
\end{align}
Here, for all $a \in V_h^{\ell}$ and $\mathbf{a}\in \mathbf{V}_h^{\ell}$, $\nabla_h a$ and $\divh{\mathbf{a}}$ denote the broken differential operators whose restrictions to each element $\kappa \in \mathcal{T}_h$ are defined as $\nabla a_{|\kappa}$ and $\div{\mathbf{a}}_{|\kappa}$, respectively. Then, the broken version of the strain tensor is defined as $\boldsymbol{\epsilon}_h(\mathbf{u}) = \left(\nabla_h \mathbf{u} + \nabla_h \mathbf{u}^T\right)/2$. We set:
\begin{gather}
    \omega_{\mu}^{\pm} = \frac{\mu^{\mp}}{\mu^{+} + \mu^{-}}, \quad
    \omega_{\delta_1}^{\pm} = \frac{(\mu\delta_1)^{\mp}}{(\mu\delta_1)^{+} + (\mu\delta_1)^{-}}, \quad
    \omega_{\lambda}^{\pm} = \frac{\lambda^{\mp}}{\lambda^{+} + \lambda^{-}},
    \\
    \omega_{\delta_2}^{\pm} = \frac{(\lambda\delta_2)^{\mp}}{(\lambda\delta_2)^{+} + (\lambda\delta_2)^{-}},\quad
    \omega_{\mathbf{D}}^{\pm} = \frac{\eta_{\mathbf{D}_n}^{\mp}}{\eta_{\mathbf{D}_n}^{+} + \eta_{\mathbf{D}_n}^{-}}\nonumber
\end{gather}
where $\eta_{\mathbf{D}_n}^{\pm} = \mathbf{n}^{{\pm}^T} \, \mathbf{D}^{\pm} \, \mathbf{n}^{{\pm}}$. The PolyDG-WSIP method requires the definition of the following stabilization functions $\sigma, \sigma_{\delta_1}, \xi, \xi_{\delta_2}, \zeta \in L^{\infty}(\mathcal{F}_h)$ are defined according to \cite{Ern2009} as:
\begin{equation}
    \label{eq:stabilization_func}
    \begin{aligned}
    \sigma &= \left\{\begin{aligned}
    &\alpha_1 \gamma_{\mu} \underset{\kappa \in \{\kappa^+,\kappa^-\}}{\mbox{max}}\left( \frac{\ell_{\kappa}^2}{h_{\kappa}}\right) \ & F \in \mathcal{F}_I,\\
    &\alpha_1 \mu_{\kappa} \ell_\kappa^2 h_{\kappa}^{-1} \ & F \in \mathcal{F}_B,\\
    \end{aligned}
    \right.
    \ \
    \sigma_{\delta_1} = && \left\{\begin{aligned}
    &\alpha_2 \gamma_{\delta_1}  \underset{\kappa \in \{\kappa^+,\kappa^-\}}{\mbox{max}}\left( \frac{\ell_{\kappa}^2}{h_{\kappa}}\right) \ & F \in \mathcal{F}_I,\\
    &\alpha_2 \mu_{\kappa} \delta_{1_{\kappa}} \ell_\kappa^2 h_{\kappa}^{-1} \ & F \in \mathcal{F}_B,\\
    \end{aligned}
    \right.
    \\
    \xi & = \left\{\begin{aligned}
    &\alpha_3 \gamma_{\lambda} \underset{\kappa \in \{\kappa^+,\kappa^-\}}{\mbox{max}}\left( \frac{\ell_\kappa^2}{h_{\kappa}}\right) \ & F \in \mathcal{F}_I,\\
    &\alpha_3 \lambda_{\kappa} \ell_\kappa^2 h_{\kappa}^{-1} \ & F \in \mathcal{F}_B,\\
    \end{aligned}
    \right.
    \ \
    \xi_{\delta_2} = &&\left\{\begin{aligned}
    &\alpha_4 \gamma_{\delta_2} \underset{\kappa \in \{\kappa^+,\kappa^-\}}{\mbox{max}}\left( \frac{\ell_\kappa^2}{h_{\kappa}}\right) \ & F \in \mathcal{F}_I,\\
    &\alpha_4 \lambda_{\kappa} \delta_{2_{\kappa}} \ell_\kappa^2 h_{\kappa}^{-1} \ & F \in \mathcal{F}_B,\\
    \end{aligned}
    \right.\\
    \zeta &= \left\{\begin{aligned}
    &\alpha_5 \gamma_{\mathbf{D}} \underset{\kappa \in \{\kappa^+,\kappa^-\}}{\mbox{max}}\left( \frac{\ell_\kappa^2}{h_{\kappa}}\right) \ & F \in \mathcal{F}_I,\\
    &\alpha_5 {\overline{\mathbf{D}}_{\kappa}}^{-1} \ell_\kappa^2 h_{\kappa}^{-1} \ & F \in \mathcal{F}_B,\\
    \end{aligned}
    \right.
    \end{aligned}   
\end{equation}
where $\alpha_1, \alpha_2, \alpha_3, \alpha_4,\alpha_5 \in \mathbb{R}$ are positive constants to be properly defined, $\ell_{\kappa}$ is the (local) polynomial degree of approximation, $h_{\kappa}$ is the diameter of the element $\kappa \in \mathcal{T}_h$, and the coefficients $\gamma_{\mu}$, $\gamma_{\delta_1}$, $\gamma_{\lambda}$, $\gamma_{\delta_2}$, and  $\gamma_{\mathbf{D}}$ are given by:
\begin{gather}
    \gamma_{\mu}^{\pm} = \frac{\mu^{+} \, \mu^{-}}{\mu^{+} + \mu^{-}}, \quad
    \gamma_{\delta_1}^{\pm} = \frac{(\mu\delta_1)^{+} \, (\mu\delta_1)^{-}}{(\mu\delta_1)^{+} + (\mu\delta_1)^{-}}, \,
    \\
    \gamma_{\lambda}^{\pm} = \frac{\lambda^{+} \, \lambda^{-}}{\lambda^{+} + \lambda^{-}}, \quad
    \gamma_{\delta_2}^{\pm} = \frac{(\lambda\delta_2)^{+} \, (\lambda\delta_2)^{-}}{(\lambda\delta_2)^{+} + (\lambda\delta_2)^{-}}, \quad
    \gamma_{\mathbf{D}} = \frac{\eta_{\mathbf{D}_n}^{+} \, \eta_{\mathbf{D}_n}^{-}}{\eta_{\mathbf{D}_n}^{+} + \eta_{\mathbf{D}_n}^{-}}.
\end{gather}
\subsection{\textit{A-priori} analysis of semi-discrete formulation}
Before writing the stability theorem in the semi-discrete version, we introduce the following auxiliary dG norms $\forall \mathbf{u} \in \mathbf{V}_h^{\ell}$, $\forall \varphi \in V_h^{\ell}$:
\begin{align}
    \|\mathbf{u}_h\|_\mathrm{dG,e}^2 = & \,\|\sqrt{2\mu}\epsilon_h(\mathbf{u}_h)\|^2+\|\sqrt{\lambda}\divh{\mathbf{u}_h}\|^2 + \sum_{F\in\mathcal{F}} (\|\sqrt{\sigma}\jump{\mathbf{u}_h}\|^2_F + \|\sqrt{\xi}\jump{\mathbf{u}_h}_\mathbf{n}\|^2_F), \\
    \|\mathbf{u}\|_{\mathrm{dG},\delta}^2 = & \,\|\sqrt{2\mu\delta_1}\epsilon_h(\mathbf{u}_h)\|^2+\|\sqrt{\lambda\delta_2}\divh{\mathbf{u}_h}\|^2 + \sum_{F\in\mathcal{F}} (\|\sqrt{\sigma_{\delta_1}}\jump{\mathbf{u}_h}\|^2_F + \|\sqrt{\xi_{\delta_2}}\jump{\mathbf{u}_h}_\mathbf{n}\|^2_F), \\
    \|\varphi_h\|_{\mathrm{dG},\varphi}^2 = & \, \|\sqrt{\mathbf{D}}\nabla_h\varphi_h\|^2 + \sum_{F\in\mathcal{F}} \|\sqrt{\zeta}\jump{\varphi_h}\|^2_F.
\end{align}
Moreover, we introduce the additional dG norms to derive the convergence estimate $\forall \mathbf{u} \in \mathbf{H}^{\ell}(\mathcal{T}_h)$, $\forall \varphi \in H^{\ell}(\mathcal{T}_h)$:
\begin{align}
    \tn\mathbf{u}_h\tn_\mathrm{dG,e}^2 = &  \,\|\mathbf{u}_h\|_\mathrm{dG,e}^2 + \sum_{F\in\mathcal{F}} \|\wavgmu{\sqrt{2\mu}\mathbf{\epsilon}_h(\mathbf{u}_h)}\|^2_F, \\
    \tn\mathbf{u}_h\tn_\mathrm{dG,\delta}^2 = & \, \|\mathbf{u}_h\|_\mathrm{dG,\delta}^2 + \sum_{F\in\mathcal{F}} \|\wavgd{\sqrt{2\mu\delta_1}\mathbf{\epsilon}_h(\mathbf{u}_h)}\|^2_F, \\
    \tn\varphi_h\tn_{\mathrm{dG},\varphi}^2 = & \,\|\varphi_h\|_{\mathrm{dG},\varphi}^2 + \sum_{F\in\mathcal{F}} \|\wavgD{\mathbf{D}\nabla_h\varphi_h}\|^2_F.
\end{align}
\begin{theorem}[Stability estimate]
\label{thm:stabdisc}
Let us consider -- for any time $t\in(0,T_f]$ -- $(\mathbf{u}_h,\varphi_h)(t)\in\mathbf{V}_h^\ell\times V^\ell_h$ to be the solution of the semi-discrete problem \eqref{eq:Biot_semidiscr} with homogeneous Dirichlet boundary conditions. Under Assumptions \ref{ass:coeff_reg}, \ref{ass:mesh_Th1}, and assuming the following additional regularity on the forcing terms $\mathbf{f}\in C^0(0,T_f; \mathbf{L}^2(\Omega))$, $g\in C^0(0,T_f; L^2(\Omega))$, the following stability estimate holds:
\begin{equation}
\begin{split}
& \|\sqrt{\tau_2\rho}\mathbf{u}_h\|^2 + \int_0^t \left(\|\tau_2\sqrt{\rho}\dot{\mathbf{u}}_h\|^2 + \|\sqrt{\rho}\mathbf{u}_h\|^2 + \left\|\dfrac{\tau_2}{2} \mathbf{u}_h\right\|^2_\mathrm{dG,e} + \|\tau_2\mathbf{u}_h\|^2_{\mathrm{dG},\delta} + \|\sqrt{\tau_2 \tau_1 d_0} \varphi_h\|^2 \right) \lesssim \\ & \|\sqrt{\tau_2\rho}\mathbf{u}_{h,0}\|^2 + t I_{h,0} + \int_0^t \left(\|\mathbf{F}\|^2 
+ \dfrac{\left\|\sqrt{\tau_2} G\right\|^2}{{\|\sqrt{\mathbf{D}}\|^2}} + t\left\|\dfrac{\tau_2\mathbf{f}}{\sqrt{\mu\delta_1}}\right\|^2 \hspace{-4pt} + t\dfrac{\|\tau_2 g\|^2}{\|\sqrt{\mathbf{D}}\|^2} + t\left\|\dfrac{\mathbf{F}}{\sqrt{\mu\delta_1}}\right\|^2 \hspace{-4pt} + t\dfrac{\|G\|^2}{\|\sqrt{\mathbf{D}}\|^2}\right)\hspace{-2pt},
\end{split}
\end{equation}
where:
\begin{equation}
\label{eq:int_forcing}
\begin{aligned}
     \mathbf{F} = & \,\int_0^t \mathbf{f} + \rho \dot{\mathbf{u}}_{h,0} - 2 \nabla\cdot\mu\delta_1\boldsymbol{\epsilon}(\mathbf{u}_{h,0}) - \nabla(\lambda\delta_2\nabla\cdot\mathbf{u}_{h,0}), \\
     G = & \, \int_0^t g + d_0 (\varphi_{h,0} + \tau_1 \dot{\varphi}_{h,0}) + \gamma(\nabla\cdot\mathbf{u}_{h,0}+\tau_2\nabla\cdot\dot{\mathbf{u}}_{h,0}), \\
     I_{h,0} = & \, \|\tau_2\sqrt{\rho}\dot{\mathbf{u}}_{h,0}\|^2 + \|\sqrt{\rho}\mathbf{u}_{h,0}\|^2 + \|\tau_2 \mathbf{u}_{h,0}\|^2_\mathrm{dG} + \|\sqrt{\tau_2 \tau_1 d_0} \varphi_{h,0}\|^2.
\end{aligned}
\end{equation}
In this theorem, the (hidden) stability constant is independent of the physical parameters.
\end{theorem}
\begin{proof}
The proof of the stability can be adapted following the same steps of Theorem \ref{thm:stab}. The only differences are related to Poincarè and Korn's inequalities, for which we use their discrete version (cf. \cite{dipietro2020,BottiDiPietro2020}).
\end{proof}
To prove the \textit{a-priori} error estimate we introduce the interpolants of the solutions $\mathbf{u}_\mathrm{I}\in\mathbf{V}_h^\ell$ and $\varphi_\mathrm{I}\in V_h^\ell$ \cite{babuska-interpolant}.
\begin{proposition}
\label{prop:interpolant}
Let Assumption \ref{ass:mesh_Th1} be fulfilled. If $d \geq 2$, then the following estimates hold: 
\begin{equation}
    \forall \mathbf{u}\in\mathbf{H}^n(\mathcal{T}_h),\; \exists \mathbf{u}_\mathrm{I}\in \mathbf{V}_h^\ell: \; \|\mathbf{u}-\mathbf{u}_\mathrm{I}\|^2 + \tn \mathbf{u}-\mathbf{u}_\mathrm{I}\tn_\mathrm{dG,e}^2 \lesssim \sum_{K\in\mathcal{T}_h} \dfrac{h_K^{s-2}}{\ell^{2n-3}} \|\mathbf{u}\|_{\mathbf{H}^n(K)}^2,
\end{equation}
\begin{equation}
    \forall \varphi\in H^n(\mathcal{T}_h),\; \exists \varphi_\mathrm{I}\in V_h^\ell: \; \|\varphi-\varphi_\mathrm{I}\|^2 + \tn \varphi-\varphi_\mathrm{I}\tn_\mathrm{dG,\varphi}^2 \lesssim \sum_{K\in\mathcal{T}_h} \dfrac{h_K^{s-2}}{\ell^{2n-3}} \|\varphi\|_{H^n(K)}^2,
\end{equation}
where $s= \min\{\ell+1,n\}$.
\end{proposition}

\begin{theorem}[\textit{A-priori} error estimate]
\label{thm:conv}
Let us consider -- for any time $t\in(0,T_f]$ -- $(\mathbf{u},\varphi)(t)\in\mathbf{V}\times V$ and $(\mathbf{u}_h,\varphi_h)(t)\in\mathbf{V}_h^\ell\times V_h^l$ to be the solution of the problems \eqref{eq:Biot_weak_form} and \eqref{eq:Biot_semidiscr}, respectively, with homogeneous Dirichlet boundary conditions, and sufficiently large penalty parameters. Under Assumption \ref{ass:coeff_reg}, then the following estimate holds:
\begin{equation}
\label{eq:conv}
    \begin{aligned} \|\mathbf{u}-\mathbf{u}_h\|^2 + & \int_0^t \left( \|\varphi-\varphi_h\|^2 + \|\dot{\mathbf{u}}-\dot{\mathbf{u}}_h\|^2 + \|\mathbf{u}-\mathbf{u}_h\|^2 +  \|\mathbf{u}-\mathbf{u}_h\|_\mathrm{dG,e}^2 \right) 
    \lesssim \\
    (1 + t e^t) & \sum_{K\in\mathcal{T}_h} \dfrac{h_K^{s-2}}{\ell^{2n-3}} \Big [\|\mathbf{u}\|_{\mathbf{H}^n(K)}^2 + \int_0^t\Big(\|\ddot{\mathbf{u}}\|_{\mathbf{H}^n(K)}^2
    + \|\dot{\mathbf{u}}\|_{\mathbf{H}^n(K)}^2
    + \|\mathbf{u}\|_{\mathbf{H}^n(K)}^2 \\ 
    & \qquad + 
    \|\mathbf{W}\|_{\mathbf{H}^n(K)}^2 
    + \|\dot{\varphi}\|_{{H}^n(K)}^2
    + \|\varphi\|_{{H}^n(K)}^2
    + \|\Psi\|_{{H}^n(K)}^2 \Big)\Big],
    \end{aligned}
\end{equation}
where $\mathbf{W} = \int_0^t\mathbf{u}(s)\mathrm{d}s$, and $\Psi(t) = \int_0^t\varphi(s)\mathrm{d}s$.
\end{theorem}
\begin{proof}
    The proof of Theorem~\ref{thm:conv} is reported in Appendix \ref{sec:appendixconv}.
\end{proof}
\begin{remark}
In the statement of Theorem~\ref{thm:conv}, we neglect the dependence on the physical parameters. However, in the proof reported in Appendix \ref{sec:appendixconv}, the result has the same robustness properties as Theorem \ref{thm:stab}, due to the use of the same inequalities in the proof of the stability estimate.
\end{remark}

\subsection{Time discretization}
\label{sec:time_discretization}

This section aims to introduce the time discretization of the semi-discrete problem \eqref{eq:Biot_semidiscr}. The time-integration scheme depends on the parameters we consider in \eqref{eq:system}. Namely, when $\tau_1 > 0$, the system is second-order hyperbolic, and we choose to use a Newmark-$\beta$ method for the whole system. When $\tau_1 = 0$, the second equation \eqref{eq:mass_cons_1} is parabolic, then we couple a Newmark-$\beta$ method for the first equation with a $\theta$-method for the second one. We denote by Newmark-$\beta$-$\theta$ the coupling of these two time-marching schemes. 
\par
By fixing a basis for the space $\mathbf{V}_h^{\ell} \times V_h^{\ell}$ and denoting by $\left[ \mathbf{U}, \boldsymbol{\Phi} \right]^T$ the vector of the expansion coefficients of the variables $(\mathbf{u}_h, \varphi_h)$, we can rewrite the semi-discrete problem \eqref{eq:Biot_semidiscr} in the equivalent form:
\begin{equation}
\label{eq:semidiscr_alg_sys}
\begin{aligned}
\begin{bmatrix}
    \mathbf{M}_{u} & 0 \\
    \mathbf{C}_{\tau_2} & \mathbf{M}_{\varphi, \tau_1}
\end{bmatrix}
\begin{bmatrix}
    \ddot{\mathbf{U}} \\
    \ddot{\boldsymbol{\Phi}} 
\end{bmatrix} &+
\begin{bmatrix}
    \mathbf{A}_{e,\delta_1} + \mathbf{A}_{\text{div},\delta_2} & 0 \\
    \mathbf{C} & \mathbf{M}_{\varphi}
\end{bmatrix}
\begin{bmatrix}
    \dot{\mathbf{U}} \\
    \dot{\boldsymbol{\Phi}} 
\end{bmatrix} +
\begin{bmatrix}
    \mathbf{A}_e + \mathbf{A}_{\text{div}} & - \mathbf{C}^T \\
    0 & \mathbf{A}_{\varphi}
\end{bmatrix}
\begin{bmatrix}
    \mathbf{U} \\
    \boldsymbol{\Phi} 
\end{bmatrix} =
\begin{bmatrix}
    \mathbf{F} \\
    \mathbf{G}
\end{bmatrix}
\end{aligned}
\end{equation}
with initial conditions $\mathbf{U}(0) = \mathbf{U}_0$, $\dot{\mathbf{U}}(0) = \mathbf{U}_1$, $\Phi(0) = \Phi_0$, and (when $\tau_1>0$) $\dot{\Phi}(0) = \Phi_1$. The vectors $\mathbf{F}, \mathbf{G}$ are representations of the linear functionals appearing on the right-hand side of \eqref{eq:Biot_semidiscr}.
\par
To integrate \eqref{eq:semidiscr_alg_sys} in time, we introduce a time-step $\Delta t = T_f/n$, with $n\in\mathbb{N}_{>0}$, discretize the interval $(0, T_f]$ as a sequence of time instants $\{ t_k \}_{0\le k\le n}$ such that $t_{k+1} - t_k = \Delta t$.

\subsubsection{Case \texorpdfstring{$\tau_1 > 0$}{} (Newmark-\texorpdfstring{$\beta$}{} method)}
We start by defining $\mathbf{X}^k = \mathbf{X}(t^k)$, with $\mathbf{X} = \left[\mathbf{U}, \boldsymbol{\Phi} \right]^T$. Next, we rewrite \eqref{eq:semidiscr_alg_sys} in a compact form as $\mathbf{A} \ddot{\mathbf{X}} + \mathbf{B} \dot{\mathbf{X}} + \mathbf{C} \mathbf{X} = \mathbf{F}$ and derive
\begin{equation}
    \label{eq:def_L_timeint}
    \ddot{\mathbf{X}} = \mathbf{A}^{-1} \left( \mathbf{F} - \mathbf{B} \dot{\mathbf{X}} - \mathbf{C} \mathbf{X} \right) = \mathbf{A}^{-1} \mathbf{F} - \mathbf{A}^{-1} \mathbf{B} \dot{\mathbf{X}} - \mathbf{A}^{-1}\mathbf{C} \mathbf{X} = \mathcal{L}(t, \mathbf{X}, \dot{\mathbf{X}}).
\end{equation}
Last, we integrate in time \eqref{eq:def_L_timeint} with the use of Newmark-$\beta$ scheme, that exploits a Taylor expansion for $\mathbf{X}$ and $\mathbf{Y} = \dot{\mathbf{X}}$:
\begin{equation}
    \label{eq:newmark}
    \left\{
    \begin{aligned}
        & \mathbf{X}^{k+1} = \mathbf{X}^{k} + \Delta t \mathbf{Y}^{k} + \Delta t^2 \left( \beta_N \mathcal{L}^{k+1} + (\frac{1}{2} - \beta_N) \mathcal{L}^{k} \right), \\
        & \mathbf{Y}^{k+1} = \mathbf{Y}^{k} + \Delta t \left( \gamma_N \mathcal{L}^{k+1} + (1 - \gamma_N) \mathcal{L}^{k} \right),
    \end{aligned}
    \right.
\end{equation}
where $\mathcal{L}^k = \mathcal{L}(t^k, \mathbf{X}^k, \dot{\mathbf{X}}^k)$ and the Newmark parameters $\beta_N, \gamma_N$ satisfy: $0 \leq 2 \beta_N \leq 1$, $0 \leq \gamma_N \leq 1$. The typical choices for the Newmark parameters, that ensure unconditional stability and second-order accuracy for the scheme, are $\beta_N = 1/4$ and $\gamma_N = 1/2$. These are the values used in all the numerical tests of Section~\ref{sec:numerical_results}.

\subsubsection{Case \texorpdfstring{$\tau_1 = 0$}{} (Coupling Newmark-\texorpdfstring{$\beta$}{} and \texorpdfstring{$\theta$}{}-methods)}
In this second case, we adopt a coupling between a Newmark-$\beta$ method for discretizing the first equation of \eqref{eq:semidiscr_alg_sys} and a $\theta$-method for the pressure equation. The complete calculations are reported in the appendix section \ref{sec:appendix}. For the sake of clarity, we report the final formulation in a compact form:
\begin{equation}
    \label{eq:full_time_discretization_matrix}
    \mathbf{L} \mathbf{X}^{n+1} = \mathbf{R} \mathbf{X}^{n} + \mathbf{S}^{n+1} \qquad n>0, 
\end{equation}
where:
\begin{equation}
    \label{eq:full_time_discretization_matrix2}
    \begin{aligned}
    & \mathbf{L} = \begin{bmatrix}
    \frac{1}{\beta_N \Delta t^2} \mathbf{M}_u + \mathbf{A}_u + \frac{\gamma_N}{\beta_N \Delta t} \mathbf{A}_{u,\delta} & -\mathbf{C}^T & 0 & 0 \\
    \mathbf{C}_u & \frac{1}{\Delta t}\mathbf{M}_{\varphi} + \theta \mathbf{A}_{\varphi} & 0 & 0 \\
    0 & 0 & \mathbf{I} & -\Delta t \gamma_N \mathbf{I} \\
    - \frac{1}{\beta_N \Delta t^2} \mathbf{I} & 0 & 0 & \mathbf{I} \\
    \end{bmatrix} \\
    & \mathbf{R} = \begin{bmatrix}
    \frac{1}{\beta_N \Delta t^2} \mathbf{M}_u + \frac{\gamma_N}{\beta_N \Delta t} \mathbf{A}_{u,\delta} & 0 & \frac{1}{\beta_N \Delta t} \mathbf{M}_u - \frac{\beta_N-\gamma_N}{\beta_N} \mathbf{A}_{u,\delta} & \frac{1 - 2\beta_N}{2\beta_N} \mathbf{M}_u - \frac{\Delta t (2\beta_N-\gamma_N)}{2\beta_N} \mathbf{A}_{u,\delta}  \\
    \mathbf{C}_u & \frac{1}{\Delta t}\mathbf{M}_{\varphi} - \Tilde{\theta}\mathbf{A}_{\varphi} & \mathbf{C}_z & \mathbf{C}_a \\
    0 & 0 & \mathbf{I} & \Delta t (1-\gamma_N) \mathbf{I} \\
    - \frac{1}{\beta_N \Delta t^2} \mathbf{I} & 0 & - \frac{1}{\beta_N \Delta t} \mathbf{I} & \frac{2\beta_N - 1}{2 \beta_N}\mathbf{I} \\
    \end{bmatrix} \\[10pt]
    & \mathbf{X}^n = \begin{bmatrix}
    \mathbf{U}^n, \ \boldsymbol{\Phi}^n, \ \mathbf{Z}^n, \ \mathbf{A}^n
    \end{bmatrix}^T 
    \qquad
    \mathbf{S}^{n+1} = \begin{bmatrix}
    \mathbf{F}^{n+1}, \ \theta \mathbf{G}^{n+1} + \Tilde{\theta} \mathbf{G}^{n}, \ 0, \ 0
    \end{bmatrix}^T
    \end{aligned}
\end{equation}

\begin{remark}
    From an operational point of view, we do not solve the whole system \eqref{eq:full_time_discretization_matrix}. For the sake of reducing the overall computational cost, we solve it just for $\mathbf{U}^{n+1}$, $\boldsymbol{\Phi}^{n+1}$. Then, we update the values of $\mathbf{A}^{n+1}$ and $\mathbf{Z}^{n+1}$.
\end{remark}

%% file: Sections/NumericalResults.tex
This section aims to assess the performance of the proposed scheme in terms of accuracy and robustness with respect to the model parameters. Then, we test the method by addressing some benchmark and literature test cases. The numerical implementation is carried out in the open-source \texttt{lymph} library \cite{Antonietti2024}, implementing the PolyDG method for multiphysics. In all the tests, the PolyDG space discretization is solved monolithically, and it is coupled with the Newmark-$\beta$ method when $\tau_1 > 0$ and with a coupled Newmark-$\beta$-$\theta$-method when $\tau_1 = 0$ (cf. Section~\ref{sec:time_discretization}) for time integration. The parameters of the Newmark-$\beta$ and $\theta$-method are $\gamma_N = 1/2$, $\beta_N = 1/4$, and $\theta = 1/2$. All the penalty coefficients $\alpha_i$, $i=1,2,\dots,5$ in \eqref{eq:stabilization_func} are set equal to $10$.

\subsection{Convergence tests}
\label{sec:convtest}
\input{Sections/ConvTest}

\subsection{Superconvergence tests}
\label{sec:superconvtest}
\input{Sections/SuperConvTest}

\subsection{Robustness tests}
\label{sec:robtest}
\input{Sections/RobTest}

\subsection{Wave propagation in thermo-elastic media}
\label{sec:carcione2018}
\input{Sections/Carcione2018}

\subsection{Fluid flow in heterogeneous poro-viscoelastic media}
\label{sec:spe10}
\input{Sections/Spe10}

%% file: Sections/ConvTest.tex
\textit{Convergence vs space discretization parameters.} 
The aim of this section is to assess the performance of the proposed scheme in terms of accuracy with respect to the space discretization parameters, i.e., the mesh size $h$ and the polynomial degree of approximation $\ell$.
\begin{table}[ht]
\begin{minipage}[b]{0.4\linewidth}
    \centering
    \includegraphics[width=0.50\linewidth]{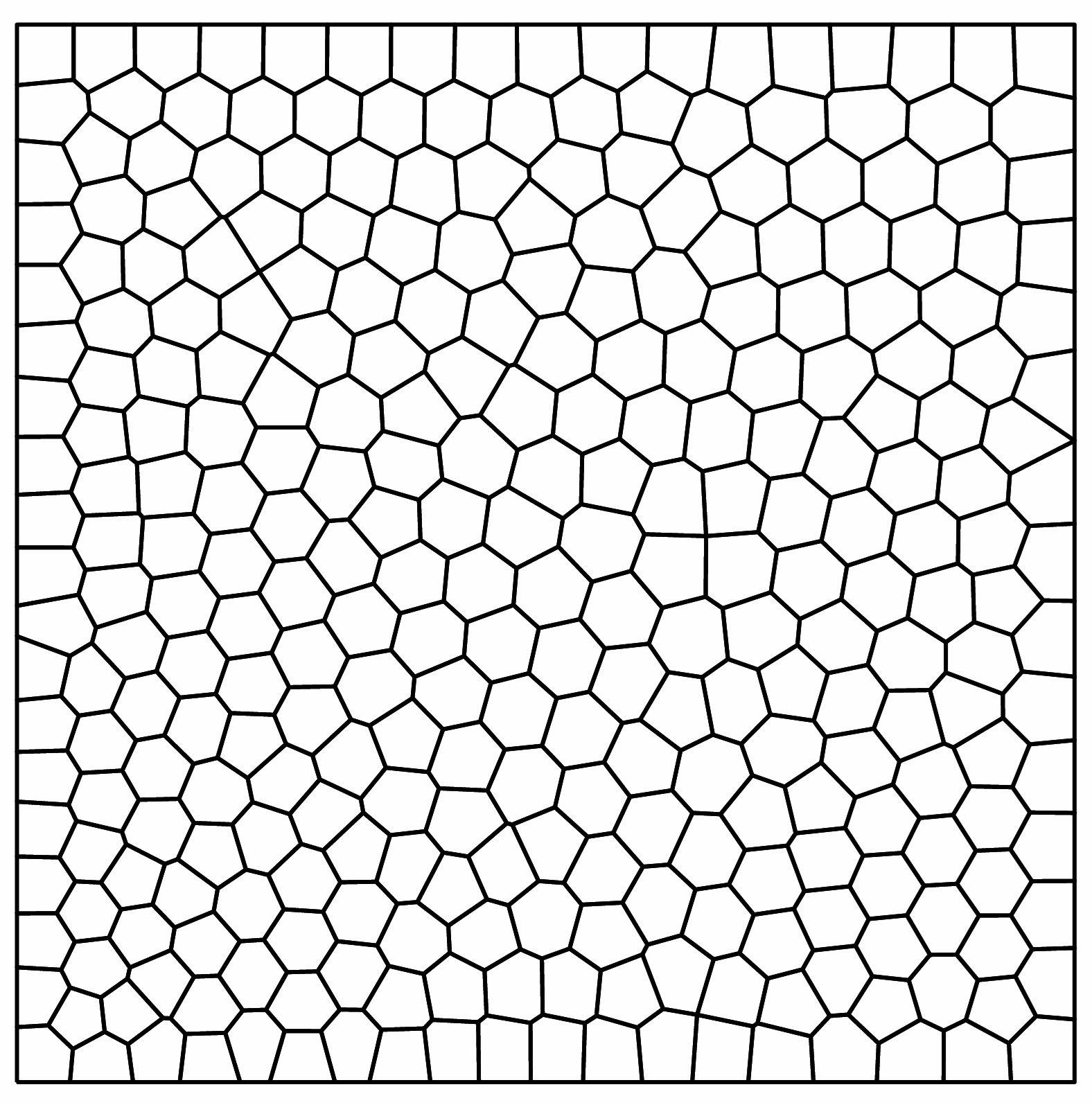}
    \vspace{-0.2cm}
    \captionof{figure}{Convergence test: example of a $2D$ Voronoi polygonal mesh made of 300 elements.}
    \label{fig:voronoi_mesh}
\end{minipage}\hfill
\begin{minipage}[b]{0.56\linewidth}
    \centering 
    \small
    \begin{tabular}{ l | l  c l | l}
    \textbf{Coefficient} & \textbf{Value} & & \textbf{Coefficient} & \textbf{Value} \T\B \\
    $\rho \ [\si{\kilogram \per \meter\cubed}]$ & 1 & & $\gamma \ [-]$ & 1 \T\B \\ 
    $\mu \ [\si{\pascal}]$ & 1 & & $\lambda \ [\si{\pascal}]$ & 1 \T\B \\ 
    $d_0 \ [\si{\per \pascal}]$ & 1 & & $\mathbf{D} \ [\si{\meter \squared \per \pascal \per \second}]$ & $\mathbf{I}$  \T\B \\
    $\delta_1 [\si{\second}]$ & 1 & & $\delta_2 [\si{\second}]$ & 1 \T\B \\
    $\tau_1 [\si{\second}]$ & 1 & & $\tau_2 [\si{\second}]$ & 1 \T\B \\
    \end{tabular}
\caption{Convergence test: problem parameters for the convergence analysis}
\label{tab:params_convtest}
\end{minipage}
\end{table}
We consider problem \eqref{eq:system} in the square domain $\Omega = (0,1)^2$ with manufactured analytical solutions:
\begin{equation}
    \begin{aligned}
    \mathbf{u} =  \sin(2 \pi t) \left[\begin{aligned}
    x^2 \sin(\pi x) \sin(\pi y) \\
    - x^2 \sin(\pi x) \sin(\pi y)
    \end{aligned} \right], \quad 
    \varphi = \sin(\sqrt{2} \pi t) \left( x^2 \cos\left(\frac{\pi x}{2}\right) \sin(\pi x) \right) .
    \end{aligned}
\end{equation}
The initial and boundary conditions and the forcing terms are inferred from the exact solutions. The model coefficients are chosen as reported in Table~\ref{tab:params_convtest}. We remark that - for completeness - we also assess the method's performance for the case $\tau_1 = 0$. 

For the $h$-convergence, a sequence of polygonal meshes in Figure~\ref{fig:voronoi_mesh} is considered, and we consider polynomial degree $\ell = 3$. At the same time, for the $\ell$-convergence we fix a computational mesh of $100$ elements and vary the polynomial degree $\ell = 1,2,\dots,5$. The time discretization parameters are $T_f = 0.1$, $\Delta t = \num[exponent-product=\ensuremath{\cdot}, print-unity-mantissa=false]{5e-5}$. 

\input{IMG/ConvH.tex}
\input{IMG/ConvH_tau0.tex}

In Figure~\ref{fig:convh}, we report the $L^2$ and $dG$-errors for the two unknowns with respect to the mesh size (\textit{log-log} scale). In agreement with the results expected by Theorem~\ref{thm:conv} (a-priori error estimate) and by the theory of the PolyDG methods (cf. \cite{CangianiDong2017,Cangiani2014}), we observe that, as we are using $\ell = 3$, the $dG$-errors show a convergence rate proportional to $h^3$. Moreover, concerning the $L^2$-errors, we observe that the errors decay as $h^{\ell + 1}$. In Figure~\ref{fig:convh_tau0}, we report the results for the same convergence test, but taking $\tau_1 = 0$ in \eqref{eq:mass_cons_1}. Also, for this configuration, the numerical results match the predicted convergence rates of the PolyDG framework.

\input{IMG/ConvP.tex}
\input{IMG/ConvP_tau0.tex}

In Figure~\ref{fig:convp}, Figure~\ref{fig:convp_tau0}, we report the results for the $\ell$-convergence test. We observe that, for both the model configurations, the errors of the two variables decay exponentially.

\vspace{0.5cm}
\textit{Convergence vs time discretization parameter.} This section aims to assess the performance of the proposed scheme in terms of accuracy with respect to the time-step $\Delta t$. To focus on the time-marching scheme, we consider problem \eqref{eq:system} in the square domain $\Omega = (0,1)^2$ with manufactured analytical solutions 
\begin{equation}
    \begin{aligned}
    \mathbf{u} = \sin(2 \pi t) \left[ \begin{aligned}
    & x + y \\
    & 3x - 5y
    \end{aligned} \right], \quad
    \varphi = \sin(\sqrt{2} \pi t) \left( 10x + 6y \right) .
    \end{aligned}
\end{equation}
The initial conditions, boundary conditions, and forcing terms are inferred from the exact solutions. The model coefficients are chosen as reported in Table~\ref{tab:params_convtest}. To test the convergence of the scheme with respect to $\Delta t$ (both for the Newmark-$\beta$ and for the coupled Newmark-$\beta$-$\theta$ method), we fix a mesh of $N = 100$ elements with $\ell = 1$, and we fix the final time $T_f = 0.1$. Then, we vary the time-step parameter.

\input{IMG/ConvDT.tex}
\input{IMG/ConvDT_tau0.tex}

In Figure~\ref{fig:convdt}, Figure~\ref{fig:convdt_tau0}, we report the results for the convergence against the time-step $\Delta t$. We observe that for both the Newmark-$\beta$ scheme ($\tau_1 \neq 0$) and for the coupled Newmark-$\beta$-$\theta$ time-marching schemes ($\tau_1 = 0$), we have second-order accuracy. The time-integration schemes' parameters we used in this test are $\gamma_N = 1/2$, $\beta_N = 1/4$, and $\theta = 1/2$, as these choices ensure unconditional stability and second-order accuracy.

%% file: IMG/ConvH.tex
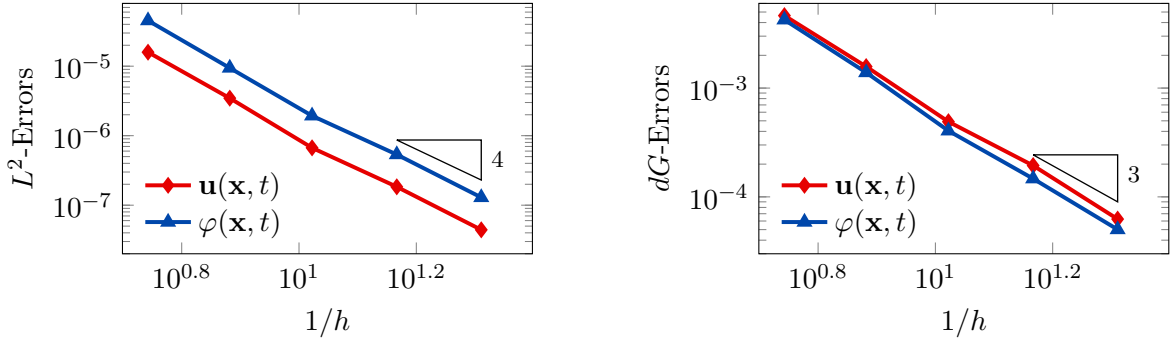
\begin{figure}[H]
\begin{subfigure}{1\textwidth}
\centering
\begin{subfigure}[b]{0.49\textwidth}
\begin{tikzpicture}
\begin{axis}[%
width=0.65\textwidth,
height=0.4\textwidth,
at={(0\textwidth,0\textwidth)},
scale only axis,
xmode=log,
xmin=5,
xmax=25,
xminorticks=true,
xlabel={$1/h$},
ymode=log,
ymin=2e-08,
ymax=8e-5,
yminorticks=true,
ylabel={$L^2$-Errors},
legend style={draw=none,fill=none,legend cell align=left},
legend pos=south west
]
\addplot [color=myred,solid,line width=1.5pt, mark=diamond*,mark options={color=myred}]
  table[row sep=crcr]{
    5.5278   1.5892e-05 \\
    7.6209   3.4643e-06 \\
    10.5252   6.6847e-07 \\
    14.679   1.8384e-07 \\
    20.4507   4.4213e-08 \\
};
\addlegendentry{$\mathbf{u}(\mathbf{x},t)$}

\addplot [color=myblue,solid,line width=1.5pt,mark=triangle*,mark options={color=myblue}]
  table[row sep=crcr]{
    5.5278   4.5291e-05 \\
    7.6209   9.4187e-06 \\
    10.5252   1.9331e-06 \\
    14.679   5.3528e-07 \\
    20.4507   1.2926e-07 \\
};
\addlegendentry{$\varphi(\mathbf{x},t)$}

\addplot [color=black,solid,line width=0.5pt]
  table[row sep=crcr]{
 14.679       8.6651e-07\\
 20.4507      8.6651e-07 \\
 20.4507      2.3e-7 \\
 14.679       8.6651e-07 \\ 
};
\node[right, align=left, text=black, font=\footnotesize]
at (axis cs:20.4507,4.5e-07) {4}; 

\end{axis}
\end{tikzpicture}
\end{subfigure}
\begin{subfigure}[b]{0.49\textwidth}
\begin{tikzpicture}
\begin{axis}[%
width=0.65\textwidth,
height=0.4\textwidth,
at={(0\textwidth,0\textwidth)},
scale only axis,
xmode=log,
xmin=5,
xmax=25,
xminorticks=true,
xlabel={$1/h$},
ymode=log,
ymin=3e-5,
ymax=0.006,
yminorticks=true,
ylabel={$dG$-Errors},
legend style={draw=none,fill=none,legend cell align=left},
legend pos=south west
]
\addplot [color=myred,solid,line width=1.5pt, mark=diamond*,mark options={color=myred}]
  table[row sep=crcr]{
    5.5278   0.0046326 \\
    7.6209   0.0015763 \\
    10.5252   0.00049093 \\
    14.679   0.000194 \\
    20.4507   6.249e-05 \\
};
\addlegendentry{$\mathbf{u}(\mathbf{x},t)$}

\addplot [color=myblue,solid,line width=1.5pt,mark=triangle*,mark options={color=myblue}]
  table[row sep=crcr]{
    5.5278   0.0042431 \\
    7.6209   0.0013903 \\
    10.5252   0.0004049 \\
    14.679   0.00014692 \\
    20.4507   5.0114e-05 \\
};
\addlegendentry{$\varphi(\mathbf{x},t)$}

\addplot [color=black,solid,line width=0.5pt]
  table[row sep=crcr]{
  14.679      2.4338e-04 \\
  20.4507      2.4338e-04 \\
  20.4507      9e-5 \\
  14.679      2.4338e-04 \\ 
};
\node[right, align=left, text=black, font=\footnotesize]
at (axis cs:20.4507,1.6e-04) {3}; 

\end{axis}
\end{tikzpicture}
\end{subfigure}
\end{subfigure}
\caption{Convergence test vs $h$ ($\tau_1 = 1$): computed errors in $L^2$-norm (left) and $dG$-norm (right) versus $1/h$ (\textit{log-log} scale). The errors are computed at the final time $T_f$. The polynomial degree of approximation is taken as $\ell = 3$.}
\label{fig:convh}
\end{figure}

%% file: IMG/ConvH_tau0.tex
\begin{figure}[H]
\begin{subfigure}{1\textwidth}
\centering
\begin{subfigure}[b]{0.49\textwidth}
\begin{tikzpicture}
\begin{axis}[%
width=0.65\textwidth,
height=0.4\textwidth,
at={(0\textwidth,0\textwidth)},
scale only axis,
xmode=log,
xmin=5,
xmax=25,
xminorticks=true,
xlabel={$1/h$},
ymode=log,
ymin=1e-08,
ymax=3e-5,
yminorticks=true,
ylabel={$L^2$-Errors},
legend style={draw=none,fill=none,legend cell align=left},
legend pos=south west
]
\addplot [color=myred,solid,line width=1.5pt, mark=diamond*,mark options={color=myred}]
  table[row sep=crcr]{
    5.5278   1.5888e-05 \\
    7.6209   3.4635e-06 \\
    10.5252   6.6831e-07 \\
    14.679   1.838e-07 \\
    20.4507   4.4203e-08 \\
};
\addlegendentry{$\mathbf{u}(\mathbf{x},t)$}

\addplot [color=myblue,solid,line width=1.5pt,mark=triangle*,mark options={color=myblue}]
  table[row sep=crcr]{
    5.5278   4.751e-06 \\
    7.6209   1.0544e-06 \\
    10.5252   2.3421e-07 \\
    14.679   5.1347e-08 \\
    20.4507   1.4405e-08 \\
};
\addlegendentry{$\varphi(\mathbf{x},t)$}

\addplot [color=black,solid,line width=0.5pt]
  table[row sep=crcr]{
 14.679       2.6372e-07\\
 20.4507      2.6372e-07 \\
 20.4507      7e-8 \\
 14.679       2.6372e-07 \\ 
};
\node[right, align=left, text=black, font=\footnotesize]
at (axis cs:20.4507,1.4e-07) {4}; 

\end{axis}
\end{tikzpicture}
\end{subfigure}
\begin{subfigure}[b]{0.49\textwidth}
\begin{tikzpicture}
\begin{axis}[%
width=0.65\textwidth,
height=0.4\textwidth,
at={(0\textwidth,0\textwidth)},
scale only axis,
xmode=log,
xmin=5,
xmax=25,
xminorticks=true,
xlabel={$1/h$},
ymode=log,
ymin=1e-5,
ymax=0.006,
yminorticks=true,
ylabel={$dG$-Errors},
legend style={draw=none,fill=none,legend cell align=left},
legend pos=south west
]
\addplot [color=myred,solid,line width=1.5pt, mark=diamond*,mark options={color=myred}]
  table[row sep=crcr]{
    5.5278   0.0046314 \\
    7.6209   0.0015759 \\
    10.5252   0.00049082 \\
    14.679   0.00019396 \\
    20.4507   6.2475e-05 \\
};
\addlegendentry{$\mathbf{u}(\mathbf{x},t)$}

\addplot [color=myblue,solid,line width=1.5pt,mark=triangle*,mark options={color=myblue}]
  table[row sep=crcr]{
    5.5278   0.00077903 \\
    7.6209   0.00032593 \\
    10.5252   0.00012019 \\
    14.679   3.3458e-05 \\
    20.4507   1.5042e-05 \\
};
\addlegendentry{$\varphi(\mathbf{x},t)$}

\addplot [color=black,solid,line width=0.5pt]
  table[row sep=crcr]{
  14.679      2.4338e-04 \\
  20.4507      2.4338e-04 \\
  20.4507      9e-5 \\
  14.679      2.4338e-04 \\ 
};
\node[right, align=left, text=black, font=\footnotesize]
at (axis cs:20.4507,1.6e-04) {3}; 

\end{axis}
\end{tikzpicture}
\end{subfigure}
\end{subfigure}
\caption{Convergence test vs $h$ ($\tau_1 = 0$): computed errors in $L^2$-norm (left) and $dG$-norm (right) versus $1/h$ (\textit{log-log} scale). The errors are computed at the final time $T_f$. The polynomial degree of approximation is taken as $\ell = 3$.}
\label{fig:convh_tau0}
\end{figure}
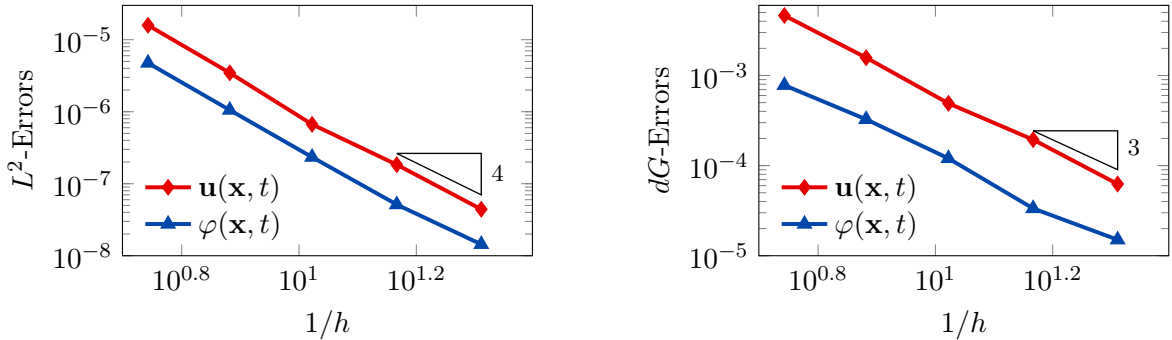

%% file: IMG/ConvP.tex
\begin{figure}[H]
\begin{subfigure}{1\textwidth}
\centering
\begin{subfigure}[b]{0.49\textwidth}
\begin{tikzpicture}
\begin{axis}[%
width=0.65\textwidth,
height=0.4\textwidth,
at={(0\textwidth,0\textwidth)},
scale only axis,
xmin=0,
xmax=5.5,
xminorticks=true,
xlabel={$\ell$},
ymode=log,
ymin=5e-09,
ymax=0.09,
yminorticks=true,
ylabel={$L^2$-Errors},
legend style={draw=none,fill=none,legend cell align=left},
legend pos=south west
]
\addplot [color=myred,solid,line width=1.5pt, mark=diamond*,mark options={color=myred}]
  table[row sep=crcr]{
1   0.013201 \\
2   0.00041881 \\
3   1.5892e-05 \\
4   6.2229e-07 \\
5   1.3908e-08 \\
};
\addlegendentry{$\mathbf{u}(\mathbf{x},t)$}

\addplot [color=myblue,solid,line width=1.5pt,mark=square*,mark options={color=myblue}]
  table[row sep=crcr]{
1   0.029957 \\
2   0.00095666 \\
3   4.5291e-05 \\
4   1.6893e-06 \\
5   4.0119e-08 \\
};
\addlegendentry{$\varphi(\mathbf{x},t)$}

\addplot [color=black,dashed,line width=1.0pt]
  table[row sep=crcr]{
    1.5 11.1090e-003
    2   2.4788e-03 \\
    3   123.4098e-06 \\
    4   6.1442e-06 \\
    5   305.9023e-09 \\
};
\addlegendentry{$e^{-3 \ell}$}

\end{axis}
\end{tikzpicture}
\end{subfigure}
\begin{subfigure}[b]{0.49\textwidth}
\begin{tikzpicture}
\begin{axis}[%
width=0.65\textwidth,
height=0.4\textwidth,
at={(0\textwidth,0\textwidth)},
scale only axis,
xmin=0,
xmax=5.5,
xminorticks=true,
xlabel={$\ell$},
ymode=log,
ymin=3e-6,
ymax=0.8,
yminorticks=true,
ylabel={$dG$-Errors},
legend style={draw=none,fill=none,legend cell align=left},
legend pos=south west
]
\addplot [color=myred,solid,line width=1.5pt, mark=diamond*,mark options={color=myred}]
  table[row sep=crcr]{
1   0.20768 \\
2   0.055365 \\
3   0.0046326 \\
4   0.00023503 \\
5   9.1285e-06 \\
};
\addlegendentry{$\mathbf{u}(\mathbf{x},t)$}

\addplot [color=myblue,solid,line width=1.5pt,mark=square*,mark options={color=myblue}]
  table[row sep=crcr]{
1   0.48053 \\
2   0.067188 \\
3   0.0042431 \\
4   0.00019592 \\
5   7.1101e-06 \\
};
\addlegendentry{$\varphi(\mathbf{x},t)$}

\addplot [color=black,dashed,line width=1.0pt]
  table[row sep=crcr]{
    1.5 1.110899653824231
    2   0.247875217666636 \\
    3   0.012340980408668 \\
    4   0.000614421235333 \\
    5   0.000030590232050 \\
};
\addlegendentry{$e^{-3 \ell}$}

\end{axis}
\end{tikzpicture}
\end{subfigure}
\end{subfigure}
\caption{Convergence test vs $\ell$ ($\tau_1 = 1$): computed errors in $L^2$-norm (left) and $dG$-norm (right) versus $\ell$ (\textit{semi-log} scale). The errors are computed at the final time $T_f$. The computational mesh is made of $100$ polygons.}
\label{fig:convp}
\end{figure}
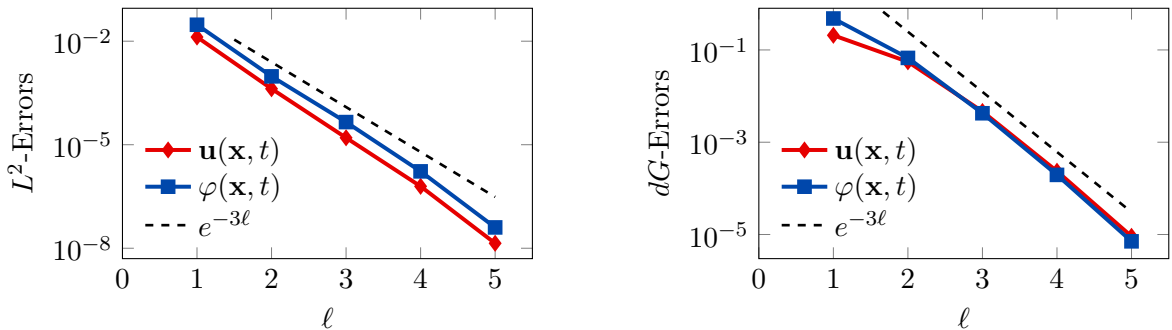

%% file: IMG/ConvP_tau0.tex
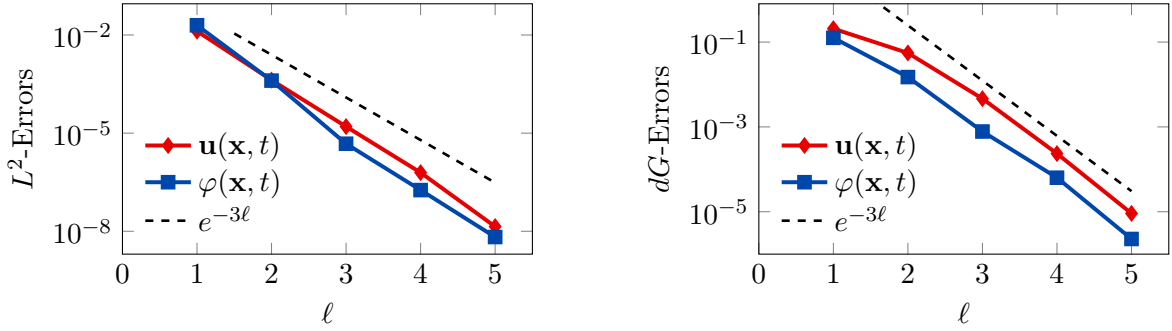
\begin{figure}[H]
\begin{subfigure}{1\textwidth}
\centering
\begin{subfigure}[b]{0.49\textwidth}
\begin{tikzpicture}
\begin{axis}[%
width=0.65\textwidth,
height=0.4\textwidth,
at={(0\textwidth,0\textwidth)},
scale only axis,
xmin=0,
xmax=5.5,
xminorticks=true,
xlabel={$\ell$},
ymode=log,
ymin=2e-09,
ymax=0.09,
yminorticks=true,
ylabel={$L^2$-Errors},
legend style={draw=none,fill=none,legend cell align=left},
legend pos=south west
]
\addplot [color=myred,solid,line width=1.5pt, mark=diamond*,mark options={color=myred}]
  table[row sep=crcr]{
1   0.013094 \\
2   0.00041841 \\
3   1.5888e-05 \\
4   6.2213e-07 \\
5   1.3904e-08 \\
};
\addlegendentry{$\mathbf{u}(\mathbf{x},t)$}

\addplot [color=myblue,solid,line width=1.5pt,mark=square*,mark options={color=myblue}]
  table[row sep=crcr]{
1   0.019847 \\
2   0.00040333 \\
3   4.751e-06 \\
4   1.8169e-07 \\
5   6.6642e-09 \\
};
\addlegendentry{$\varphi(\mathbf{x},t)$}

\addplot [color=black,dashed,line width=1.0pt]
  table[row sep=crcr]{
    1.5 11.1090e-003
    2   2.4788e-03 \\
    3   123.4098e-06 \\
    4   6.1442e-06 \\
    5   305.9023e-09 \\
};
\addlegendentry{$e^{-3 \ell}$}

\end{axis}
\end{tikzpicture}
\end{subfigure}
\begin{subfigure}[b]{0.49\textwidth}
\begin{tikzpicture}
\begin{axis}[%
width=0.65\textwidth,
height=0.4\textwidth,
at={(0\textwidth,0\textwidth)},
scale only axis,
xmin=0,
xmax=5.5,
xminorticks=true,
xlabel={$\ell$},
ymode=log,
ymin=1e-6,
ymax=0.8,
yminorticks=true,
ylabel={$dG$-Errors},
legend style={draw=none,fill=none,legend cell align=left},
legend pos=south west
]
\addplot [color=myred,solid,line width=1.5pt, mark=diamond*,mark options={color=myred}]
  table[row sep=crcr]{
1   0.20709 \\
2   0.055349 \\
3   0.0046314 \\
4   0.00023497 \\
5   9.1262e-06 \\
};
\addlegendentry{$\mathbf{u}(\mathbf{x},t)$}

\addplot [color=myblue,solid,line width=1.5pt,mark=square*,mark options={color=myblue}]
  table[row sep=crcr]{
1   0.126 \\
2   0.014964 \\
3   0.00077903 \\
4   6.3451e-05 \\
5   2.2788e-06 \\
};
\addlegendentry{$\varphi(\mathbf{x},t)$}

\addplot [color=black,dashed,line width=1.0pt]
  table[row sep=crcr]{
    1.5 1.110899653824231
    2   0.247875217666636 \\
    3   0.012340980408668 \\
    4   0.000614421235333 \\
    5   0.000030590232050 \\
};
\addlegendentry{$e^{-3 \ell}$}

\end{axis}
\end{tikzpicture}
\end{subfigure}
\end{subfigure}
\caption{Convergence test vs $\ell$ ($\tau_1 = 0$): computed errors in $L^2$-norm (left) and $dG$-norm (right) versus $\ell$ (\textit{semi-log} scale). The errors are computed at the final time $T_f$. The computational mesh is made of $100$ polygons.}
\label{fig:convp_tau0}
\end{figure}

%% file: IMG/ConvDT.tex
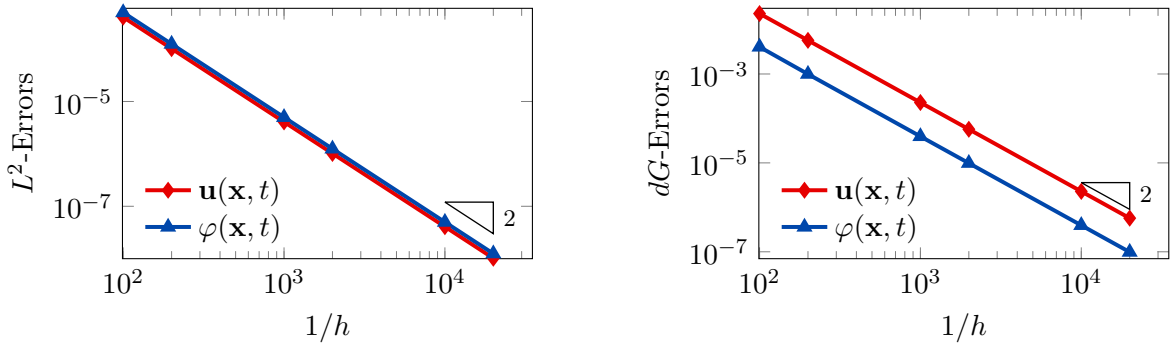
\begin{figure}[H]
\begin{subfigure}{1\textwidth}
\centering
\begin{subfigure}[b]{0.49\textwidth}
\begin{tikzpicture}
\begin{axis}[%
width=0.65\textwidth,
height=0.4\textwidth,
at={(0\textwidth,0\textwidth)},
scale only axis,
xmode=log,
xmin=99,
xmax=35000,
xminorticks=true,
xlabel={$1/h$},
ymode=log,
ymin=1e-8,
ymax=6e-4,
yminorticks=true,
ylabel={$L^2$-Errors},
legend style={draw=none,fill=none,legend cell align=left},
legend pos=south west
]
\addplot [color=myred,solid,line width=1.5pt, mark=diamond*,mark options={color=myred}]
  table[row sep=crcr]{
100   0.00040811 \\
200   0.00010202 \\
1000   4.0807e-06 \\
2000   1.0202e-06 \\
10000   4.0807e-08 \\
20000   1.0202e-08 \\
};
\addlegendentry{$\mathbf{u}(\mathbf{x},t)$}

\addplot [color=myblue,solid,line width=1.5pt,mark=triangle*,mark options={color=myblue}]
  table[row sep=crcr]{
100   0.00049922 \\
200   0.00012478 \\
1000   4.9907e-06 \\
2000   1.2477e-06 \\
10000   4.9907e-08 \\
20000   1.2477e-08 \\
};
\addlegendentry{$\varphi(\mathbf{x},t)$}

\addplot [color=black,solid,line width=0.5pt]
  table[row sep=crcr]{
 10000      1.2e-07 \\
 20000      1.2e-07 \\
 20000      3e-8 \\
 10000      1.2e-07 \\ 
};
\node[right, align=left, text=black, font=\footnotesize]
at (axis cs:20000,6e-08) {2}; 

\end{axis}
\end{tikzpicture}
\end{subfigure}
\begin{subfigure}[b]{0.49\textwidth}
\begin{tikzpicture}
\begin{axis}[%
width=0.65\textwidth,
height=0.4\textwidth,
at={(0\textwidth,0\textwidth)},
scale only axis,
xmode=log,
xmin=99,
xmax=35000,
xminorticks=true,
xlabel={$1/h$},
ymode=log,
ymin=7e-8,
ymax=0.03,
yminorticks=true,
ylabel={$dG$-Errors},
legend style={draw=none,fill=none,legend cell align=left},
legend pos=south west
]
\addplot [color=myred,solid,line width=1.5pt, mark=diamond*,mark options={color=myred}]
  table[row sep=crcr]{
100   0.022822 \\
200   0.0057051 \\
1000   0.0002282 \\
2000   5.7049e-05 \\
10000   2.282e-06 \\
20000   5.7049e-07 \\
};
\addlegendentry{$\mathbf{u}(\mathbf{x},t)$}

\addplot [color=myblue,solid,line width=1.5pt,mark=triangle*,mark options={color=myblue}]
  table[row sep=crcr]{
100   0.0040608 \\
200   0.00099301 \\
1000   3.942e-05 \\
2000   9.853e-06 \\
10000   3.9409e-07 \\
20000   9.8524e-08 \\
};
\addlegendentry{$\varphi(\mathbf{x},t)$}

\addplot [color=black,solid,line width=0.5pt]
  table[row sep=crcr]{
 10000      3.6e-06 \\
 20000      3.6e-06 \\
 20000      9e-7 \\
 10000      3.6e-06 \\ 
};
\node[right, align=left, text=black, font=\footnotesize]
at (axis cs:20000,2e-06) {2};  

\end{axis}
\end{tikzpicture}
\end{subfigure}
\end{subfigure}
\caption{Convergence test vs $\Delta t$ ($\tau_1 = 1$): computed errors in $L^2$-norm (left) and $dG$-norm (right) versus $1/h$ (\textit{log-log} scale). The errors are computed at the final time $T_f$. The time-marching scheme is the Newmark-$\beta$ method.}
\label{fig:convdt}
\end{figure}

%% file: IMG/ConvDT_tau0.tex
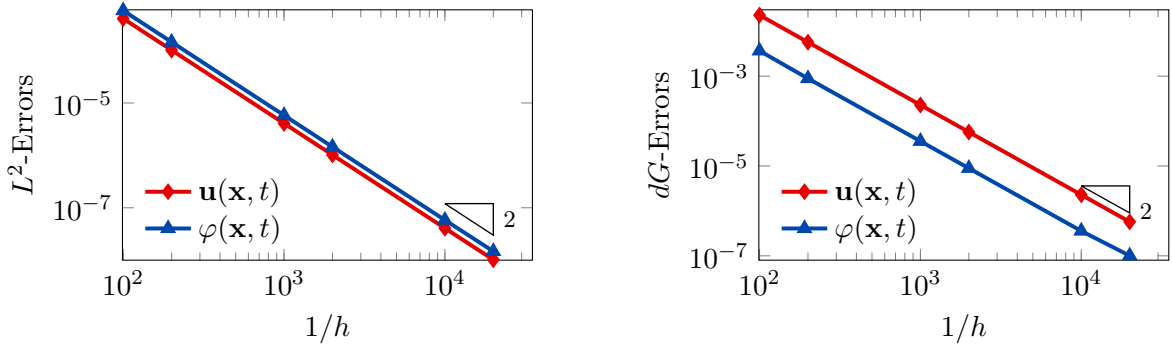
\begin{figure}[H]
\begin{subfigure}{1\textwidth}
\centering
\begin{subfigure}[b]{0.49\textwidth}
\begin{tikzpicture}
\begin{axis}[%
width=0.65\textwidth,
height=0.4\textwidth,
at={(0\textwidth,0\textwidth)},
scale only axis,
xmode=log,
xmin=99,
xmax=35000,
xminorticks=true,
xlabel={$1/h$},
ymode=log,
ymin=1e-8,
ymax=6e-4,
yminorticks=true,
ylabel={$L^2$-Errors},
legend style={draw=none,fill=none,legend cell align=left},
legend pos=south west
]
\addplot [color=myred,solid,line width=1.5pt, mark=diamond*,mark options={color=myred}]
  table[row sep=crcr]{
100   0.00040817 \\
200   0.00010204 \\
1000   4.0814e-06 \\
2000   1.0203e-06 \\
10000   4.0813e-08 \\
20000   1.021e-08 \\
};
\addlegendentry{$\mathbf{u}(\mathbf{x},t)$}

\addplot [color=myblue,solid,line width=1.5pt,mark=triangle*,mark options={color=myblue}]
  table[row sep=crcr]{
100   0.00058079 \\
200   0.00014532 \\
1000   5.8145e-06 \\
2000   1.4536e-06 \\
10000   5.8254e-08 \\
20000   1.4741e-08 \\
};
\addlegendentry{$\varphi(\mathbf{x},t)$}

\addplot [color=black,solid,line width=0.5pt]
  table[row sep=crcr]{
 10000      1.2e-07 \\
 20000      1.2e-07 \\
 20000      3e-8 \\
 10000      1.2e-07 \\ 
};
\node[right, align=left, text=black, font=\footnotesize]
at (axis cs:20000,6e-08) {2}; 

\end{axis}
\end{tikzpicture}
\end{subfigure}
\begin{subfigure}[b]{0.49\textwidth}
\begin{tikzpicture}
\begin{axis}[%
width=0.65\textwidth,
height=0.4\textwidth,
at={(0\textwidth,0\textwidth)},
scale only axis,
xmode=log,
xmin=99,
xmax=35000,
xminorticks=true,
xlabel={$1/h$},
ymode=log,
ymin=8e-8,
ymax=0.03,
yminorticks=true,
ylabel={$dG$-Errors},
legend style={draw=none,fill=none,legend cell align=left},
legend pos=south west
]
\addplot [color=myred,solid,line width=1.5pt, mark=diamond*,mark options={color=myred}]
  table[row sep=crcr]{
100   0.022822 \\
200   0.0057051 \\
1000   0.0002282 \\
2000   5.705e-05 \\
10000   2.282e-06 \\
20000   5.7052e-07 \\
};
\addlegendentry{$\mathbf{u}(\mathbf{x},t)$}

\addplot [color=myblue,solid,line width=1.5pt,mark=triangle*,mark options={color=myblue}]
  table[row sep=crcr]{
100   0.0036683 \\
200   0.00088697 \\
1000   3.5603e-05 \\
2000   8.9007e-06 \\
10000   3.5708e-07 \\
20000   1.0045e-07 \\
};
\addlegendentry{$\varphi(\mathbf{x},t)$}

\addplot [color=black,solid,line width=0.5pt]
  table[row sep=crcr]{
 10000      3.6e-06 \\
 20000      3.6e-06 \\
 20000      9e-7 \\
 10000      3.6e-06 \\ 
};
\node[right, align=left, text=black, font=\footnotesize]
at (axis cs:20000,1e-06) {2}; 

\end{axis}
\end{tikzpicture}
\end{subfigure}
\end{subfigure}
\caption{Convergence test vs $\Delta t$ ($\tau_1 = 0$): computed errors in $L^2$-norm (left) and $dG$-norm (right) versus $1/h$ (\textit{log-log} scale). The errors are computed at the final time $T_f$. The time-marching scheme is the coupled Newmark-$\beta$ and $\theta$-method.}
\label{fig:convdt_tau0}
\end{figure}

%% file: Sections/SuperConvTest.tex
The aim of this section is to prove that the PolyDG scheme proposed for this problem is not only optimal convergence, but it also shows some superconvergence properties. To this aim, we consider the following exact solutions:
\begin{equation}
    \begin{aligned}
    \mathbf{u}(x,y,t) & \ =  \nu_u \, \sin(2 \pi t) \left( \begin{aligned}
    & x^2 \sin(\pi x) \sin(\pi y) \\
    & - x^2 \sin(\pi x) \sin(\pi y)
    \end{aligned} \right), \\
    \varphi(x,y,t) & \ = \nu_{\varphi} \, \sin(\sqrt{2} \pi t) \left( x^2 \cos\left(\frac{\pi x}{2}\right) \sin(\pi x) \right),
    \end{aligned}
\end{equation}
from which we infer initial and boundary conditions, as well as forcing terms. The parameters $\nu_u$, $\nu_{\varphi}$ control the magnitude of the displacement and the generalized pressure, respectively. The model coefficients are reported in Table~\ref{tab:params_convtest} (chosen as in the convergence test). We observe that the optimal convergence property has been already proven in Section~\ref{sec:convtest} by setting $\nu_u = \nu_{\varphi} = 1$.

To observe the better robustness of the scheme with respect to large pressures, we propose two different tests. In the first, we set $\nu_u = 0.1$, $\nu_{\varphi} = \num[exponent-product=\ensuremath{\cdot}, print-unity-mantissa=false]{1e+4}$, we consider the same sequence of meshes as in Section~\ref{sec:convtest} and polynomial degree $\ell = 2$. For the second test, we fix $\nu_u = 0.1$, the mesh, and the polynomial degree of approximation; then we vary the values of $\nu_{\varphi} = \left[ 1, \num[exponent-product=\ensuremath{\cdot}, print-unity-mantissa=false]{1e+1}, \num[exponent-product=\ensuremath{\cdot}, print-unity-mantissa=false]{1e+2}, \num[exponent-product=\ensuremath{\cdot}, print-unity-mantissa=false]{1e+3}, \num[exponent-product=\ensuremath{\cdot}, print-unity-mantissa=false]{1e+4}, \num[exponent-product=\ensuremath{\cdot}, print-unity-mantissa=false]{1e+5}, \num[exponent-product=\ensuremath{\cdot}, print-unity-mantissa=false]{1e+6} \right]$. We consider the following discretization parameters for the second test: $N = 400, \ \ell = 2$.

\begin{table}[ht]
    \centering 
    \begin{tabular}{c  c  c  c  c  c  c c c}
    $1/h$ & $\|\mathbf{e}^{u}\|_{L^2}$  & $\text{roc}^{u}_{L^2}$ & $\|\mathbf{e}^{u}\|_{dG}$  & $\text{roc}^{u}_{dG}$& $\|e^{p}\|_{L^2}$  & $\text{roc}^{p}_{L^2}$ & $\|e^{p}\|_{dG}$  & $\text{roc}^{p}_{dG}$  \T\B\B \\
    \hline
    5.53 & \num[exponent-product=\ensuremath{\cdot}]{4.49e-3} & - & 0.15 & - & 8.47 & - & 1411.33 & - \T\B \\
    7.62 & \num[exponent-product=\ensuremath{\cdot}]{1.06e-3} & 4.49 & 0.06 & 3.11 & 2.83 & 3.42 & 779.86 & 1.85 \T\B \\
    10.53 & \num[exponent-product=\ensuremath{\cdot}]{2.86e-4} & 4.06 & 0.02 & 3.48 & 0.85 & 3.73 & 373.68 & 2.28 \T\B \\
    14.68 & \num[exponent-product=\ensuremath{\cdot}]{9.30e-05} & 3.38 & \num[exponent-product=\ensuremath{\cdot}]{7.47e-03} & 2.69 & 0.36 & 2.62 & 196.18 & 1.94 \T\B \\
    20.45 & \num[exponent-product=\ensuremath{\cdot}]{3.12e-05} & 3.30 & \num[exponent-product=\ensuremath{\cdot}]{2.70e-3} & 3.07 & 0.12 & 3.19 & 96.92 & 2.13 \T\B \\
    \end{tabular}
    \\[5pt]
\caption{Superconvergence test: computed errors and convergence rates in $L^2$- and $dG$-norms versus $h$ using as polynomial degree of approximation $\ell = 2$.}
\label{tab:superconvh}
\end{table}

\input{IMG/Superconv_Nu.tex}

By looking at Table~\ref{tab:superconvh} we observe the superconvergence phenomenon for the displacement field. Indeed, we observe that using a polynomial degree of approximation equal to $\ell$, then the error of the displacement in $dG$-norm converges with order $\ell+1$ (we remark that the expected order is $\ell$ in this case Theorem~\ref{thm:conv}, \cite{CangianiDong2017,Cangiani2014} and this rate is observed for the generalized-pressure). Moreover, for what concerns the error in $L^2$-norm, we observe $(\ell+1)+1$ convergence rate for the first refinements.

In Figure~\ref{fig:superconvnu}, we observe the behavior of the errors with respect to increasing values of $\nu_{\varphi}$. We see that both the $L^2$- and $dG$-errors of the displacements are way lower than the errors of the generalized pressure (even for not too big values of $\nu_{\varphi}$). It is interesting to notice that, for the first tested values of $\nu_{\varphi}$, the displacement errors remain almost constant while the errors of the generalized pressure start growing.

%% file: IMG/Superconv_Nu.tex
\begin{center}
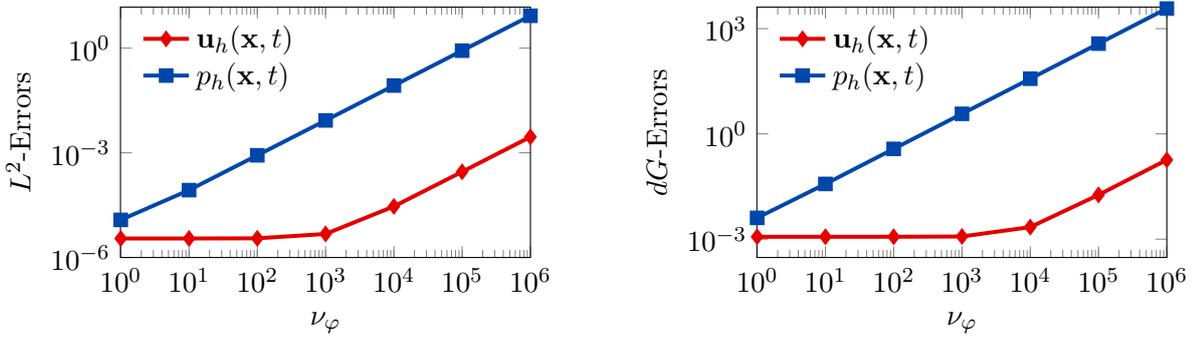
\begin{figure}[H]
\begin{subfigure}[b]{0.49\textwidth}
\begin{tikzpicture}
\begin{axis}[%
width=0.65\textwidth,
height=0.4\textwidth,
at={(0\textwidth,0\textwidth)},
scale only axis,
xmode=log,
xmin=0.99,
xmax=1000001,
xminorticks=true,
xlabel={$\nu_{\varphi}$},
ymode=log,
ymin=1e-6,
ymax=15,
yminorticks=true,
ylabel={$L^2$-Errors},
legend style={draw=none,fill=none,legend cell align=left},
legend pos=north west
]
\addplot [color=myred,solid,line width=1.5pt, mark=diamond*,mark options={color=myred}]
  table[row sep=crcr]{
1   3.5054e-06 \\
10   3.5079e-06 \\
100   3.5437e-06 \\
1000   4.7303e-06 \\
10000   2.9151e-05 \\
100000   0.00028644 \\
1000000   0.0028612 \\
};
\addlegendentry{$\mathbf{u}_h(\mathbf{x},t)$}

\addplot [color=myblue,solid,line width=1.5pt,mark=square*,mark options={color=myblue}]
  table[row sep=crcr]{
1   1.1964e-05 \\
10   8.5295e-05 \\
100   0.00084852 \\
1000   0.0084844 \\
10000   0.084844 \\
100000   0.84844 \\
1000000   8.4844 \\
};
\addlegendentry{$\varphi_h(\mathbf{x},t)$}
\end{axis}
\end{tikzpicture}
\end{subfigure}
\begin{subfigure}[b]{0.49\textwidth}
\begin{tikzpicture}
\begin{axis}[%
width=0.65\textwidth,
height=0.4\textwidth,
at={(0\textwidth,0\textwidth)},
scale only axis,
xmode=log,
xmin=0.99,
xmax=1000001,
xminorticks=true,
xlabel={$\nu_{\varphi}$},
ymode=log,
ymin=3e-4,
ymax=4100,
yminorticks=true,
ylabel={$dG$-Errors},
legend style={draw=none,fill=none,legend cell align=left},
legend pos=north west
]
\addplot [color=myred,solid,line width=1.5pt, mark=diamond*,mark options={color=myred}]
  table[row sep=crcr]{
1   0.0011687 \\
10   0.0011687 \\
100   0.0011696 \\
1000   0.0011909 \\
10000   0.0022048 \\
100000   0.018256 \\
1000000   0.1817 \\
};
\addlegendentry{$\mathbf{u}_h(\mathbf{x},t)$}

\addplot [color=myblue,solid,line width=1.5pt,mark=square*,mark options={color=myblue}]
  table[row sep=crcr]{
1   0.0040383 \\
10   0.037449 \\
100   0.37373 \\
1000   3.7368 \\
10000   37.3677 \\
100000   373.6762 \\
1000000   3736.7617 \\
};
\addlegendentry{$\varphi_h(\mathbf{x},t)$}
\end{axis}
\end{tikzpicture}
\end{subfigure}
\caption{Superconvergence test: computed errors in $L^2$-norm (left) and $dG$-norm (right) versus $\nu_{\varphi}$ (\textit{log-log} scale) using as polynomial degrees of approximation and number of elements: $\ell = 2, N = 400$.}
\label{fig:superconvnu}
\end{figure}
\end{center}

%% file: Sections/RobTest.tex
In this section, we address the problem of testing the scheme's robustness with respect to some of the physical parameters appearing in \eqref{eq:system}. In particular, we are interested in testing the robustness of our method with respect to the coefficients for which the stability estimate in Theorem~\eqref{thm:stabdisc} is not robust: low values of the permeability $\mathbf{D}$, low values of the second Lamé parameter $\mu$, and low values of the viscoelasticity retardation time $\delta_1$. In addition to that, we are interested in assessing the performance of the scheme with respect to high values of secondary consolidation $\delta_2 \lambda$, possibly considering low values of the storage coefficient $d_0$. Last, we test our method in the case of strong coupling between the momentum equation and the mass/energy one, by considering degenerating values of the first Lamé parameter $\lambda$. The values of the other model parameters are as in Table~\ref{tab:params_convtest}. The time discretization parameters are $T_f = 0.1$, $\Delta t = \num[exponent-product=\ensuremath{\cdot}, print-unity-mantissa=false]{5e-5}$. For the robustness analysis, a sequence of polygonal meshes as the one in Figure~\ref{fig:voronoi_mesh} is considered, and we consider polynomial degree $\ell = 3$.

\begin{table}[H]
    \centering 
    \small
    \begin{tabular}{c  c  c  c  c  c  c c c}
    $1/h$ & $\|\mathbf{e}^{u}\|_{L^2}$  & $\text{roc}^{u}_{L^2}$ & $\|\mathbf{e}^{u}\|_{dG}$  & $\text{roc}^{u}_{dG}$& $\|e^{p}\|_{L^2}$  & $\text{roc}^{p}_{L^2}$ & $\|e^{p}\|_{dG}$  & $\text{roc}^{p}_{dG}$  \T\B\B \\
    \hline
    5.53 & \num[exponent-product=\ensuremath{\cdot}]{1.58e-05} & - & \num[exponent-product=\ensuremath{\cdot}]{4.63e-3} & - & \num[exponent-product=\ensuremath{\cdot}]{5.98e-4} & - & \num[exponent-product=\ensuremath{\cdot}]{1.77e-4} & - \T\B \\
    7.62 & \num[exponent-product=\ensuremath{\cdot}]{3.46e-06} & 4.74 & \num[exponent-product=\ensuremath{\cdot}]{1.58e-3} & 3.36 & \num[exponent-product=\ensuremath{\cdot}]{1.89e-4} & 3.59 & \num[exponent-product=\ensuremath{\cdot}]{7.22e-05} & 2.80 \T\B \\
    10.53 & \num[exponent-product=\ensuremath{\cdot}]{6.67e-07} & 5.10 & \num[exponent-product=\ensuremath{\cdot}]{4.91e-4} & 3.61 & \num[exponent-product=\ensuremath{\cdot}]{5.19e-05} & 4.01 & \num[exponent-product=\ensuremath{\cdot}]{2.95e-05} & 2.77 \T\B \\
    14.68 & \num[exponent-product=\ensuremath{\cdot}]{1.84e-07} & 3.88 & \num[exponent-product=\ensuremath{\cdot}]{1.94e-4} & 2.79 & \num[exponent-product=\ensuremath{\cdot}]{2.00e-05} & 2.87 & \num[exponent-product=\ensuremath{\cdot}]{1.45e-05} & 2.14 \T\B \\
    20.45 & \num[exponent-product=\ensuremath{\cdot}]{4.42e-08} & 4.30 & \num[exponent-product=\ensuremath{\cdot}]{6.25e-05} & 3.42 & \num[exponent-product=\ensuremath{\cdot}]{6.58e-06} & 3.35 & \num[exponent-product=\ensuremath{\cdot}]{6.90e-06} & 2.24 \T\B \\
    \end{tabular}
    \\[5pt]
\caption{Robustness test vs conductivity: polynomial degree $\ell = 3$. The parameters are chosen as in Table~\ref{tab:params_convtest} and the conductivity is set as $\mathbf{D} = \num[exponent-product=\ensuremath{\cdot}, print-unity-mantissa=false]{1e-06}\mathbf{I}$.}
\label{tab:rob_vs_conductivity}
\end{table}

\begin{table}[H]
    \centering 
    \small
    \begin{tabular}{c  c  c  c  c  c  c c c}
    $1/h$ & $\|\mathbf{e}^{u}\|_{L^2}$  & $\text{roc}^{u}_{L^2}$ & $\|\mathbf{e}^{u}\|_{dG}$  & $\text{roc}^{u}_{dG}$& $\|e^{p}\|_{L^2}$  & $\text{roc}^{p}_{L^2}$ & $\|e^{p}\|_{dG}$  & $\text{roc}^{p}_{dG}$  \T\B\B \\
    \hline
    5.53 & \num[exponent-product=\ensuremath{\cdot}]{3.02e-05} & - & \num[exponent-product=\ensuremath{\cdot}]{2.184e-3} & - & \num[exponent-product=\ensuremath{\cdot}]{1.95e-05} & - & \num[exponent-product=\ensuremath{\cdot}]{2.73e-3} & - \T\B \\
    7.62 & \num[exponent-product=\ensuremath{\cdot}]{1.13e-05} & 3.06 & \num[exponent-product=\ensuremath{\cdot}]{7.59e-4} & 3.29 & \num[exponent-product=\ensuremath{\cdot}]{5.54e-06} & 3.92 & \num[exponent-product=\ensuremath{\cdot}]{1.02e-3} & 3.05 \T\B \\
    10.53 & \num[exponent-product=\ensuremath{\cdot}]{3.59e-06} & 3.55 & \num[exponent-product=\ensuremath{\cdot}]{2.33e-4} & 3.66 & \num[exponent-product=\ensuremath{\cdot}]{1.27e-06} & 4.56 & \num[exponent-product=\ensuremath{\cdot}]{3.25e-4} & 3.55 \T\B \\
    14.68 & \num[exponent-product=\ensuremath{\cdot}]{9.33e-07} & 4.05 & \num[exponent-product=\ensuremath{\cdot}]{9.06e-05} & 2.84 & \num[exponent-product=\ensuremath{\cdot}]{3.61e-07} & 3.78 & \num[exponent-product=\ensuremath{\cdot}]{1.16e-4} & 3.11 \T\B \\
    20.45 & \num[exponent-product=\ensuremath{\cdot}]{3.53e-07} & 2.93 & \num[exponent-product=\ensuremath{\cdot}]{2.92e-05} & 3.41 & \num[exponent-product=\ensuremath{\cdot}]{8.31e-08} & 4.43 & \num[exponent-product=\ensuremath{\cdot}]{3.97e-05} & 3.22
    \end{tabular}
    \\[5pt]
\caption{Robustness test vs Lamé parameter $\mu$: polynomial degree $\ell = 3$. The parameters are chosen as in Table~\ref{tab:params_convtest}, the second Lamé parameter is set as $\mu = \num[exponent-product=\ensuremath{\cdot}, print-unity-mantissa=false]{1e-06}$.}
\label{tab:rob_vs_mu}
\end{table}

\begin{table}[H]
    \centering 
    \small
    \begin{tabular}{c  c  c  c  c  c  c c c}
    $1/h$ & $\|\mathbf{e}^{u}\|_{L^2}$  & $\text{roc}^{u}_{L^2}$ & $\|\mathbf{e}^{u}\|_{dG}$  & $\text{roc}^{u}_{dG}$& $\|e^{p}\|_{L^2}$  & $\text{roc}^{p}_{L^2}$ & $\|e^{p}\|_{dG}$  & $\text{roc}^{p}_{dG}$  \T\B\B \\
    \hline
    5.53 & \num[exponent-product=\ensuremath{\cdot}]{2.11e-05} & - & \num[exponent-product=\ensuremath{\cdot}]{4.76e-3} & - & \num[exponent-product=\ensuremath{\cdot}]{2.30e-05} & - & \num[exponent-product=\ensuremath{\cdot}]{2.94e-3} & - \T\B \\
    7.62 & \num[exponent-product=\ensuremath{\cdot}]{4.43e-06} & 4.86 & \num[exponent-product=\ensuremath{\cdot}]{1.61e-3} & 3.38 & \num[exponent-product=\ensuremath{\cdot}]{5.68e-06} & 4.35 & \num[exponent-product=\ensuremath{\cdot}]{1.04e-3} & 3.25 \T\B \\
    10.53 & \num[exponent-product=\ensuremath{\cdot}]{8.96e-07} & 4.95 & \num[exponent-product=\ensuremath{\cdot}]{5.02e-4} & 3.61 & \num[exponent-product=\ensuremath{\cdot}]{1.32e-06} & 4.52 & \num[exponent-product=\ensuremath{\cdot}]{3.33e-4} & 3.51 \T\B \\
    14.68 & \num[exponent-product=\ensuremath{\cdot}]{2.27e-07} & 4.13 & \num[exponent-product=\ensuremath{\cdot}]{1.97e-4} & 2.81 & \num[exponent-product=\ensuremath{\cdot}]{3.69e-07} & 3.82 & \num[exponent-product=\ensuremath{\cdot}]{1.17e-4} & 3.14 \T\B \\
    20.45 & \num[exponent-product=\ensuremath{\cdot}]{5.55e-08} & 4.24 & \num[exponent-product=\ensuremath{\cdot}]{6.39e-05} & 3.40 & \num[exponent-product=\ensuremath{\cdot}]{8.50e-08} & 4.43 & \num[exponent-product=\ensuremath{\cdot}]{4.03e-05} & 3.22 \T\B \\
    \end{tabular}
    \\[5pt]
\caption{Robustness test vs viscoelasticity retardation time $\delta_1$: polynomial degree $\ell = 3$. The parameters are chosen as in Table~\ref{tab:params_convtest}, the viscoelasticity retardation time is set as $\delta_1 = \num[exponent-product=\ensuremath{\cdot}, print-unity-mantissa=false]{1e-06}$.}
\label{tab:rob_vs_delta1}
\end{table}

First, comparing the results for low values of $\mathbf{D}$, $\mu$, and $\delta_1$ (cf. Table~\ref{tab:rob_vs_conductivity}, Table~\ref{tab:rob_vs_mu}, Table~\ref{tab:rob_vs_delta1}, respectively) with the error estimate presented in Theorem~\ref{thm:conv} and with classical error estimates of the PolyDG methods (see also Section~\ref{sec:convtest}), we observe that our scheme is robust with respect to these configuration. We observe that we lose $\ell+1$ accuracy in $L^2$-norm for the generalized pressure in the case $\mathbf{D} \ll 1$ and for the displacement in the case $\mu \ll 1$, but we still observe convergence of order $\ell$ in both cases. 

\begin{table}[H]
    \centering 
    \small
    \begin{tabular}{c  c  c  c  c  c  c c c}
    $1/h$ & $\|\mathbf{e}^{u}\|_{L^2}$  & $\text{roc}^{u}_{L^2}$ & $\|\mathbf{e}^{u}\|_{dG}$  & $\text{roc}^{u}_{dG}$& $\|e^{p}\|_{L^2}$  & $\text{roc}^{p}_{L^2}$ & $\|e^{p}\|_{dG}$  & $\text{roc}^{p}_{dG}$  \T\B\B \\
    \hline
    5.53 & \num[exponent-product=\ensuremath{\cdot}]{3.74e-4} & - & 0.02 & - & \num[exponent-product=\ensuremath{\cdot}]{2.00e-05} & - & \num[exponent-product=\ensuremath{\cdot}]{2.72e-3} & - \T\B \\
    7.62 & \num[exponent-product=\ensuremath{\cdot}]{1.31e-4} & 3.28 & \num[exponent-product=\ensuremath{\cdot}]{8.99e-3} & 1.90 & \num[exponent-product=\ensuremath{\cdot}]{5.63e-06} & 3.94 & \num[exponent-product=\ensuremath{\cdot}]{1.03e-3} & 3.02 \T\B \\
    10.53 & \num[exponent-product=\ensuremath{\cdot}]{3.36e-05} & 4.21 & \num[exponent-product=\ensuremath{\cdot}]{3.20e-3} & 3.20 & \num[exponent-product=\ensuremath{\cdot}]{1.28e-06} & 4.58 & \num[exponent-product=\ensuremath{\cdot}]{3.28e-4} & 3.55 \T\B \\
    14.68 & \num[exponent-product=\ensuremath{\cdot}]{6.66e-06} & 4.86 & \num[exponent-product=\ensuremath{\cdot}]{9.37e-4} & 3.69 & \num[exponent-product=\ensuremath{\cdot}]{3.69e-07} & 3.75 & \num[exponent-product=\ensuremath{\cdot}]{1.18e-4} & 3.09 \T\B \\
    20.45 & \num[exponent-product=\ensuremath{\cdot}]{2.11e-06} & 3.47 & \num[exponent-product=\ensuremath{\cdot}]{3.98e-4} & 2.579 & \num[exponent-product=\ensuremath{\cdot}]{8.46e-08} & 4.44 & \num[exponent-product=\ensuremath{\cdot}]{4.04e-05} & 3.22 \T\B \\
    \end{tabular}
    \\[5pt]
\caption{Robustness test vs secondary consolidation: polynomial degree $\ell = 3$. The parameters are chosen as in Table~\ref{tab:params_convtest} and the secondary consolidation coefficient is set as $\delta_2 \lambda = \num[exponent-product=\ensuremath{\cdot}, print-unity-mantissa=false]{1e+06}$.}
\label{tab:rob_vs_secondaryconsolidation}
\end{table}

\begin{table}[H]
    \centering 
    \small
    \begin{tabular}{c  c  c  c  c  c  c c c}
    $1/h$ & $\|\mathbf{e}^{u}\|_{L^2}$  & $\text{roc}^{u}_{L^2}$ & $\|\mathbf{e}^{u}\|_{dG}$  & $\text{roc}^{u}_{dG}$& $\|e^{p}\|_{L^2}$  & $\text{roc}^{p}_{L^2}$ & $\|e^{p}\|_{dG}$  & $\text{roc}^{p}_{dG}$  \T\B\B \\
    \hline
    5.53 & \num[exponent-product=\ensuremath{\cdot}]{3.74e-4} & - & 0.02 & - & \num[exponent-product=\ensuremath{\cdot}]{0.01} & - & 1.06 & - \T\B \\
    7.62 & \num[exponent-product=\ensuremath{\cdot}]{1.31e-4} & 3.28 & \num[exponent-product=\ensuremath{\cdot}]{8.99e-3} & 1.90 & \num[exponent-product=\ensuremath{\cdot}]{3.00e-3} & 4.65 & 0.33 & 3.59 \T\B \\
    10.53 & \num[exponent-product=\ensuremath{\cdot}]{3.36e-05} & 4.21 & \num[exponent-product=\ensuremath{\cdot}]{3.20e-3} & 3.20 & \num[exponent-product=\ensuremath{\cdot}]{4.62e-4} & 5.92 & 0.07 & 4.71 \T\B \\
    14.68 & \num[exponent-product=\ensuremath{\cdot}]{6.68e-06} & 4.85 & \num[exponent-product=\ensuremath{\cdot}]{9.36e-4} & 3.69 & \num[exponent-product=\ensuremath{\cdot}]{1.07e-4} & 4.41 & 0.02 & 3.59 \T\B \\
    20.45 & \num[exponent-product=\ensuremath{\cdot}]{2.20e-06} & 3.35 & \num[exponent-product=\ensuremath{\cdot}]{3.98e-4} & 2.58 & \num[exponent-product=\ensuremath{\cdot}]{6.10e-05} & 1.68 & \num[exponent-product=\ensuremath{\cdot}]{5.00e-3} & 4.45 \T\B \\
    \end{tabular}
    \\[5pt]
\caption{Robustness test vs secondary consolidation: polynomial degree $\ell = 3$. The parameters are chosen as in Table~\ref{tab:params_convtest}, the secondary consolidation coefficient is set as $\delta_2 \lambda = \num[exponent-product=\ensuremath{\cdot}, print-unity-mantissa=false]{1e+06}$, and the storage coefficient is set as $d_0 = 10^{-6}$.}
\label{tab:rob_vs_secondaryconsolidation_d0}
\end{table}

For what concerns the $dG$-error analysis, in the case of small Lamé parameter $\mu$ and small viscoelasticity retardation time $\delta_1$ the errors decay as expected, while in the case of small conductivity $\mathbf{D}$, we observe a decrease of the $dG$-errors of the pressure when refining the mesh. However, we observe that their values are affected by the fact that the conductivity value enters the norm's definition. Then, we observe the performance of the method for high values of the secondary consolidation coefficient $\delta_2 \lambda$, considering both the cases in which we have not-degenerating (cf. Table~\ref{tab:rob_vs_secondaryconsolidation}) and degenerating storage coefficient $d_0$ (cf. Table~\ref{tab:rob_vs_secondaryconsolidation_d0}). The method is robust in both regimes; however -- as expected -- it behaves slightly worse when the storage coefficient is $\ll 1$. Finally, we observe that our method is robust also in the case of strong coupling between the two equations of the model problem, cf. Table~\ref{tab:rob_vs_couplingstrength}.

\begin{table}[H]
    \centering 
    \small
    \begin{tabular}{c  c  c  c  c  c  c c c}
    $1/h$ & $\|\mathbf{e}^{u}\|_{L^2}$  & $\text{roc}^{u}_{L^2}$ & $\|\mathbf{e}^{u}\|_{dG}$  & $\text{roc}^{u}_{dG}$& $\|e^{p}\|_{L^2}$  & $\text{roc}^{p}_{L^2}$ & $\|e^{p}\|_{dG}$  & $\text{roc}^{p}_{dG}$  \T\B\B \\
    \hline
    5.53 & \num[exponent-product=\ensuremath{\cdot}]{1.56e-05} & - & \num[exponent-product=\ensuremath{\cdot}]{4.14e-3} & - & \num[exponent-product=\ensuremath{\cdot}]{4.97e-05} & - & \num[exponent-product=\ensuremath{\cdot}]{4.57e-3} & - \T\B \\
    7.62 & \num[exponent-product=\ensuremath{\cdot}]{3.38e-06} & 4.76 & \num[exponent-product=\ensuremath{\cdot}]{1.40e-3} & 3.38 & \num[exponent-product=\ensuremath{\cdot}]{1.03e-05} & 4.92 & \num[exponent-product=\ensuremath{\cdot}]{1.49e-3} & 3.50 \T\B \\
    10.53 & \num[exponent-product=\ensuremath{\cdot}]{6.46e-07} & 5.12 & \num[exponent-product=\ensuremath{\cdot}]{4.38e-4} & 3.60 & \num[exponent-product=\ensuremath{\cdot}]{2.10e-06} & 4.91 & \num[exponent-product=\ensuremath{\cdot}]{4.33e-4} & 3.82 \T\B \\
    14.68 & \num[exponent-product=\ensuremath{\cdot}]{1.79e-07} & 3.86 & \num[exponent-product=\ensuremath{\cdot}]{1.73e-4} & 2.79 & \num[exponent-product=\ensuremath{\cdot}]{5.74e-07} & 3.90 & \num[exponent-product=\ensuremath{\cdot}]{1.57e-4} & 3.04 \T\B \\
    20.45 & \num[exponent-product=\ensuremath{\cdot}]{4.30e-08} & 4.30 & \num[exponent-product=\ensuremath{\cdot}]{5.58e-05} & 3.41 & \num[exponent-product=\ensuremath{\cdot}]{1.39e-07} & 4.28 & \num[exponent-product=\ensuremath{\cdot}]{5.35e-05} & 3.26 \T\B \\
    \end{tabular}
    \\[5pt]
\caption{Robustness test vs coupling strength: polynomial degree $\ell = 3$. The parameters are chosen as in Table~\ref{tab:params_convtest}, the first Lamè parameter has been set as $\lambda = \num[exponent-product=\ensuremath{\cdot}, print-unity-mantissa=false]{1e-06}$.}
\label{tab:rob_vs_couplingstrength}
\end{table}

\begin{remark}
In Theorem~\ref{thm:stabdisc}, we observe that the stability estimate is not robust with respect to $\mathbf{D}$ and the product between $\mu$ and $\delta_1$. However, in Table~\ref{tab:rob_vs_conductivity}, Table~\ref{tab:rob_vs_mu}, and Table~\ref{tab:rob_vs_delta1}, we have shown that the proposed method is robust with respect to low values of these parameters. 
\end{remark}

%% file: Sections/Carcione2018.tex
This section considers a wave propagation problem in thermoelastic media inspired by \cite{Carcione2018}. The aim of this simulation is to prove that the proposed scheme can reproduce known results present in the literature for the thermoelastic framework and can give \textit{physically-sound} results. 

\bigskip
\noindent
\textit{Vertical source term.} We consider a domain $\Omega = (0\,\mathrm{m}, 2310\, \mathrm{m})^2$ with the thermoelastic properties reported in Table~\ref{tab:Carcione2018_test1} \cite{Carcione2018}.
\begin{table}[ht]
    \centering 
    \begin{tabular}{ l | l  c l | l }
    \textbf{Coefficient} & \textbf{Value} & & \textbf{Coefficient} & \textbf{Value} \T\B \\
    $\rho \ [\si{\kilogram \per \meter\cubed}]$ & 2650 & & $\beta \ [\si{\pascal\per\kelvin}]$ & 79200 \T\B \\ 
    $\mu \ [\si{\pascal}]$ & \num[exponent-product=\ensuremath{\cdot}]{6e+9} & & $\lambda \ [\si{\pascal}]$ & \num[exponent-product=\ensuremath{\cdot}]{4e+9} \T\B \\ 
    $a_0 \ [\si{\pascal \per \kelvin \squared}]$ & 117 & & $\boldsymbol{\Theta} \ [\si{\metre\squared \pascal \per \kelvin\squared \per \second}]$ & 10.5 $\mathbf{I}$  \T\B \\
    $\tau_1 [\si{\second}]$ & \num[exponent-product=\ensuremath{\cdot}]{1.49e-8} & & $\tau_2 [\si{\second}]$ & \num[exponent-product=\ensuremath{\cdot}]{1.49e-8} \T\B \\
    $\delta_1 [\si{\second}]$ & 0 & & $\delta_2 [\si{\second}]$ & 0 \T\B \\
    \end{tabular}
    \\[10pt]
    \caption{Wave propagation in thermoelastic media: homogeneous medium properties}
    \label{tab:Carcione2018_test1}
\end{table}
\par
In the first test case, we consider a vertical source term $\mathbf{f}$ in \eqref{eq:momentum_cons_1}. The source is set in the computational domain's center and multiplied by a time-history function $h(t)$. In our case, the time evolution is given by \cite{Carcione2018} $h(t) = A_0 \cos\left[ 2 \pi (t - t_0) f_0 \right] \exp\left[ -2 (t - t_0)^2 f_0^2 \right]$, where $A_0 = \num[exponent-product=\ensuremath{\cdot}, print-unity-mantissa=false]{1e+4}$ \si{\meter} is the amplitude, $f_0 = 5$ \si{\hertz} is the peak-frequency, and $t_0 = 3/(2 f_0) = 0.3$ \si{\second} is the time-shift. We adopt a polygonal mesh with mesh size $h \sim 76$ \si{\meter} (3500 elements) and polynomial degree $\ell = 4$. As a time-stepping scheme we employ the Newmark-$\beta$ scheme, with $\Delta t = \num[exponent-product=\ensuremath{\cdot}, print-unity-mantissa=false]{5e-4}$ \si{\second} and $T_f = 0.5$ \si{\second}. Finally, we complete our problem with homogeneous Dirichlet boundary conditions and null initial conditions.
In the following, we denote by $\mathbf{v}_h$ the solid velocity (i.e. $\dot{ \mathbf{u}}_h$) and by $v_{h,y}$ its vertical component.

\begin{figure}[ht]
\begin{subfigure}[b]{.33\textwidth}
    \centering
    \includegraphics[width=1\textwidth]{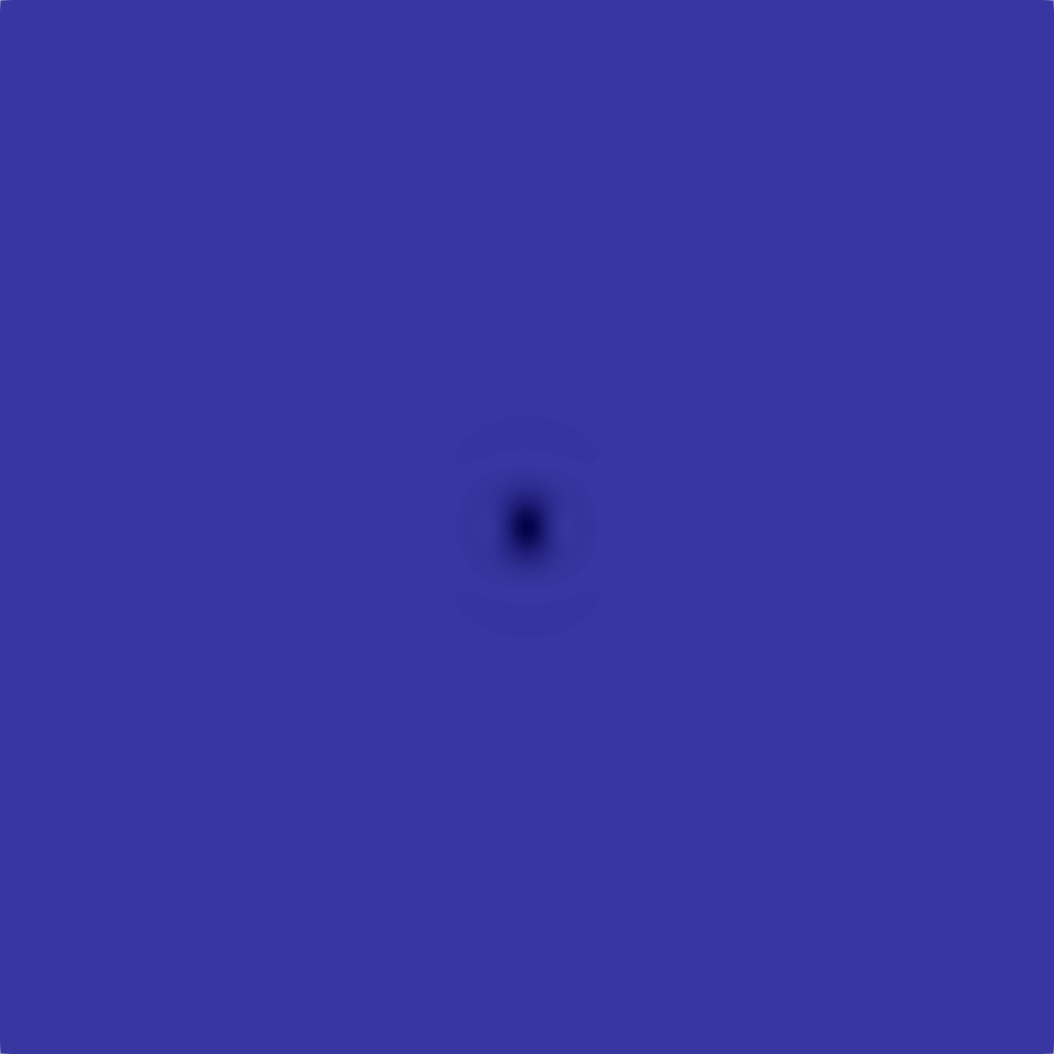}
    \label{fig:Carcione2018_v_01_fy}
\end{subfigure}
\begin{subfigure}[b]{.33\textwidth}
    \centering
    \includegraphics[width=1\textwidth]{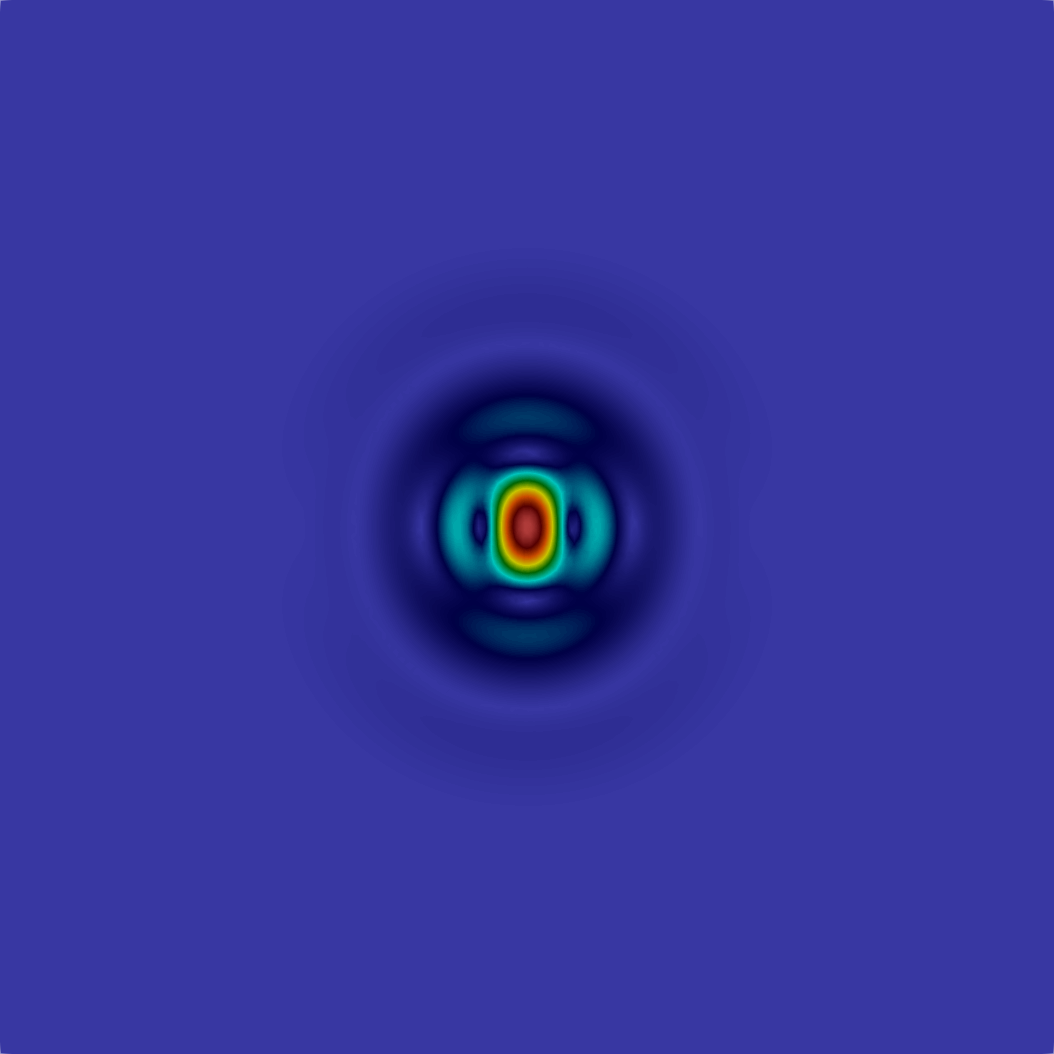}
    \label{fig:Carcione2018_v_03_fy}
\end{subfigure}
\begin{subfigure}[b]{.33\textwidth}
    \centering
    \includegraphics[width=1\textwidth]{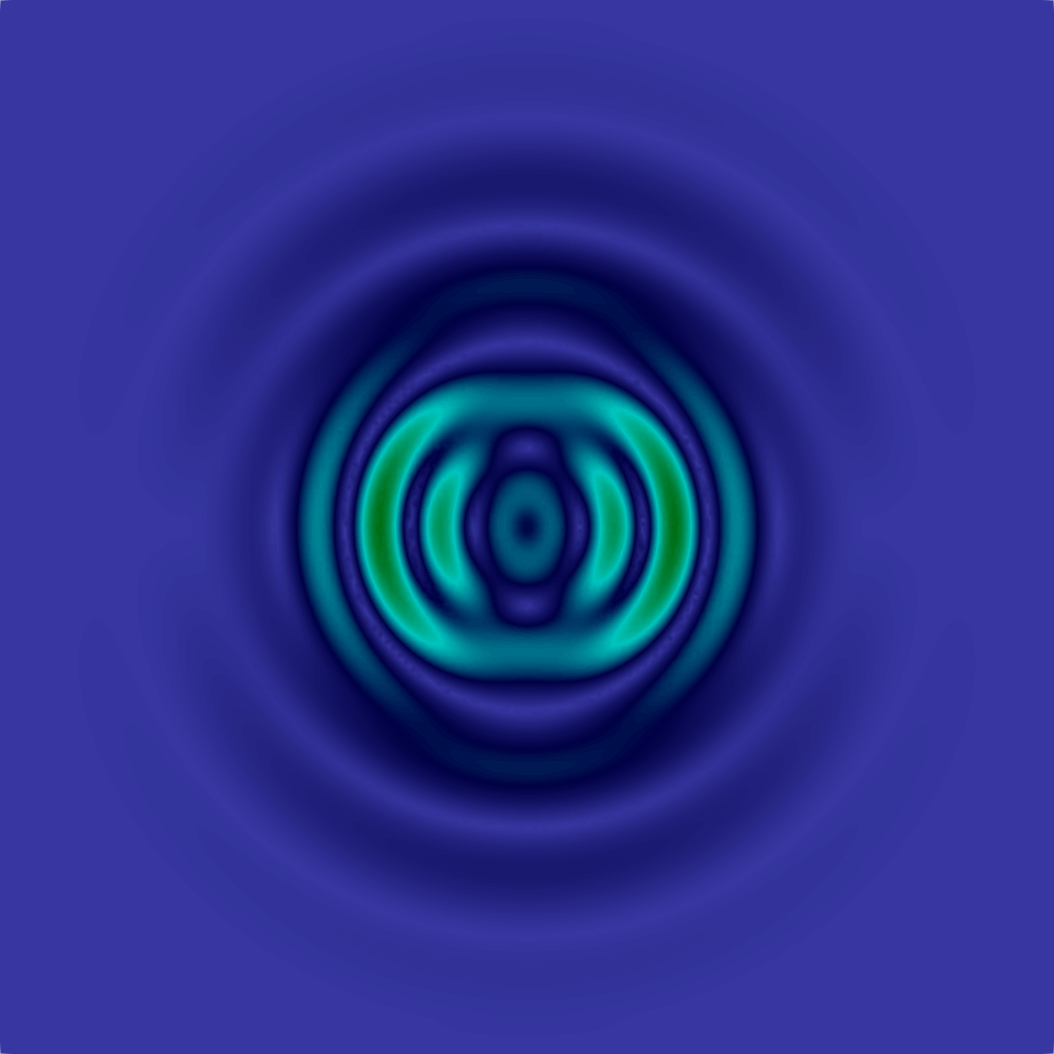}
    \label{fig:Carcione2018_v_05_fy}
\end{subfigure}

\vspace{-0.4cm}
\centering
\includegraphics[width=0.3\textwidth]{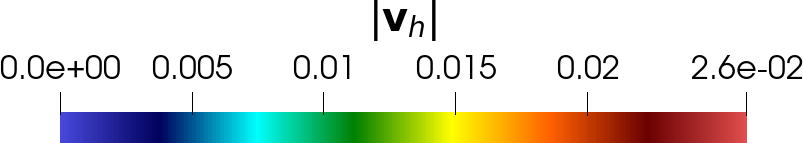}

\caption{Wave propagation in thermoelastic media with vertical source term: computed velocity field $|\mathbf{v}_h|$ at the time instants $t=0.1 \si{\second}$ (left), $t=0.3 \si{\second}$ (center), $t=0.5 \si{\second}$ (right).}
\label{fig:Carcione2018_v}
\end{figure}

\begin{figure}[ht]
\begin{subfigure}[b]{.33\textwidth}
    \centering
    \includegraphics[width=1\textwidth]{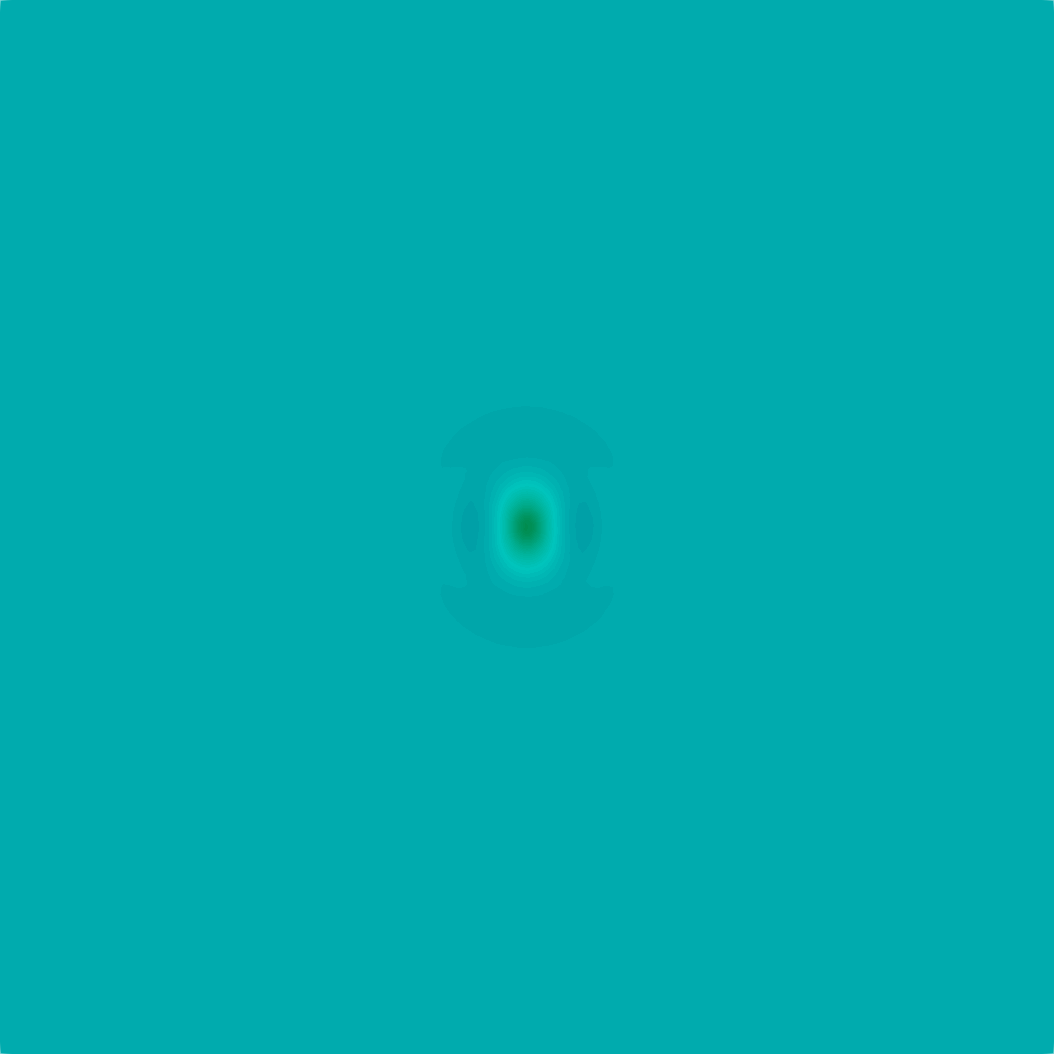}
    \label{fig:Carcione2018_vy_01_fy}
\end{subfigure}
\begin{subfigure}[b]{.33\textwidth}
    \centering
    \includegraphics[width=1\textwidth]{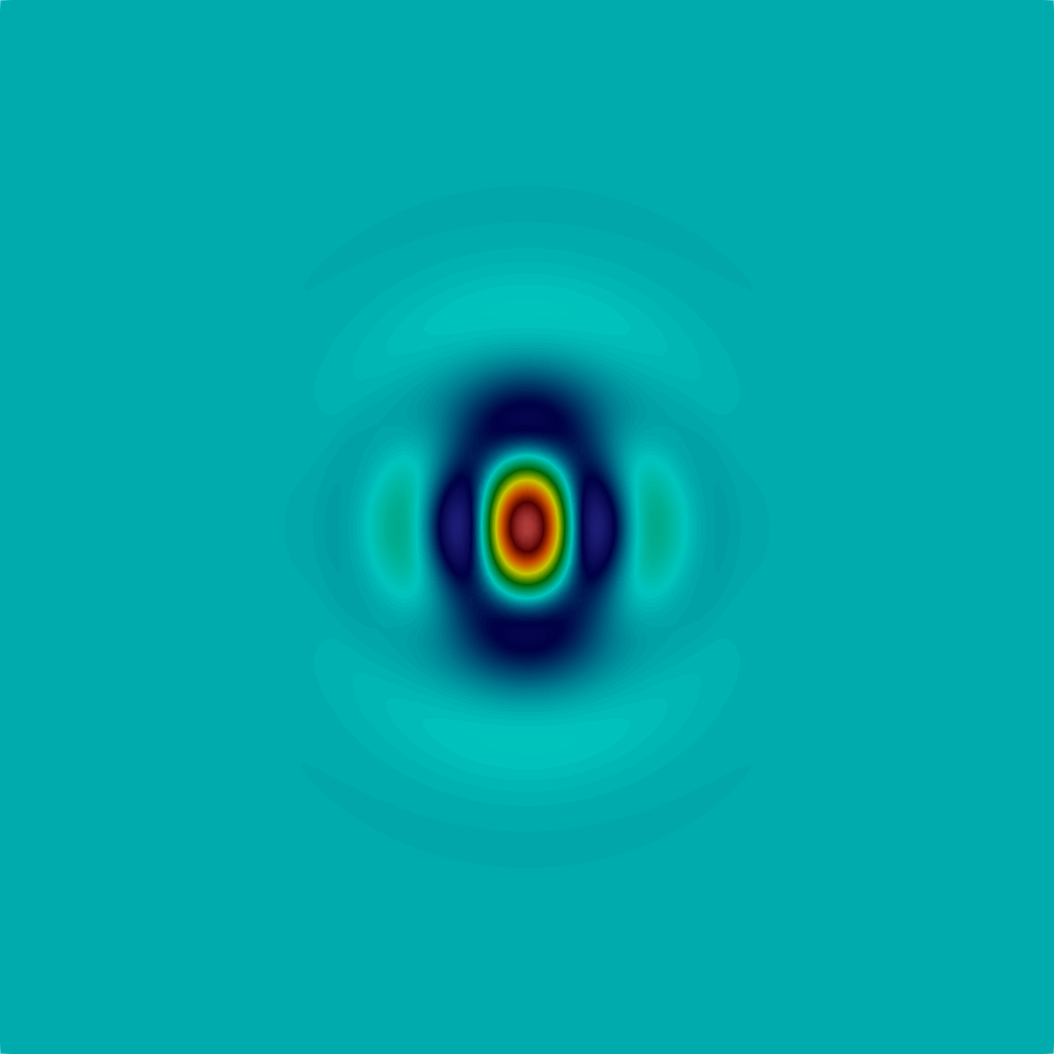}
    \label{fig:Carcione2018_vy_03_fy}
\end{subfigure}
\begin{subfigure}[b]{.33\textwidth}
    \centering
    \includegraphics[width=1\textwidth]{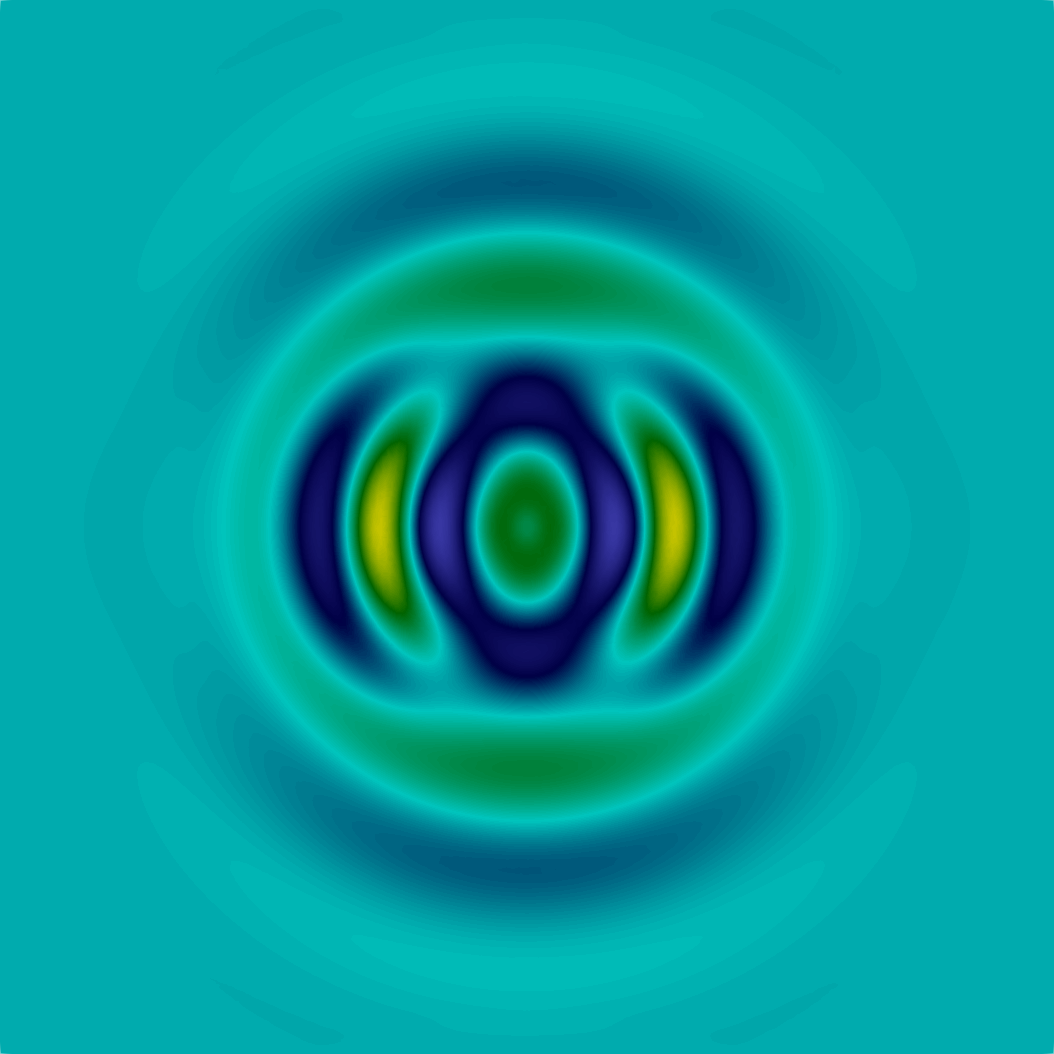}
    \label{fig:Carcione2018_vy_05_fy}
\end{subfigure}

\vspace{-0.4cm}
\centering
\includegraphics[width=0.3\textwidth]{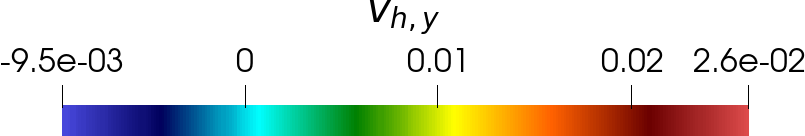}

\caption{Wave propagation in thermoelastic media with vertical source term: computed vertical component of the velocity field $v_{h,y}$ at the time instants $t=0.1 \si{\second}$ (left), $t=0.3 \si{\second}$ (center), $t=0.5 \si{\second}$ (right).}
\label{fig:Carcione2018_vy}
\end{figure}

\begin{figure}[ht]
\begin{subfigure}[b]{.33\textwidth}
    \centering
    \includegraphics[width=1\textwidth]{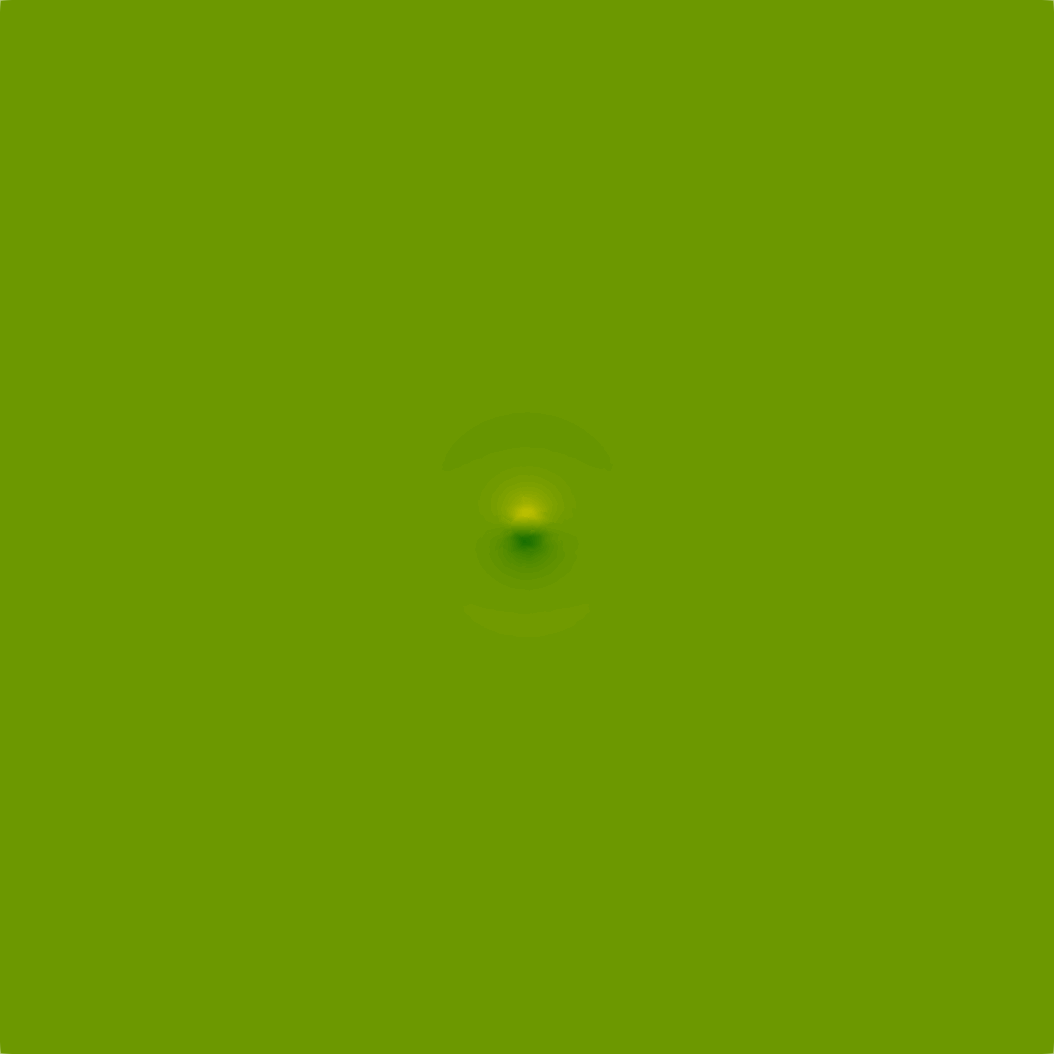}
    \label{fig:Carcione2018_T_01_fy}
\end{subfigure}
\begin{subfigure}[b]{.33\textwidth}
    \centering
    \includegraphics[width=1\textwidth]{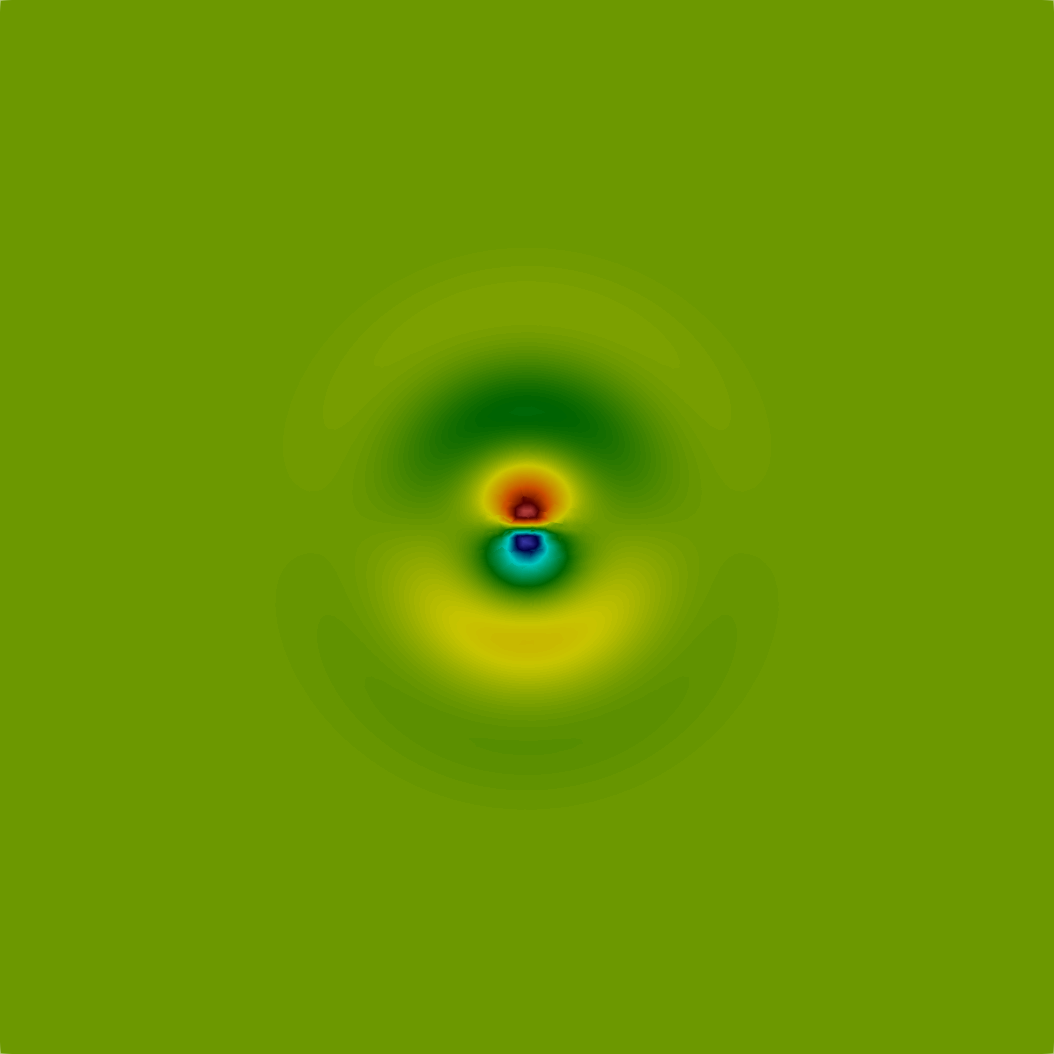}
    \label{fig:Carcione2018_T_03_fy}
\end{subfigure}
\begin{subfigure}[b]{.33\textwidth}
    \centering
    \includegraphics[width=1\textwidth]{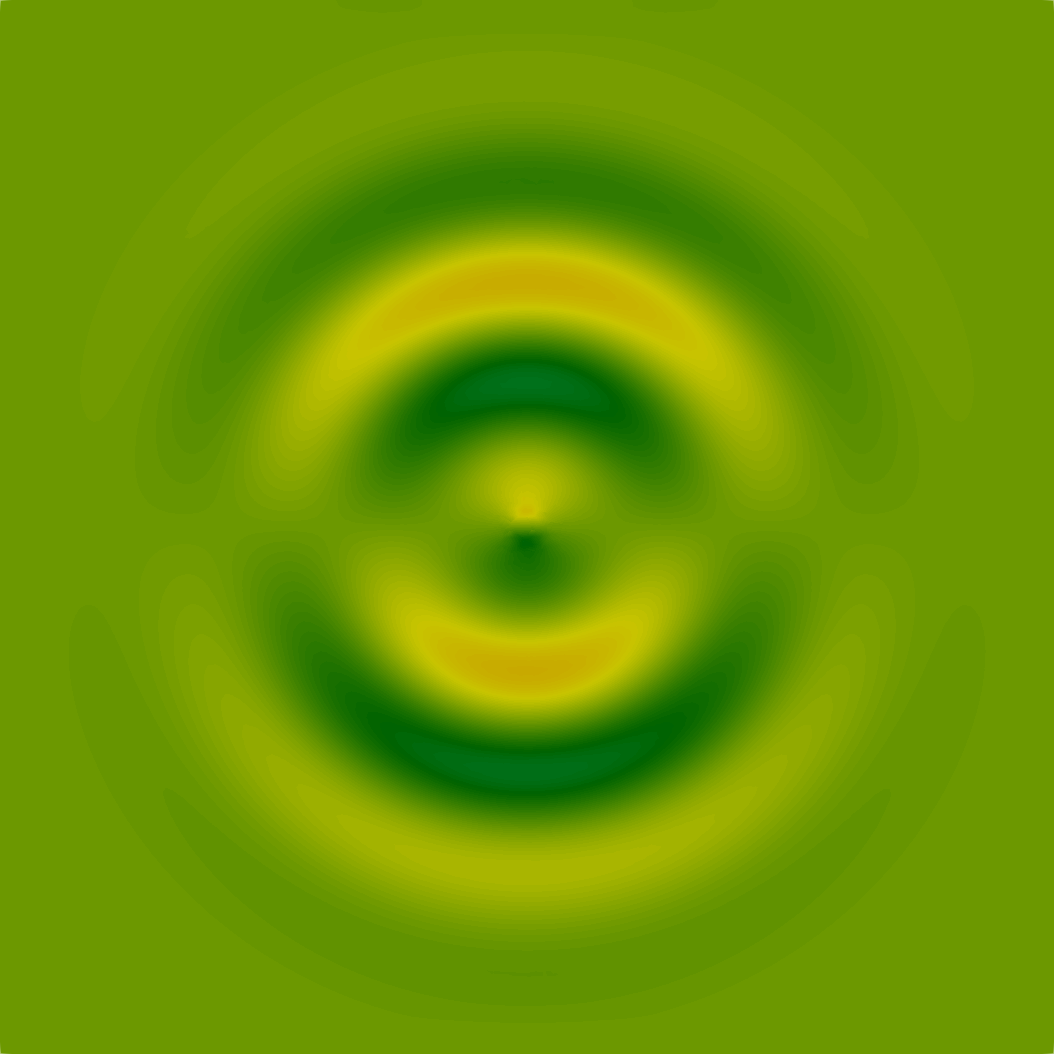}
    \label{fig:Carcione2018_T_05_fy}
\end{subfigure}

\vspace{-0.4cm}
\centering
\includegraphics[width=0.3\textwidth]{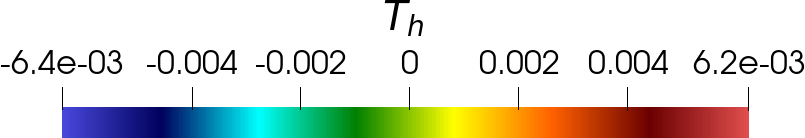}

\caption{Wave propagation in thermoelastic media with vertical source term: computed temperature field $T_h$ at the time instants $t=0.1 \si{\second}$ (left), $t=0.3 \si{\second}$ (center), $t=0.5 \si{\second}$ (right).}
\label{fig:Carcione2018_T}
\end{figure}

From the results of Figure~\ref{fig:Carcione2018_v} we observe a
symmetric wavefront that detaches from the center of the domain symmetric to the $x$- and $y$-axes. This behavior is correct due to the form of the forcing term we are imposing. Looking at the three snapshots, we can observe the presence of the elastic $E$-wave captured by our scheme and the $S$-wave. The presence of these two waves is more evident by looking at Figure~\ref{fig:Carcione2018_vy}, where the computed vertical velocity field is reported. Indeed, it is possible to observe the propagation of the $ E$ wave, traveling faster and along the $y$-direction, and the propagation of the $S$-waves that are slower and propagate along the $x$-direction. Last, in Figure~\ref{fig:Carcione2018_T}, we observe the propagation of the diffusive $T$-wave, that is originated by a vertical source term in the momentum conservation equation. A thermal source term in the energy conservation equation is not considered here. In conclusion, considering the different central peak frequencies $f_0$, we can see a good agreement between our results and those presented in \cite{Carcione2019}.

\bigskip
\noindent
\textit{Shear source term.} In the second test case, we model the forcing term of the momentum equation $\mathbf{f}$ as a shear source. The forcing term is modeled as $\mathbf{f} = - \mathbf{M} \, \div{\delta}(\mathbf{x} - \mathbf{x}_s) \, h(t)$ \cite{Morency2008}, where $\mathbf{M}$ is a moment tensodr, $\mathbf{x}_s$ is the point-source location, $\delta(\mathbf{x} - \mathbf{x}_s)$ is the Kronecker delta located in $\mathbf{x}_s$, and $h(t)$ is the aforementioned time-history function. We consider $\mathbf{M}$ a tensor with zero components on the diagonal and not zero components off-diagonal; this choice induces the presence of shear waves and generally has a strong connection with the wave patterns we observe. A momentum source of this shape is often used in the context of earthquakes. The properties of the medium and the discretization parameters are chosen as in the vertical source term test.

\begin{figure}[ht]
\begin{subfigure}[b]{.33\textwidth}
    \centering
    \includegraphics[width=1\textwidth]{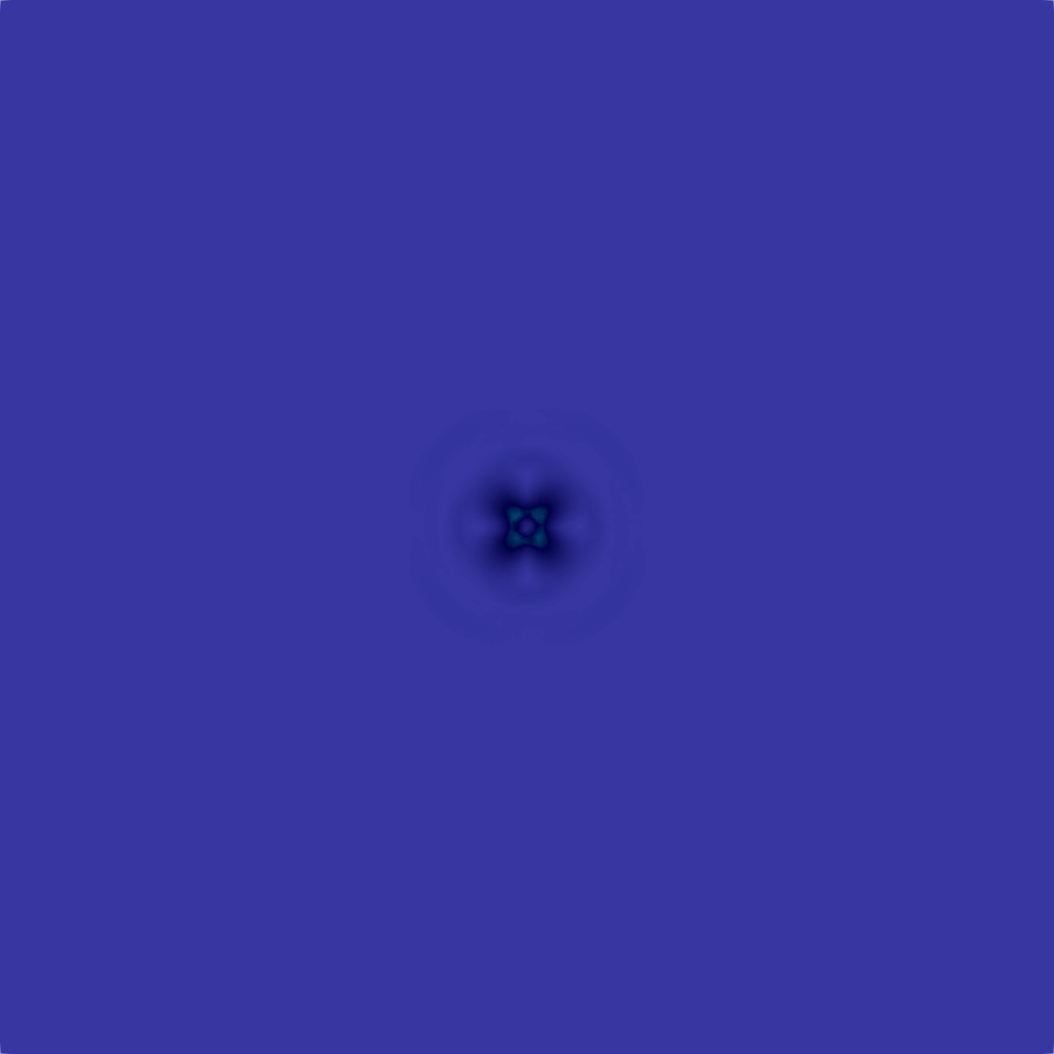}
    \label{fig:Carcione2018_v_01_fs}
\end{subfigure}
\begin{subfigure}[b]{.33\textwidth}
    \centering
    \includegraphics[width=1\textwidth]{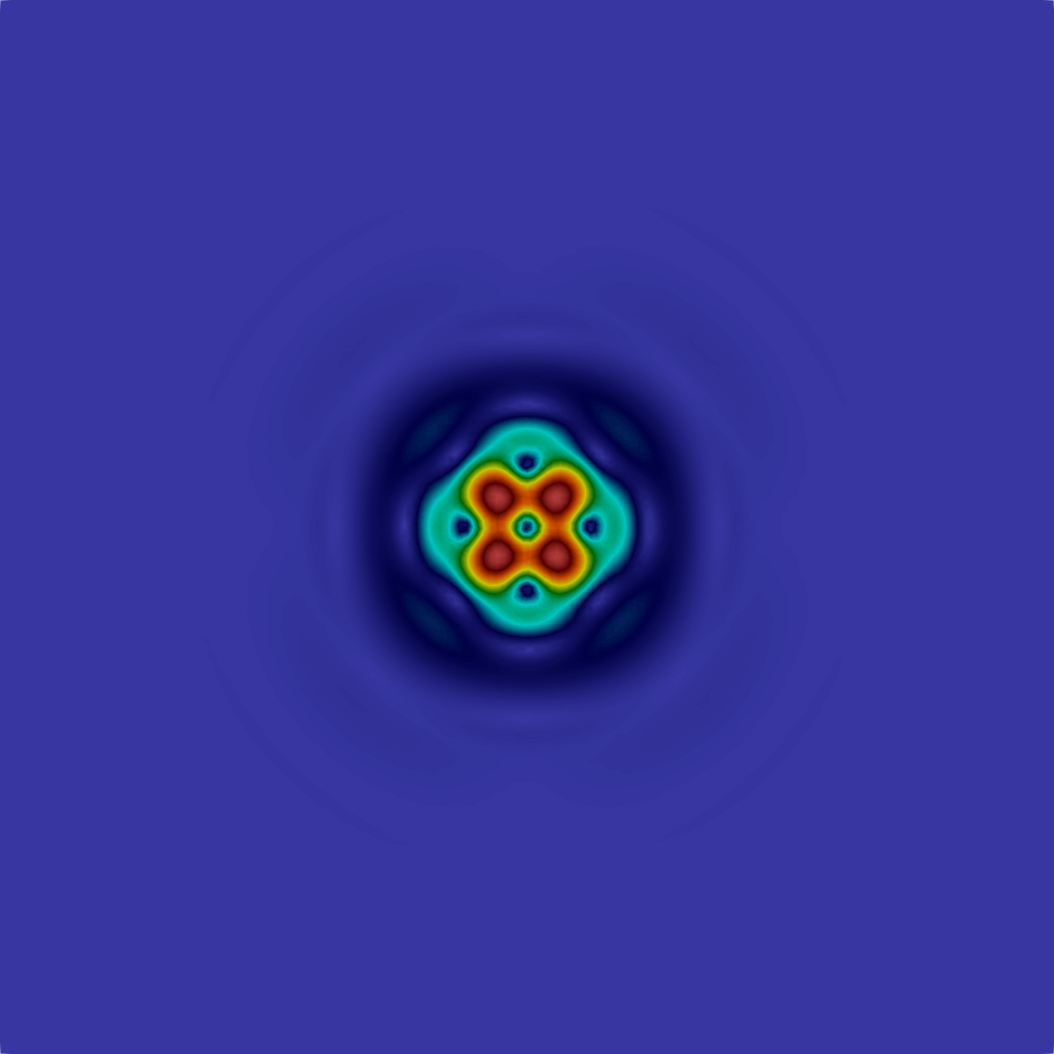}
    \label{fig:Carcione2018_v_03_fs}
\end{subfigure}
\begin{subfigure}[b]{.33\textwidth}
    \centering
    \includegraphics[width=1\textwidth]{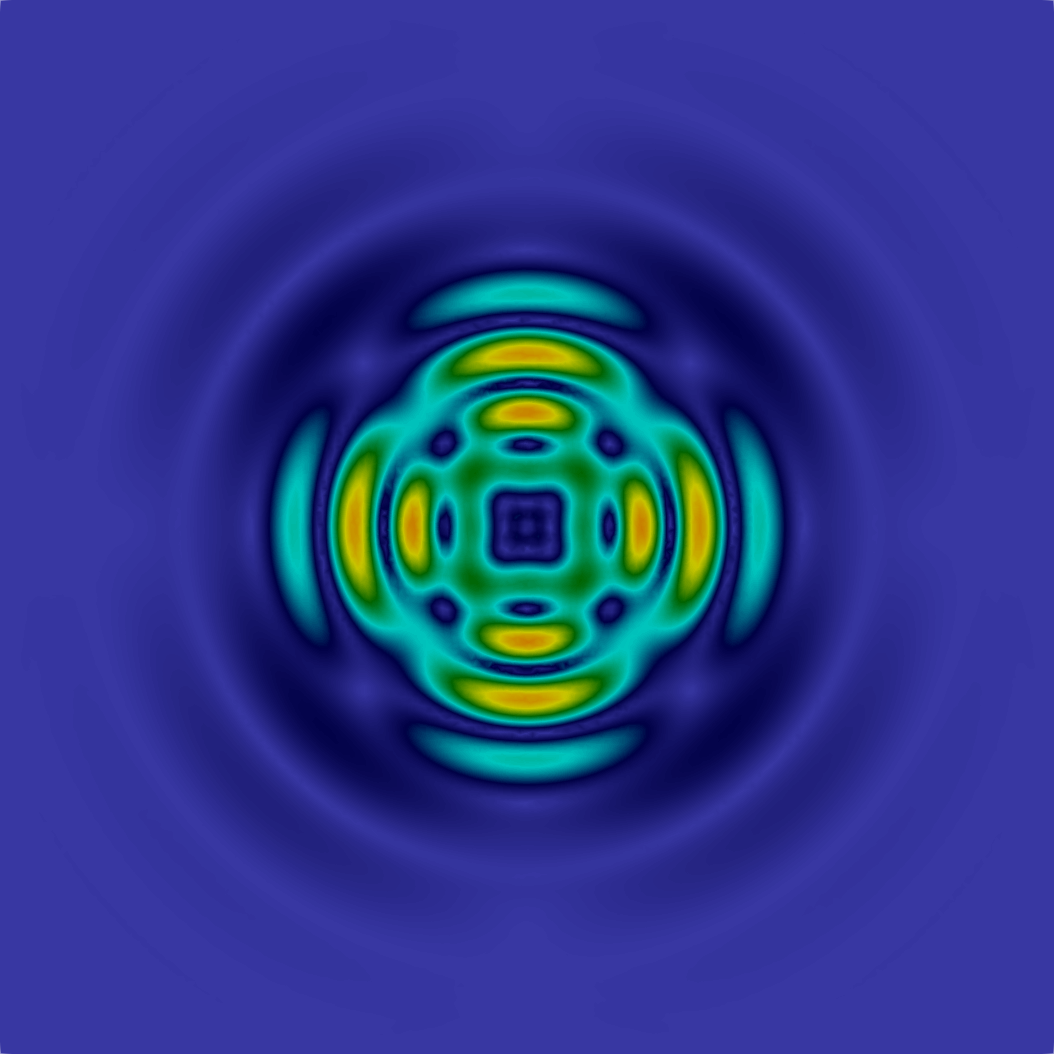}
    \label{fig:Carcione2018_v_05_fs}
\end{subfigure}

\vspace{-0.4cm}
\centering
\includegraphics[width=0.3\textwidth]{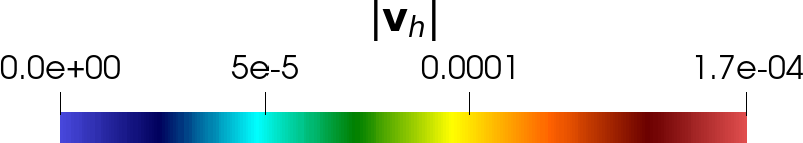}

\caption{Wave propagation in thermoelastic media with shear source term: computed velocity field $|\mathbf{v}_h|$ at the time instants $t=0.1 \si{\second}$ (left), $t=0.3 \si{\second}$ (center), $t=0.5 \si{\second}$ (right)}
\label{fig:Carcione2018_v_fs}
\end{figure}

\begin{figure}[ht]
\begin{subfigure}[b]{.33\textwidth}
    \centering
    \includegraphics[width=1\textwidth]{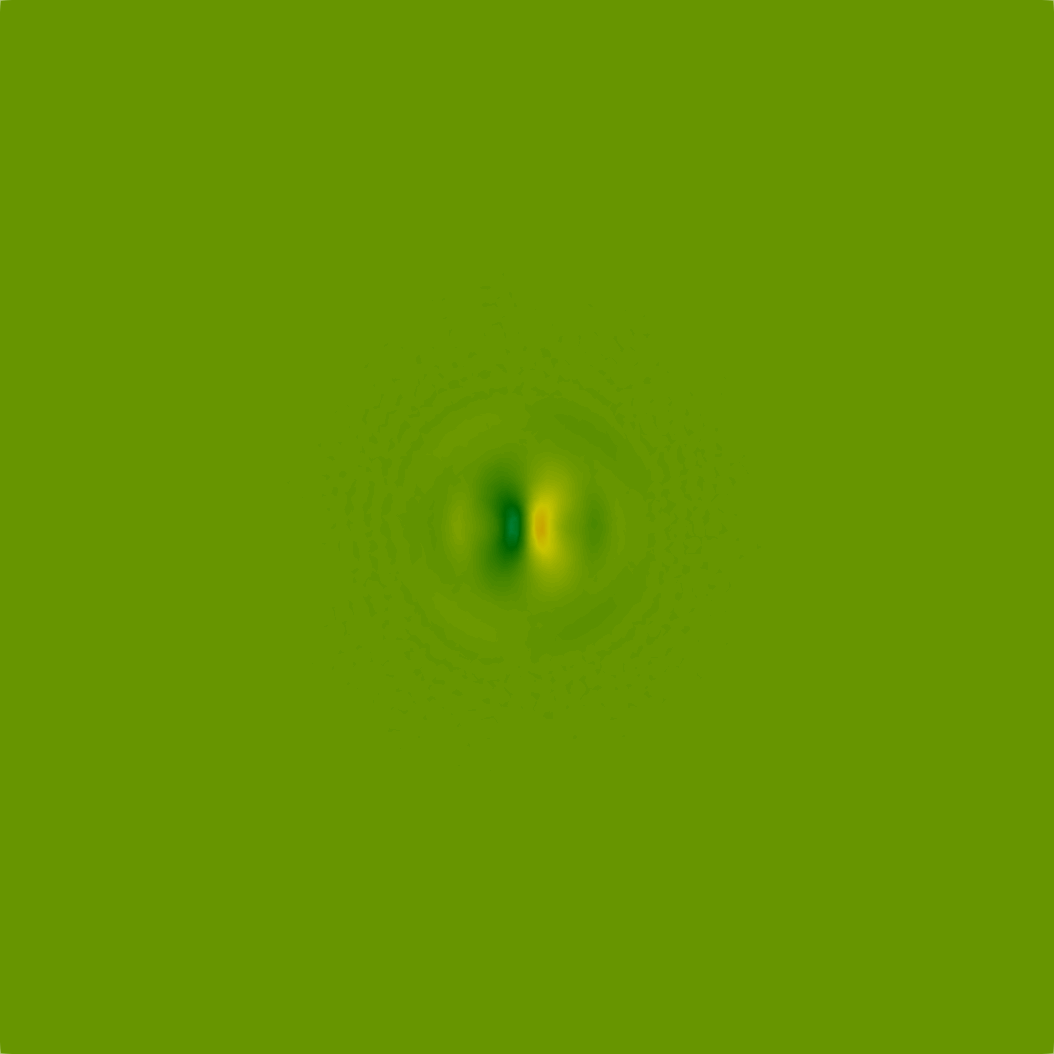}
    \label{fig:Carcione2018_vy_01_fs}
\end{subfigure}
\begin{subfigure}[b]{.33\textwidth}
    \centering
    \includegraphics[width=1\textwidth]{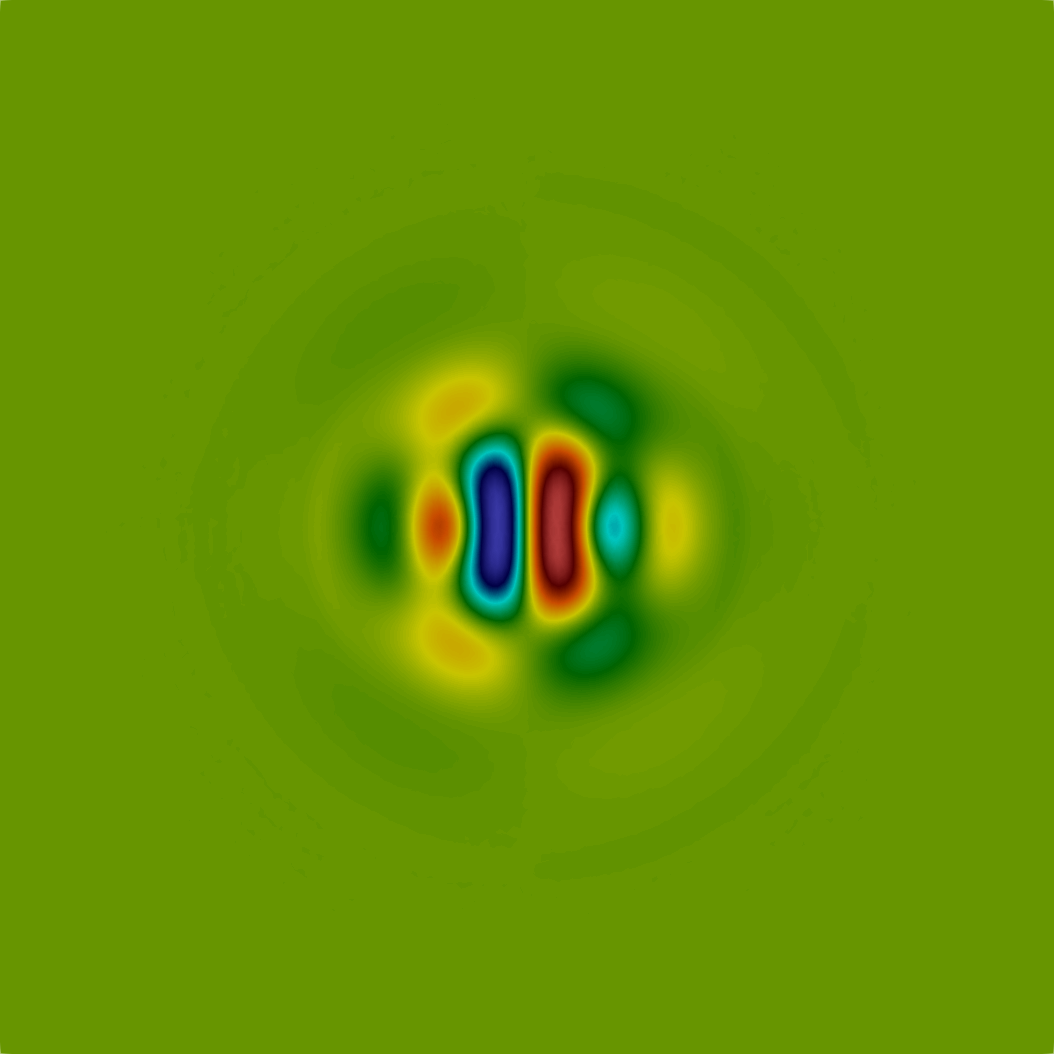}
    \label{fig:Carcione2018_vy_03_fs}
\end{subfigure}
\begin{subfigure}[b]{.33\textwidth}
    \centering
    \includegraphics[width=1\textwidth]{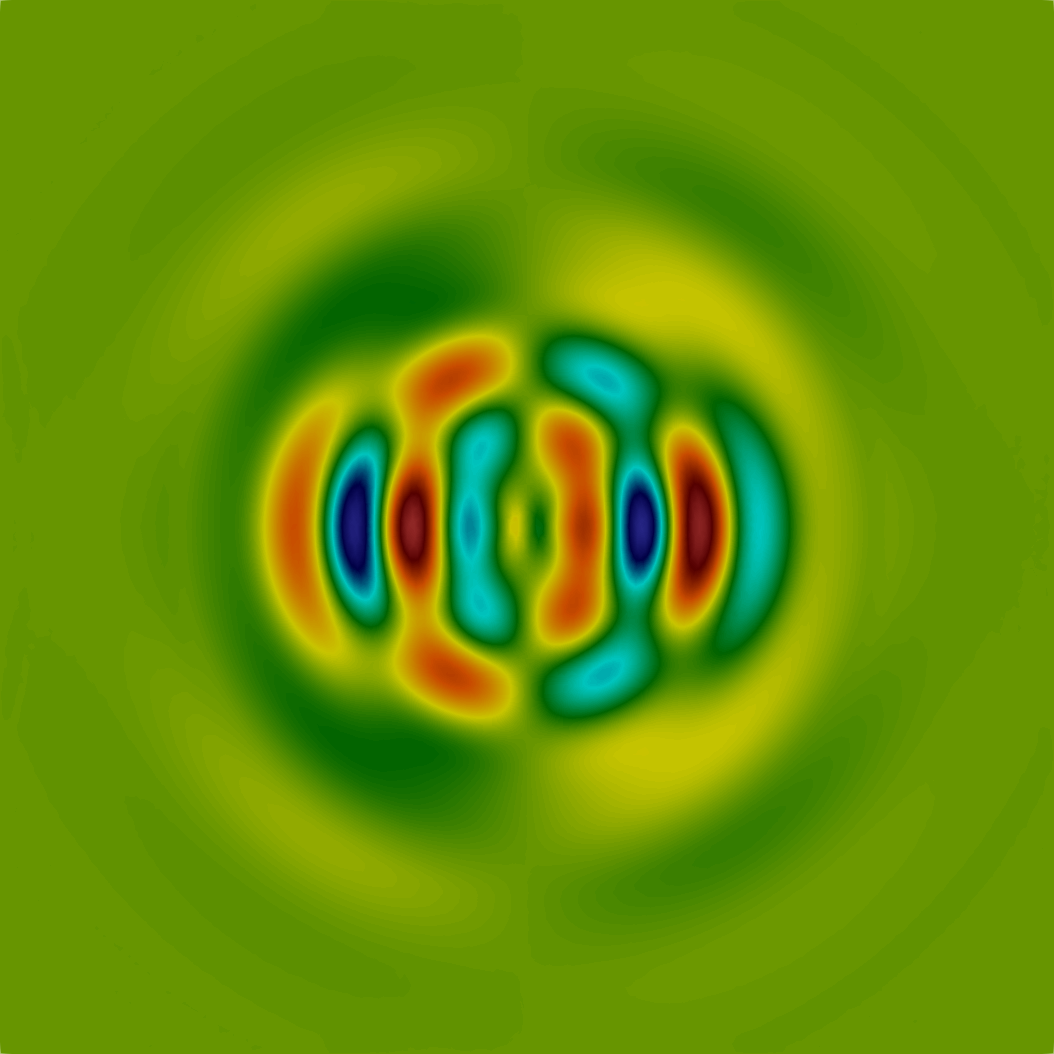}
    \label{fig:Carcione2018_vy_05_fs}
\end{subfigure}

\vspace{-0.4cm}
\centering
\includegraphics[width=0.3\textwidth]{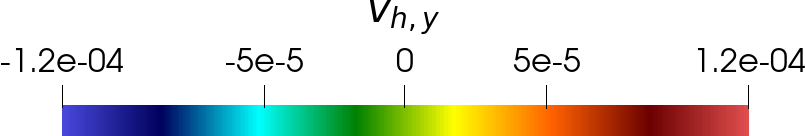}

\caption{Wave propagation in thermoelastic media with shear source term: computed vertical component of the velocity field $v_{h,y}$ at the time instants $t=0.1 \si{\second}$ (left), $t=0.3 \si{\second}$ (center), $t=0.5 \si{\second}$ (right).}
\label{fig:Carcione2018_vy_fs}
\end{figure}

\begin{figure}[ht]
\begin{subfigure}[b]{.33\textwidth}
    \centering
    \includegraphics[width=1\textwidth]{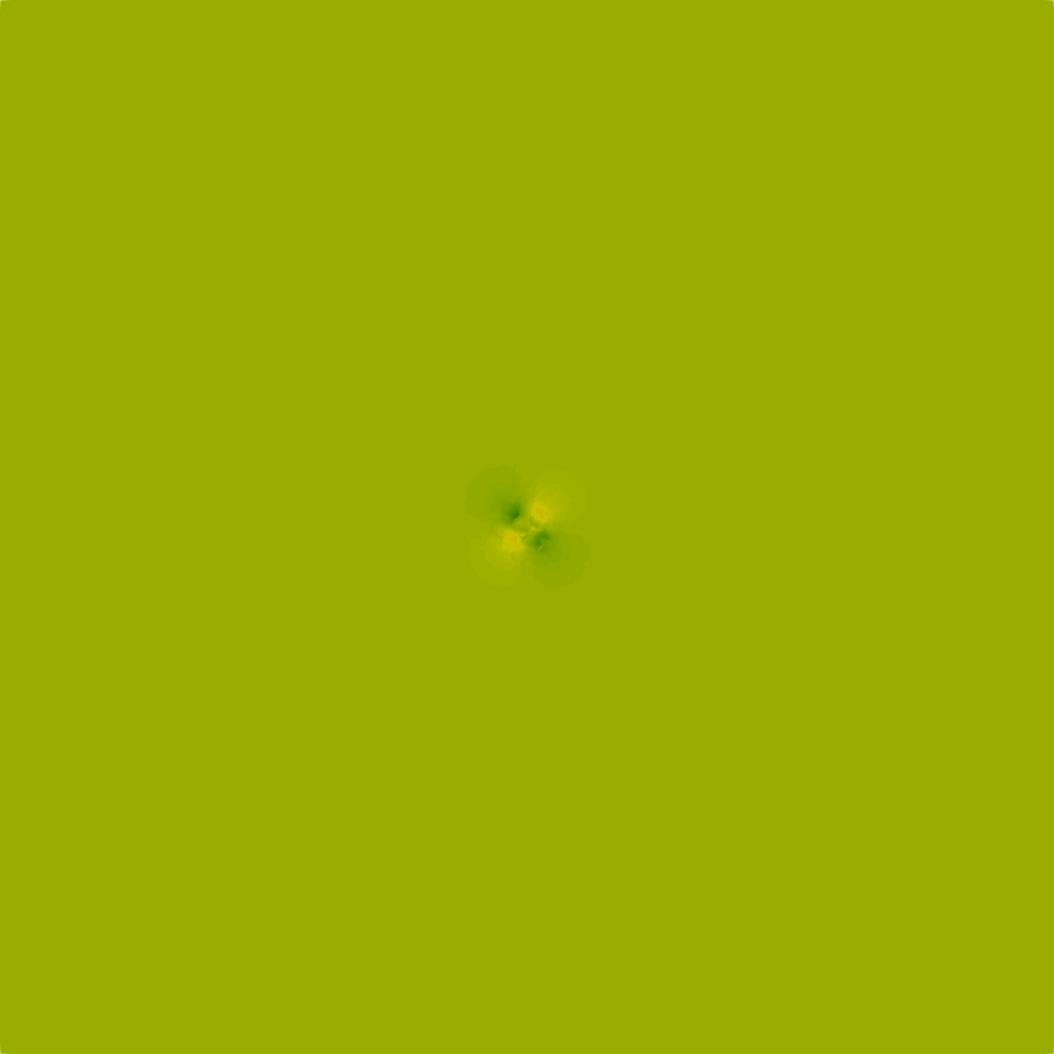}
    \label{fig:Carcione2018_T_01_fs}
\end{subfigure}
\begin{subfigure}[b]{.33\textwidth}
    \centering
    \includegraphics[width=1\textwidth]{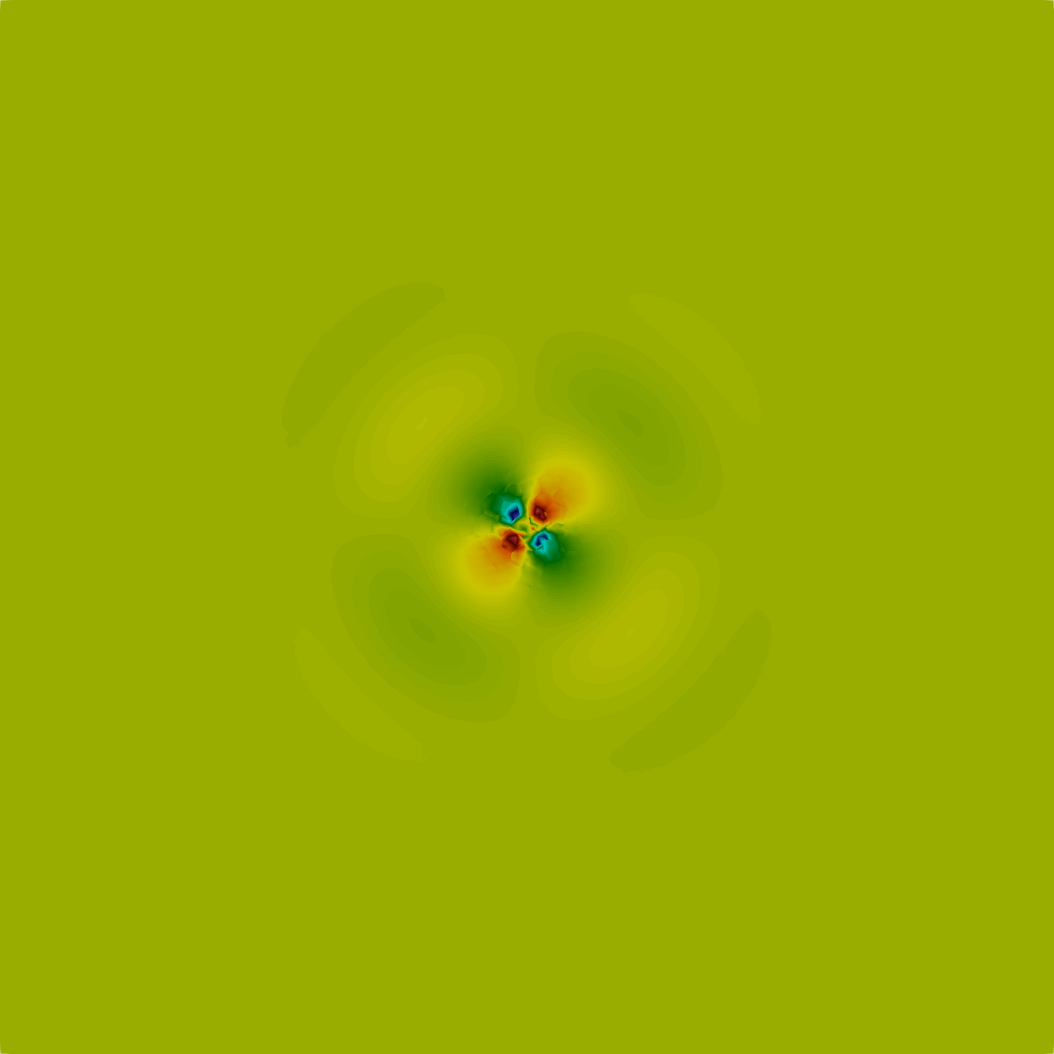}
    \label{fig:Carcione2018_T_03_fs}
\end{subfigure}
\begin{subfigure}[b]{.33\textwidth}
    \centering
    \includegraphics[width=1\textwidth]{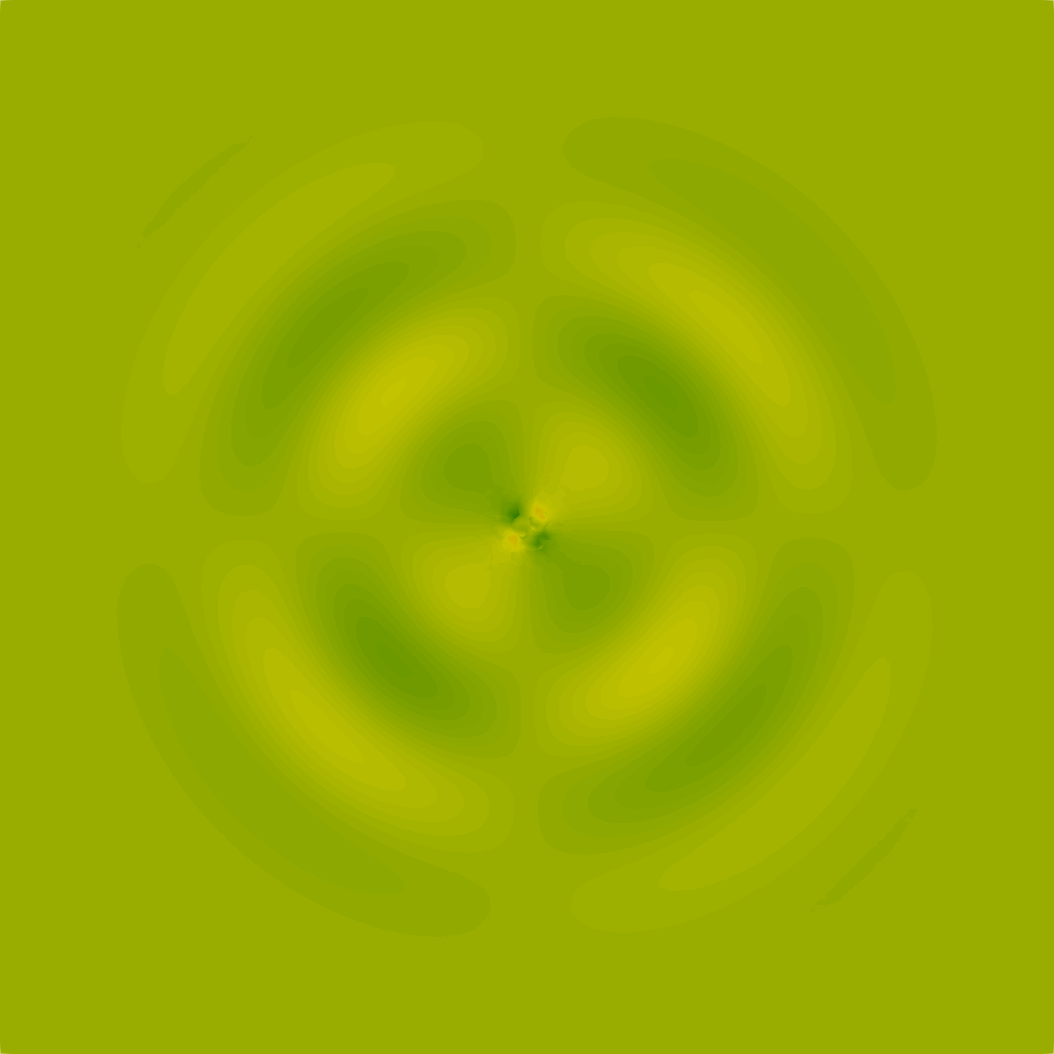}
    \label{fig:Carcione2018_T_05_fs}
\end{subfigure}

\vspace{-0.4cm}
\centering
\includegraphics[width=0.3\textwidth]{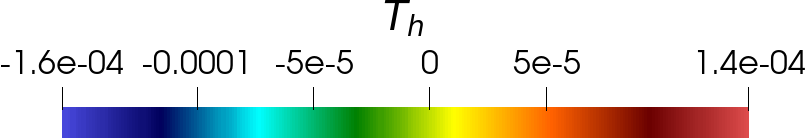}

\caption{Wave propagation in thermoelastic media with shear source term: computed temperature field $T_h$ at the time instants $t=0.1 \si{\second}$ (left), $t=0.3 \si{\second}$ (center), $t=0.5 \si{\second}$ (right).}
\label{fig:Carcione2018_T_fs}
\end{figure}

From the results of Figure~\ref{fig:Carcione2018_v_fs} we notice a 
symmetric wavefront that detaches from the center of the domain. The velocity field has a symmetric pattern with respect to the diagonals of the domain (due to the choice of the forcing term). Looking at the three snapshots, we observe the propagation of the elastic wave. From the results in Figure~\ref{fig:Carcione2018_vy_fs}, we notice the presence of shear waves and the anti-symmetric pattern of the wave fronts with respect to the $y$-axis. In Figure~\ref{fig:Carcione2018_T_fs}, we can see the presence of the diffusive thermal $T$-wave, which has a symmetric pattern with respect to the diagonals of the domain and antisymmetric with respect to the $x$- and $y$-axis. The results we obtain for this test are coherent with \cite{Carcione2018} (thermoelasticity) and \cite{Carcione2019,Bonetti2023} (thermo-poroelasticity).

\bigskip
\noindent
\textit{Shear source term in a heterogeneous media.} In this last thermoelastic test case, we model the thermoelastic wave propagation in a heterogeneous media. We split the domain $\Omega = (-1155\,\mathrm{m}, 1155\,\mathrm{m}) \times (0\,\mathrm{m}, 2310\,\mathrm{m})$ into two vertical layers. The left layer is characterized by the same thermoelastic properties of the test case \textit{Shear source term}, while in the right layer, we consider the following: cf. Table~\ref{tab:Carcione2018_test2} (the parameters that are not listed there are taken as in Table~\ref{tab:Carcione2018_test1}). This final test case aims to investigate how the heterogeneity of the media can affect the wave propagation phenomena. The forcing terms, and discretizations parameters are the same of \textit{Shear source term}-test case.
\begin{table}[H]
       \centering 
    \begin{tabular}{ l | l  c l | l}
    \textbf{Coefficient} & \textbf{Value} & & \textbf{Coefficient} & \textbf{Value} \T\B \\
    $\mu \ [\si{\pascal}]$ & \num[exponent-product=\ensuremath{\cdot},  print-unity-mantissa=false]{1e+9} & & $\lambda \ [\si{\pascal}]$ & \num[exponent-product=\ensuremath{\cdot}]{6.5e+8} \T\B \\ 
    $a_0 \ [\si{\pascal \per \kelvin \squared}]$ & 20 & & $\boldsymbol{\Theta} \ [\si{\metre\squared \pascal \per \kelvin\squared \per \second}]$ & 5 $\mathbf{I}$  \T\B \\
    $\tau_1 [\si{\second}]$ & \num[exponent-product=\ensuremath{\cdot}]{2.5e-7} & & $\tau_2 [\si{\second}]$ & \num[exponent-product=\ensuremath{\cdot}]{2.5e-7} \T\B \\
    \end{tabular}
    \\[10pt]
    \caption{Wave propagation in thermoelastic media: heterogeneous media properties (right layer). The thermoelastic properties of the left layer are reported in Table~\ref{tab:Carcione2018_test1}.}
    \label{tab:Carcione2018_test2}
\end{table}
\begin{figure}[ht]
\begin{subfigure}[b]{.33\textwidth}
    \centering
    \includegraphics[width=1\textwidth]{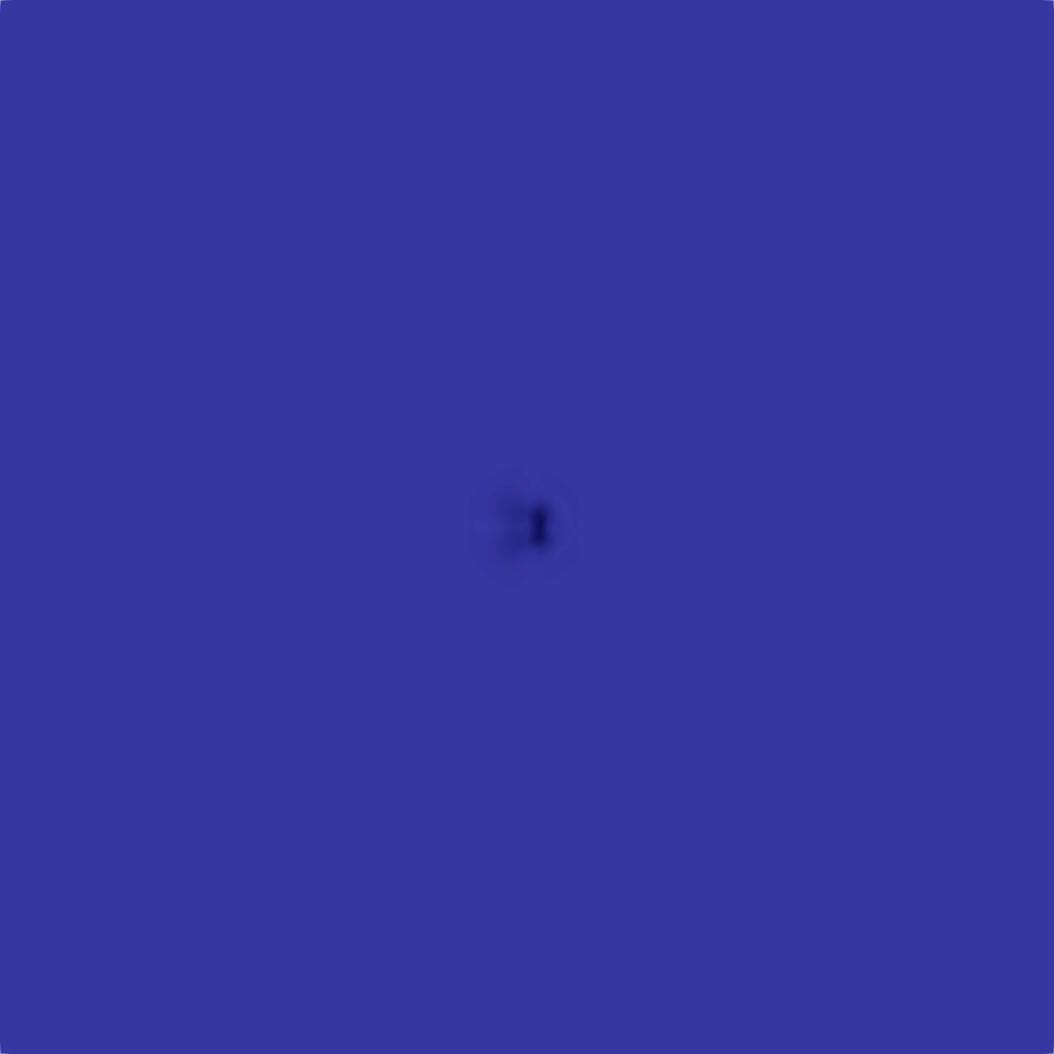}
    \label{fig:Carcione2018_v_01_fs_bidomain}
\end{subfigure}
\begin{subfigure}[b]{.33\textwidth}
    \centering
    \includegraphics[width=1\textwidth]{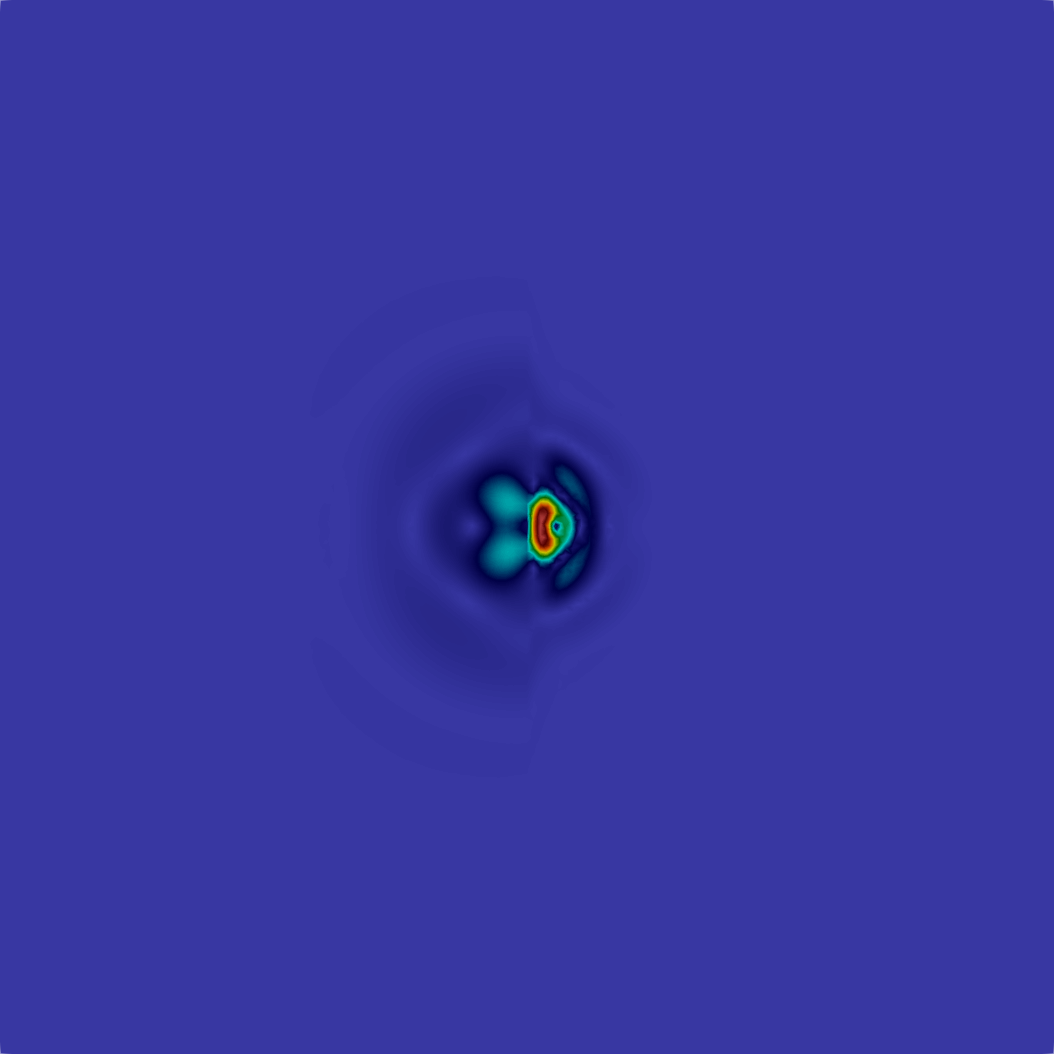}
    \label{fig:Carcione2018_v_03_fs_bidomain}
\end{subfigure}
\begin{subfigure}[b]{.33\textwidth}
    \centering
    \includegraphics[width=1\textwidth]{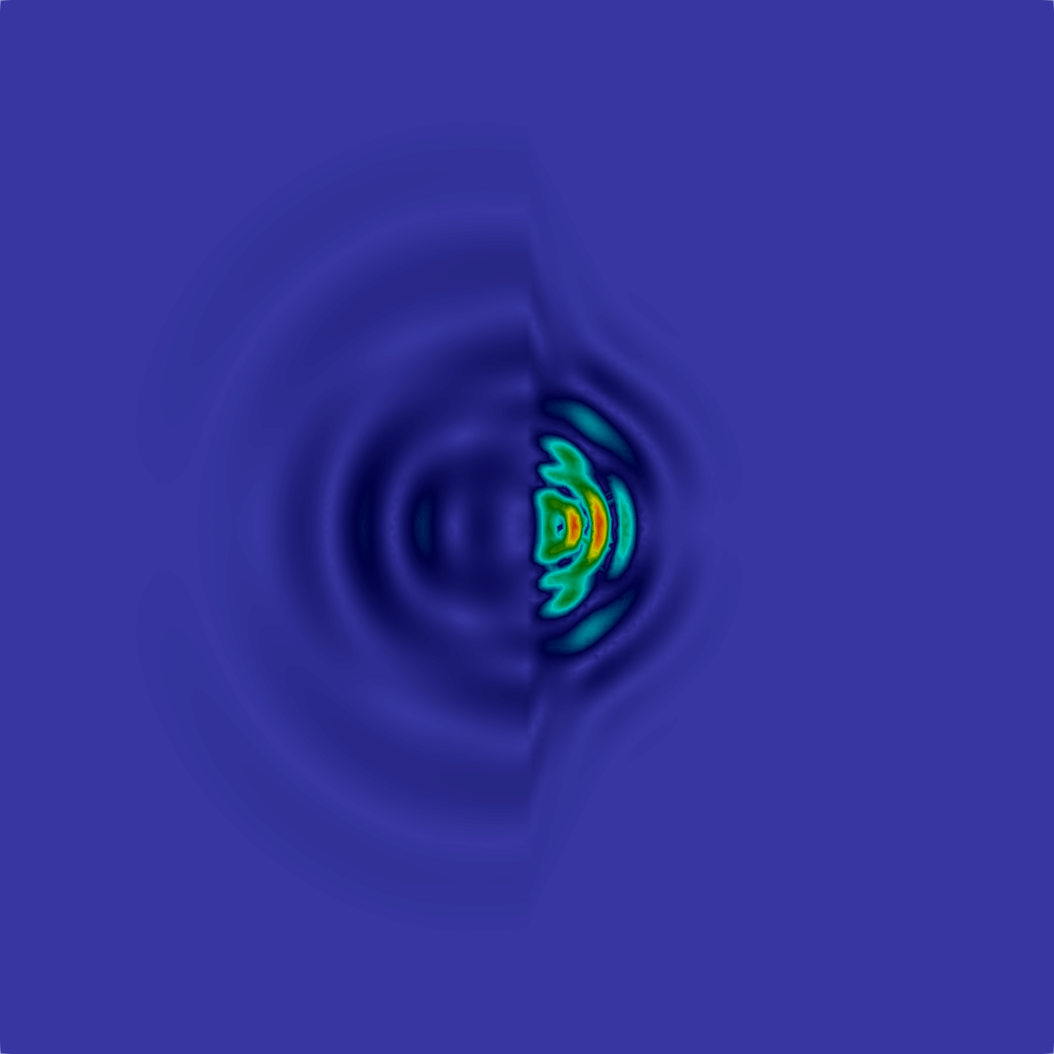}
    \label{fig:Carcione2018_v_05_fs_bidomain}
\end{subfigure}

\vspace{-0.4cm}
\centering
\includegraphics[width=0.3\textwidth]{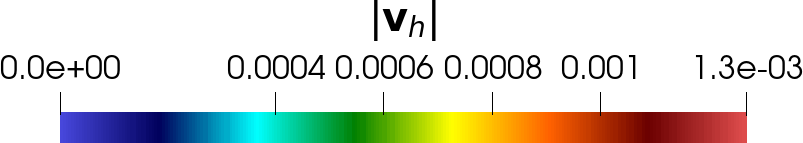}

\caption{Wave propagation in thermoelastic media with shear source term: computed velocity field $|\mathbf{v}_h|$ at the time instants $t=0.1 \si{\second}$ (left), $t=0.3 \si{\second}$ (center), $t=0.5 \si{\second}$ (right)}
\label{fig:Carcione2018_v_fs_bidomain}
\end{figure}

\begin{figure}[ht]
\begin{subfigure}[b]{.33\textwidth}
    \centering
    \includegraphics[width=1\textwidth]{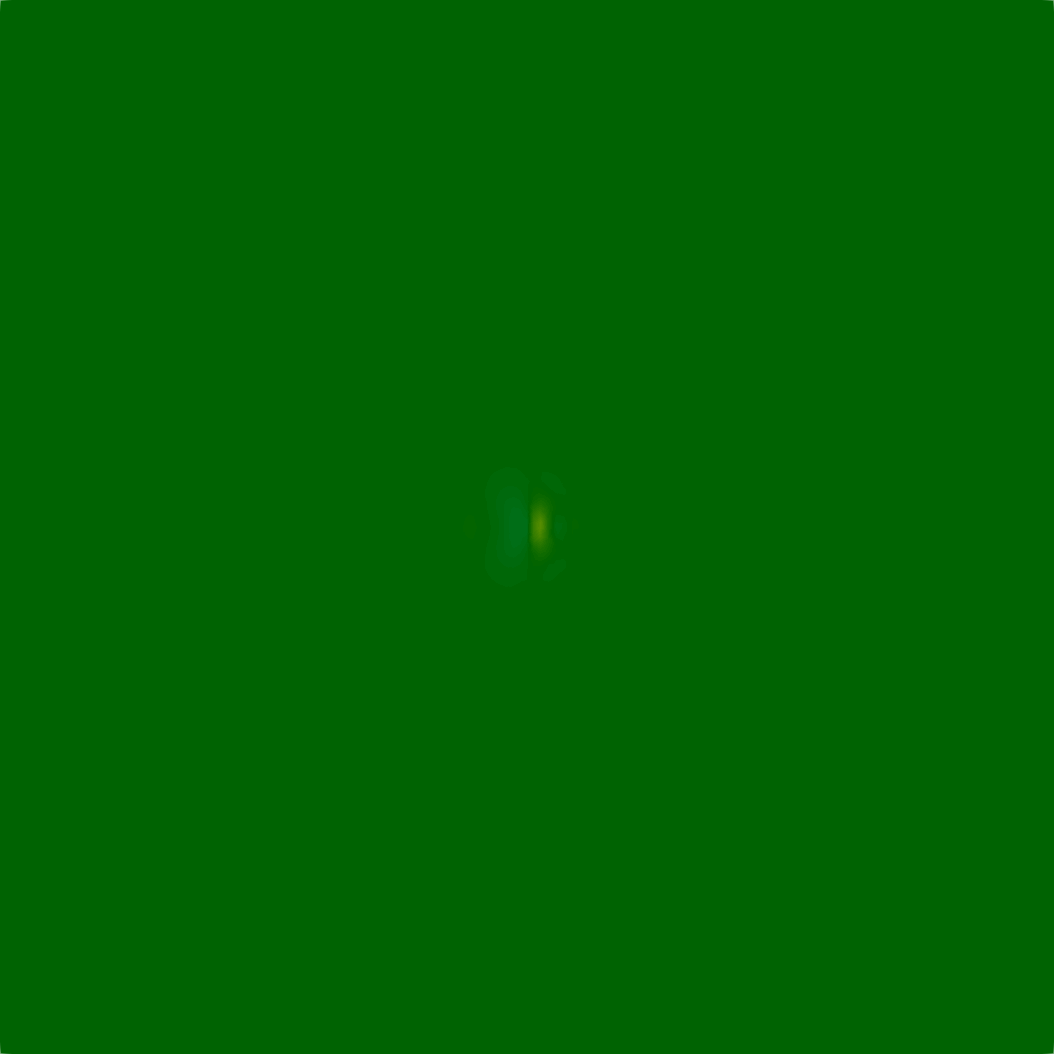}
    \label{fig:Carcione2018_vy_01_fs_bidomain}
\end{subfigure}
\begin{subfigure}[b]{.33\textwidth}
    \centering
    \includegraphics[width=1\textwidth]{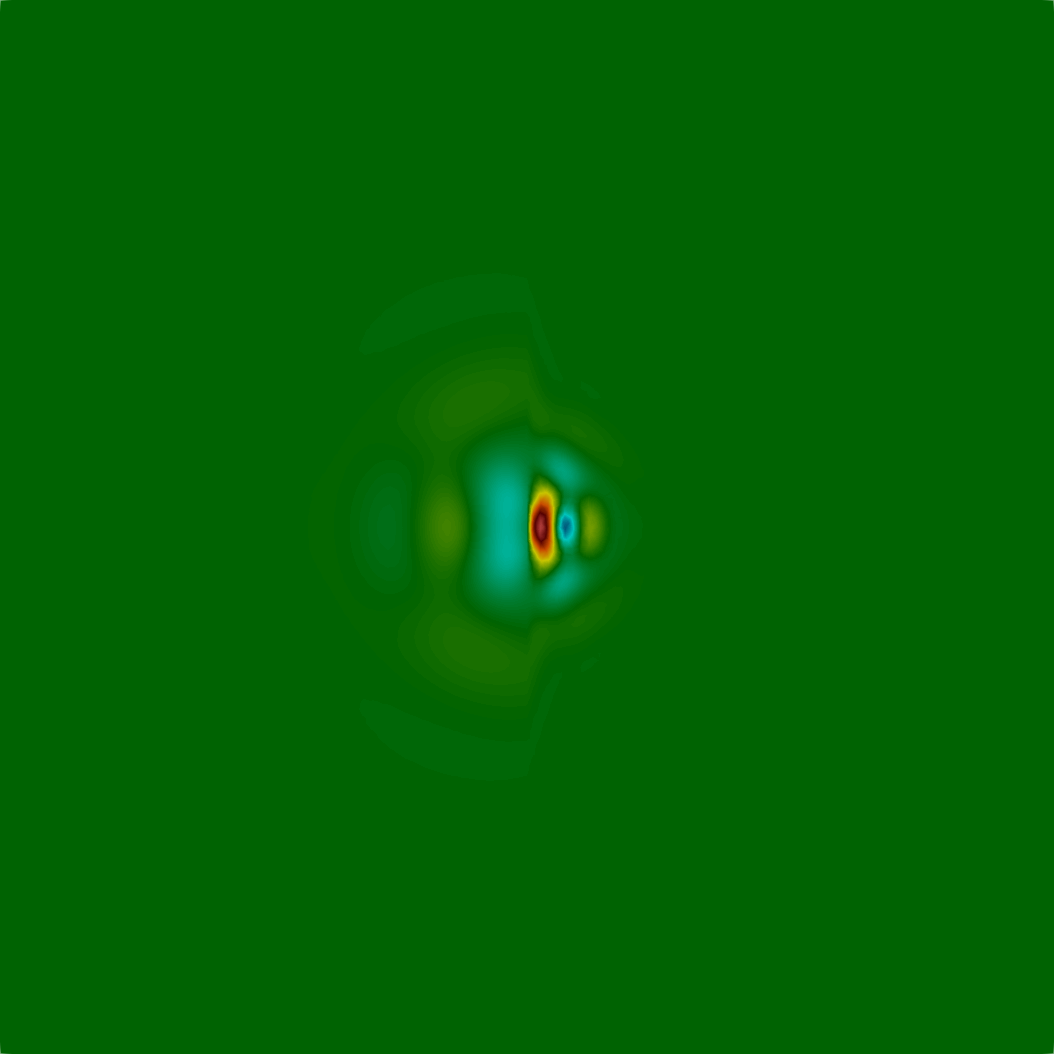}
    \label{fig:Carcione2018_vy_03_fs_bidomain}
\end{subfigure}
\begin{subfigure}[b]{.33\textwidth}
    \centering
    \includegraphics[width=1\textwidth]{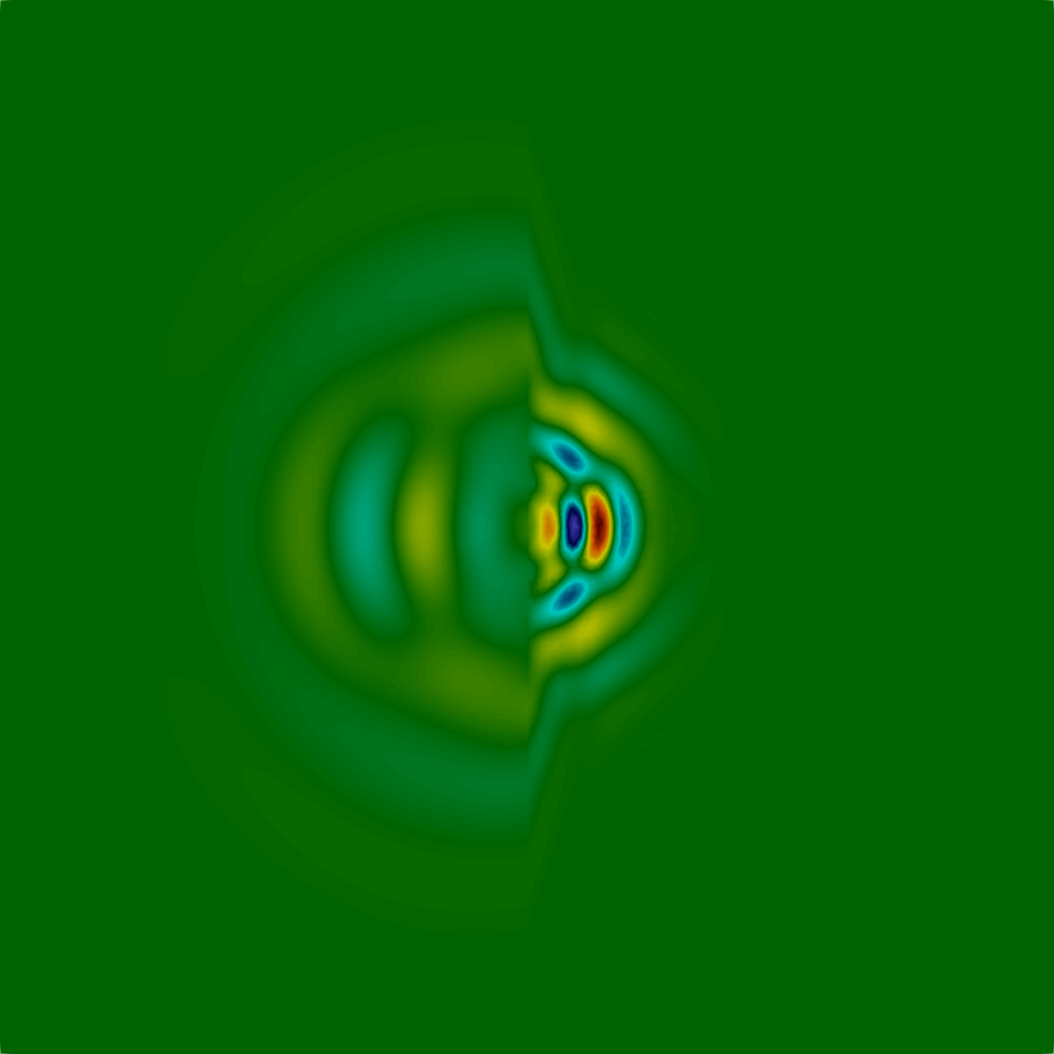}
    \label{fig:Carcione2018_vy_05_fs_bidomain}
\end{subfigure}

\vspace{-0.4cm}
\centering
\includegraphics[width=0.3\textwidth]{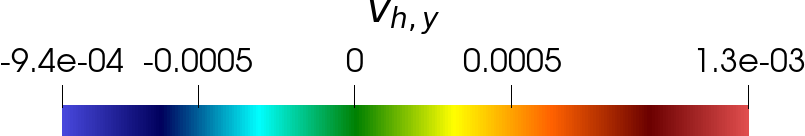}

\caption{Wave propagation in thermoelastic media with shear source term: computed vertical component of the velocity field $v_{h,y}$ at the time instants $t=0.1 \si{\second}$ (left), $t=0.3 \si{\second}$ (center), $t=0.5 \si{\second}$ (right).}
\label{fig:Carcione2018_vy_fs_bidomain}
\end{figure}

\begin{figure}[ht]
\begin{subfigure}[b]{.33\textwidth}
    \centering
    \includegraphics[width=1\textwidth]{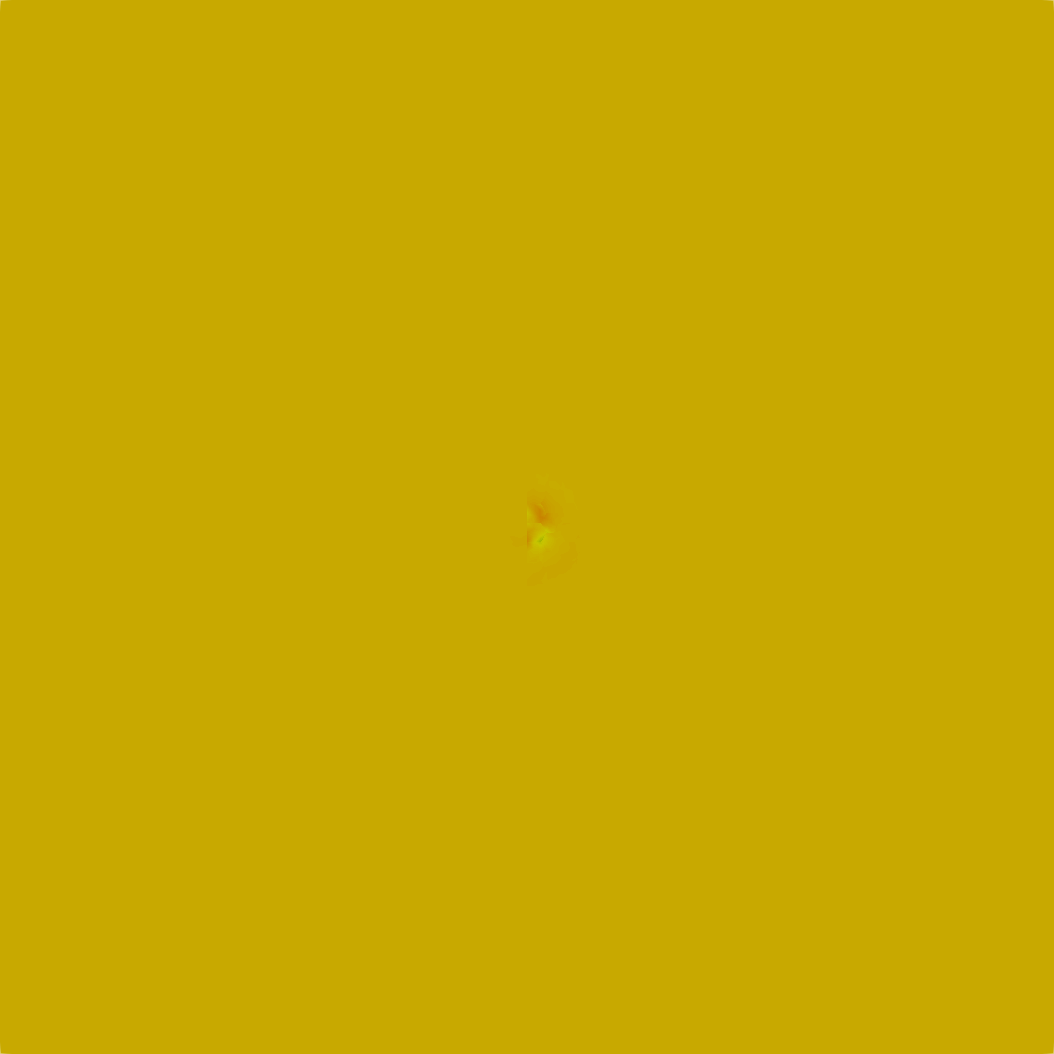}
    \label{fig:Carcione2018_T_01_fs_bidomain}
\end{subfigure}
\begin{subfigure}[b]{.33\textwidth}
    \centering
    \includegraphics[width=1\textwidth]{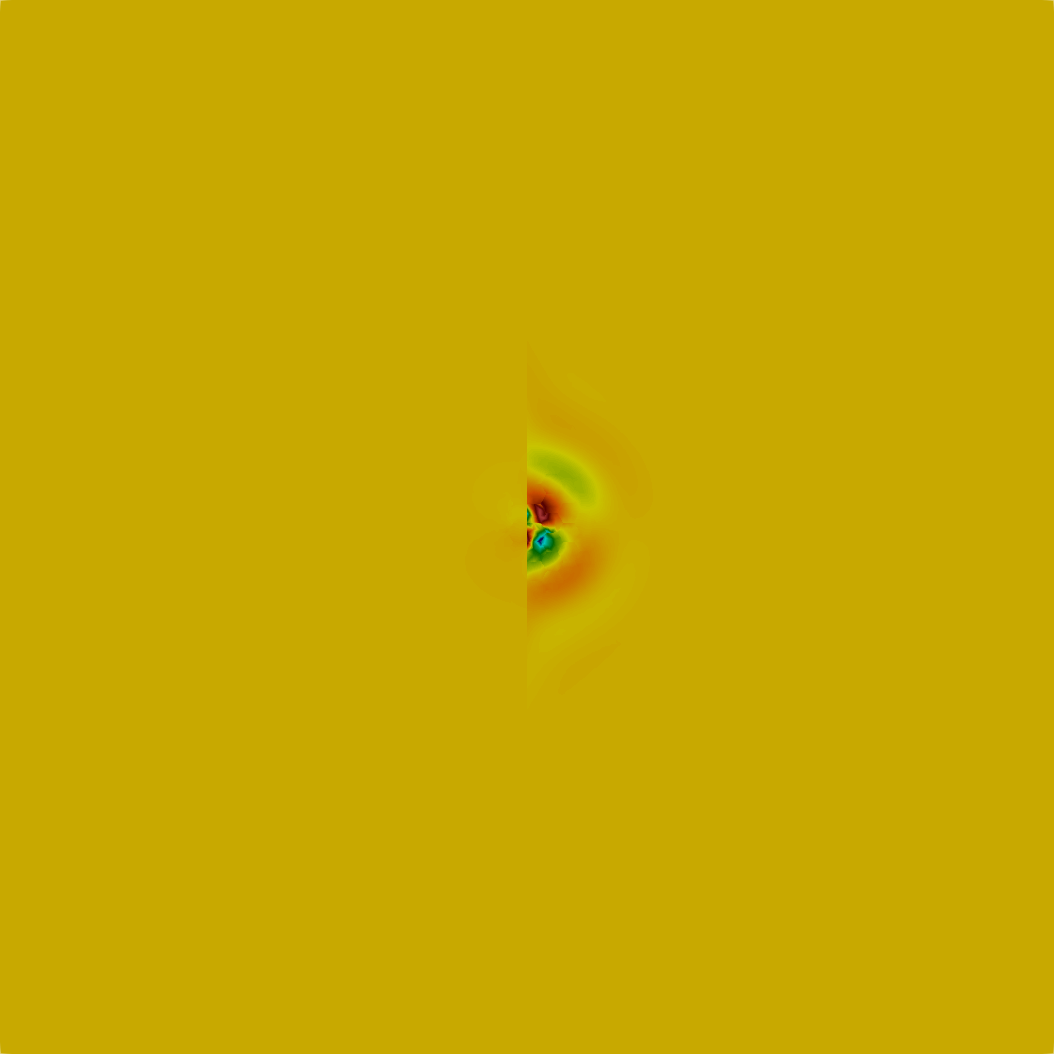}
    \label{fig:Carcione2018_T_03_fs_bidomain}
\end{subfigure}
\begin{subfigure}[b]{.33\textwidth}
    \centering
    \includegraphics[width=1\textwidth]{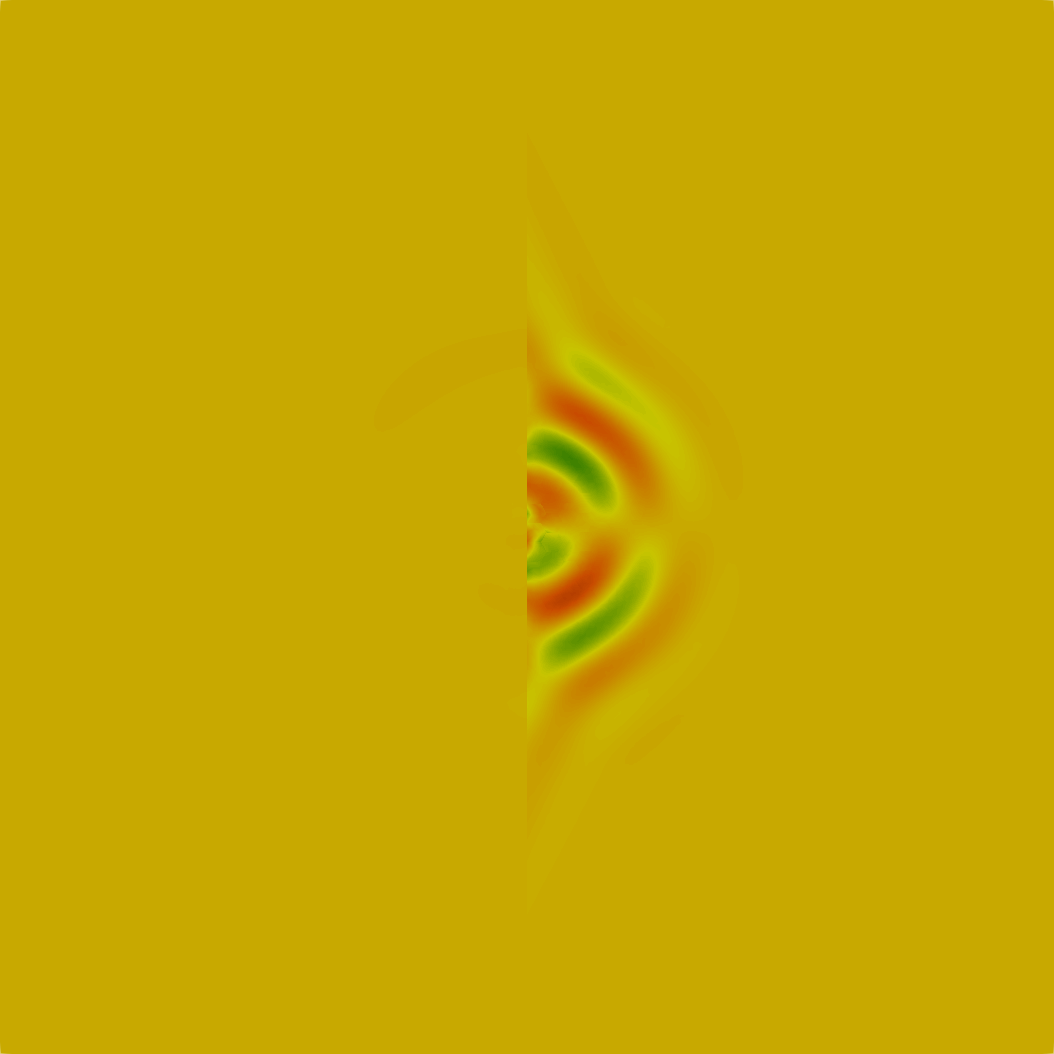}
    \label{fig:Carcione2018_T_05_fs_bidomain}
\end{subfigure}

\vspace{-0.4cm}
\centering
\includegraphics[width=0.3\textwidth]{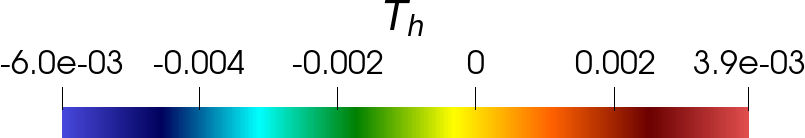}

\caption{Wave propagation in thermoelastic media with shear source term: computed temperature field $T_h$ at the time instants $t=0.1 \si{\second}$ (left), $t=0.3 \si{\second}$ (center), $t=0.5 \si{\second}$ (right).}
\label{fig:Carcione2018_T_fs_bidomain}
\end{figure}

In terms of the velocity field -- displayed in Figure~\ref{fig:Carcione2018_v_fs_bidomain} and Figure~\ref{fig:Carcione2018_vy_fs_bidomain} -- the main difference with respect to the homogeneous case is the presence of the head waves that connect the wavefronts of the two sub-domains. The head waves are particularly evident by looking at the vertical component of the velocity field (cf. last frame of Figure~\ref{fig:Carcione2018_vy_fs_bidomain}). For what concerns the temperature field $T_h$, whose behavior is shown in Figure~\ref{fig:Carcione2018_T_fs_bidomain}, we observe how the increasing value in the thermal conductivity affects the profile of the $T$-wave and, even in this case, we can observe the presence of wavelets that connect the wavefronts in the two layers. In general, by looking at the results for all three fields, we can observe that the heterogeneities in the media properties are well resolved by our method and that all the symmetric and anti-symmetric patterns observed in the \textit{Shear source term}-test case are preserved. However, the wave propagation in the two sub-layers is non-identical due to the different physical properties of the two media. We can observe a qualitatively good agreement with the results presented in \cite{Carcione2018} for the thermoelastic wave propagation and in \cite{Carcione2019,Bonetti2023} for the thermo-poroelastic wave propagation.

%% file: Sections/Spe10.tex
In this section, we consider the injection of a fluid in a heterogeneous poro-viscoelastic medium inspired by the SPE10 benchmark \cite{christie_tenth_2001}. This simulation aims to analyze the impact of the different possible modeling choices on the fluid flows in geophysical applications. We increase the initial permeability values by four orders to highlight the impact in high-permeable channels. 
\par
\begin{figure}[t]
\centering
\includegraphics[width=1\textwidth]{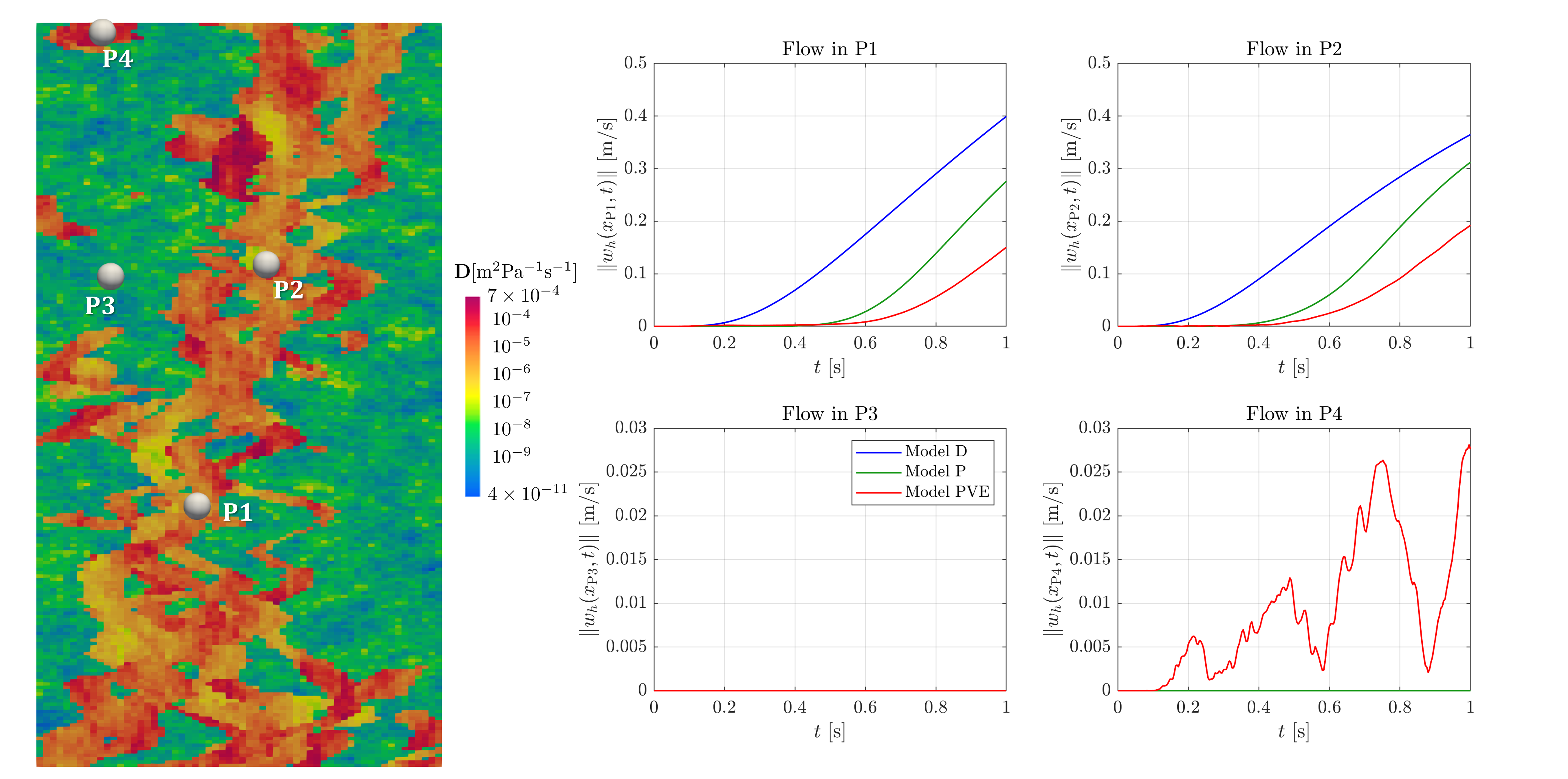}
\caption{Fluid flow in heterogeneous porous-viscoelastic medium: computed filtration velocity field  $\|\mathbf{w}_h\|$ measured in specific points of the domain with the three different models.}
\label{fig:Spe10flow}
\end{figure}

We consider as domain a horizontal slice (number 35) of the SPE10 benchmark $\Omega = (0\,\mathrm{m}, 366\,\mathrm{m})\times(0\,\mathrm{m}, 671\,\mathrm{m})$. The permeability values associated with the simulations can be observed in Figure \ref{fig:Spe10flow} (left). It can be observed that there is a channel of higher permeability that we expect to transport most of the fluid flow inside the domain. The simulation is associated with an injection of fluid in two sources with reabsorption in the middle of the domain. To reproduce this phenomenon, we set a forcing term: 
\begin{equation}
    g(\mathbf{x},t) = \dfrac{\tanh(5 t)}{10} \left(e^{-\frac{(x-190)^2+(y-550)^2}{500}}+e^{-\frac{(x-130)^2+(y-120)^2}{500}}-e^{-\frac{(x-175)^2+(y-360)^2}{500}}\right).
\end{equation}
With this choice, the injection and absorption increase in time until $t=0.5 \mathrm{s}$ when they reach a plateau. Concerning initial conditions, we set pressure, displacement, and velocity equal to $0$. Moreover, we set homogeneous Dirichlet boundary conditions for the displacement and homogeneous Neumann boundary conditions for the pressure. Concerning the space discretization, we use the cartesian mesh of 13200 elements provided by the SPE10 benchmark ($h=6.81\,\mathrm{m}$). This choice lets us maintain the initial benchmark's tensor $\mathbf{D}$ refinement level. We use the polynomial degree $\ell=2$ and a timestep $\Delta t=\num[exponent-product=\ensuremath{\cdot}, print-unity-mantissa=false]{4e-4}$.
\par
\begin{table}[H]
       \centering 
    \begin{tabular}{| r l | c | c | c | c |}
    \hline
    \multicolumn{2}{|c|}{\textbf{Parameter}} & \textbf{Model PVE} & \textbf{Model P} & \textbf{Model D} & \textbf{Reference} \T\B \\ \hline 
    $\mu$ & $[\si{\pascal}]$ & $10^9$ & - & - & \cite{Bonetti2023} \T\\ 
    $\lambda$ & $[\si{\pascal}]$ & $4\times 10^8$ & - & - & \cite{Bonetti2023} \T\\
    $\delta_1$ & $[\si{\second}]$ & $8\times10^{-5}$ & - & -  & \cite{Porenca2022} \T\\ 
    $\delta_2$ & $[\si{\second}]$ & $8\times10^{-5}$ & - & - & \cite{Porenca2022} \T\\
    $\gamma$ & $[-]$ & 1 & 0 & 0 & \cite{Bonetti2023} \T\\ 
    $d_0$ & $[\si{\per\pascal}]$ & $10^{-9}$ & $10^{-9}$ & $10^{-9}$ & \cite{Bonetti2023} \T\\
    $\tau_1$ & $[\si{\second}]$ & $\rho_f \phi^{-1} \mathbf{D}$ & $\rho_f \phi^{-1} \mathbf{D}$ & 0 & \cite{Nield2017} \T\\ 
    $\tau_2$ & $[\si{\second}]$ & $\rho_f \phi^{-1} \mathbf{D}$ & $\rho_f \phi^{-1} \mathbf{D}$ & 0 & \cite{Nield2017} \T\\
    $\rho_f$ & $[\si{\kilogram \per \metre \cubed}]$ & $1.025\times10^3$ & $1.025\times10^3$ & $1.025\times10^3$ & \cite{Fumagalli2024} \T\\
    $\phi$ & $[-]$ & 0.1 & 0.1 & 0.1 & \cite{Nield2017} \T\\
    \hline
    \end{tabular}
    \\[10pt]
    \caption{Fluid flow in heterogeneous poro-viscoelastic medium: physical parameters in the three simulation settings.}
    \label{tab:Spe10}
\end{table}
Starting from this general setting, which is common to all the simulations in this subsection, we construct three different models. In particular, we compare the proposed poro-viscoelastic model (PVE) with two porous media models: the first comprehensive of the acceleration term of the Darcy law (model P) and the second with a classical static Darcy law (model D). All these models can be obtained in our general framework by taking different values of the physical parameters. The values used in this simulation with the specific references are reported in Table \ref{tab:Spe10}. 
\begin{figure}[ht!]
\centering
\includegraphics[width=1\textwidth]{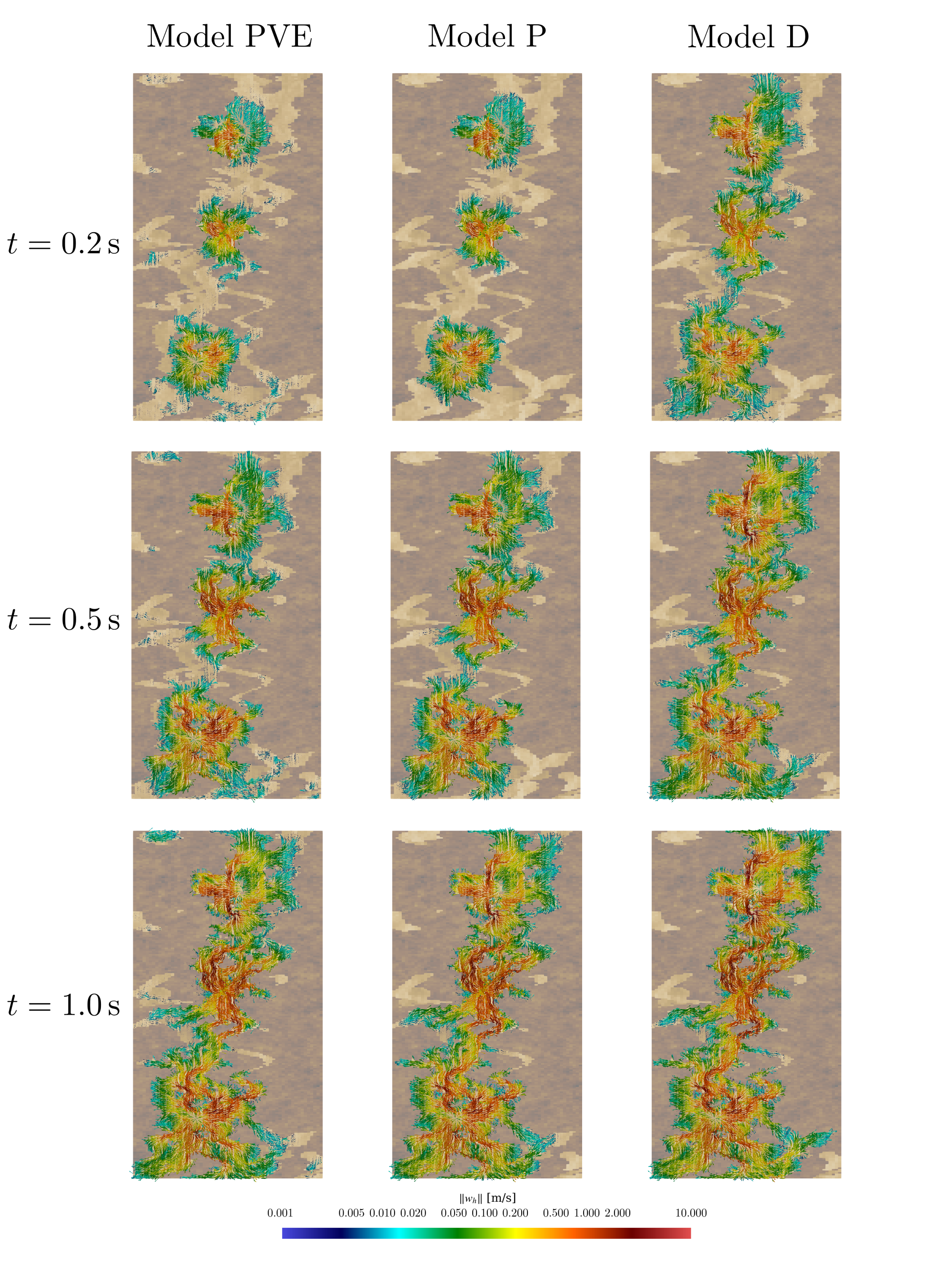}
\caption{Fluid flow in heterogeneous porous-viscoelastic medium: computed filtration velocity field  $\|\mathbf{w}_h\|$ at the time instants $t=0.2 \,\si{\second}$ (top), $t=0.5\,\si{\second}$ (middle), and $t=1.0 \,\si{\second}$ (bottom). Glyphs are not present wherever the flow is absent. Three different models are reported: model PVE (left), model P (middle), and model D (right).}
\label{fig:spe10wh}
\end{figure}
\par
In Figure \ref{fig:Spe10flow}, we report the result of the three models regarding the magnitude of the fluid filtration $\|\mathbf{w}_h\| = \| \mathbf{D} \nabla p_h\|$. In particular, we report the values measured at four specific points inside the domain. Concerning points P1$=(146.30\,\mathrm{m},236.22\,\mathrm{m})$ and P2$=(207.26\,\mathrm{m},452.63\,\mathrm{m})$, we can notice that the magnitude increases faster in model D than in the others. This is coherent with the modeling of the Darcy law, which does not consider fluid acceleration. This choice penalizes the fluid flow's continuity, reducing the previous timesteps' effect on the flow. Observing the model PVE, we highlight that the viscoelasticity reduces the magnitude of the flow in points P1 and P2, creating an additional delay in the fluid flow development. Concerning point P3$=(67.05\,\mathrm{m},441.69\,\mathrm{m})$, where the permeability is lower than in the other recording points, the fluid flow is approximately null in all the models (in the PVE model, the values are $\simeq 10^{-6} \mathrm{m/s}$). Finally, in P4$=(60.96\,\mathrm{m},661.42\,\mathrm{m})$, the elastic deformation of the soil induces a fluid flow in the model PVE, although the area is quite far from the channel in which the fluid injection is modeled.
\par
In Figure \ref{fig:spe10wh}, we report the complete fluid filtration field at three-time instants ($t=0.2,0.5,1.0\,\mathrm{s}$) and for the three different models. In particular, the \texttt{glyphs} are not reported wherever the flow is null. We can observe that the gap between models PVE and P against model D reduces during the simulation after a significant distance at the first times ($t=0.2\,\mathrm{s}$). In all the cases, the fluid flow is located along the high-permeability channel inside the domain. However, a significant difference is fluid flow induction inside non-connected high-permeability regions in model PVE, which increases over time. 
\par
\begin{figure}[ht]
\centering
\includegraphics[width=1\textwidth]{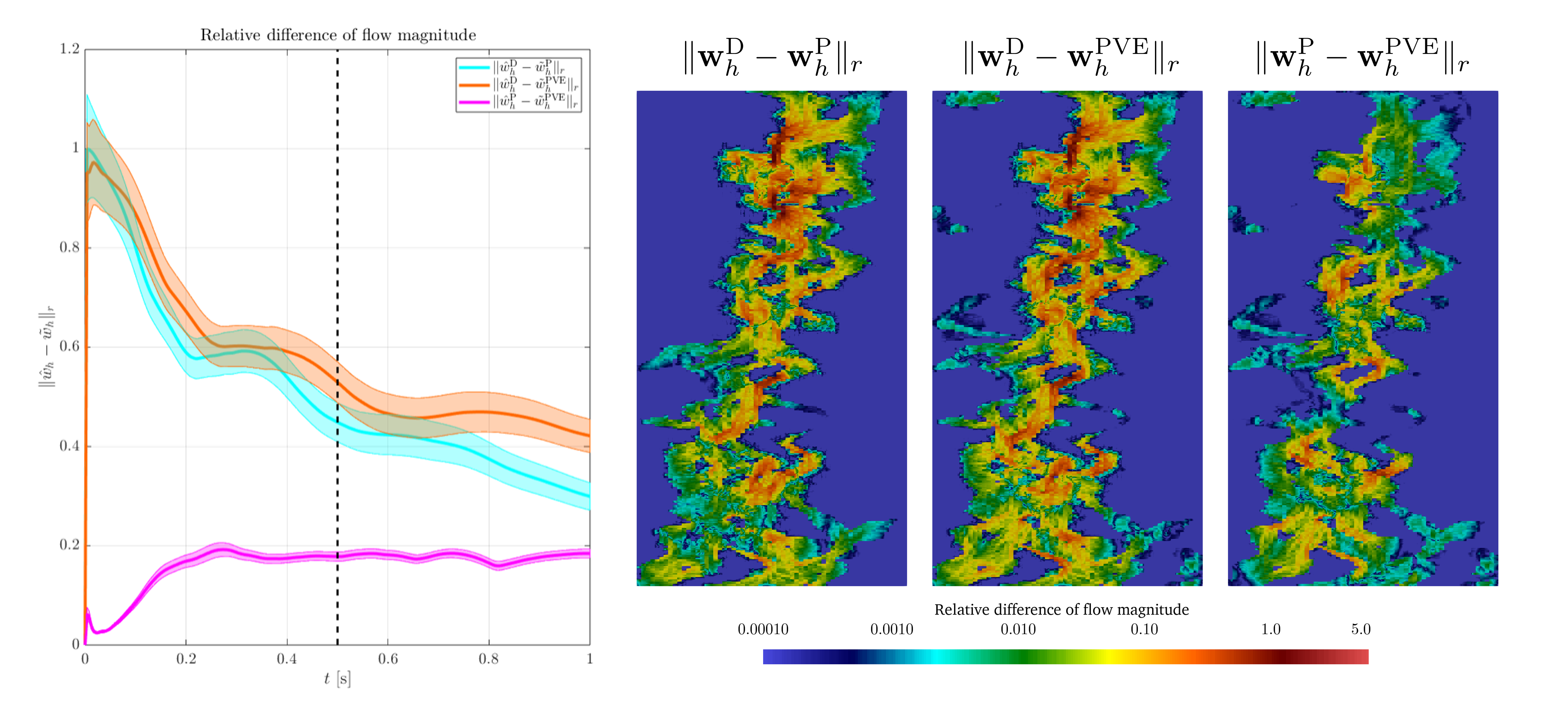}
\caption{Relative difference of flow magnitude: mean values and confidence bounds for each time (left) and detailed differences at time $0.5$ s.}
\label{fig:spe10whmean}
\end{figure}
To quantify the differences between the modeling choices, we construct a relative measure of filtration distance. Indeed, considering the models $*$ and $'$, we introduce for each mesh element $K$:
\begin{equation}
    \|\mathbf{w}_h^*-\mathbf{w}_h'\|_r = \dfrac{\|\mathbf{w}_h^*-\mathbf{w}_h'\|_K}{\|\mathbf{w}_h^\mathrm{D}\|},
\end{equation}
where the choice of rescaling with $\|\mathbf{w}_h^\mathrm{D}\|$ is guided by the fact this is the \textit{state-of-the-art} model in literature. In Figure \ref{fig:spe10whmean} (on the right), we report the relative differences in flow magnitude between the three models. As we can observe at time $t=0.5\,\mathrm{s}$, in points with high permeability values, the Darcy model overestimates the flow by more than $100\%$. In Figure \ref{fig:spe10whmean} (on the left), we report the mean in space (with confidence bounds computed with the standard deviation) for each time. We can observe that the difference between D and P, as well as between D and PVE, decreases in time. This is coherent with the additional inertia introduced by the terms $\tau_1$ and $\tau_2$. On the contrary, the difference between P and PVE is smaller than the others and almost constant in time. Indeed, after an initial time, when the models provide almost the same solution, the difference is highly guided by the capability of PVE to catch the flow in high permeability regions detached from the central channel. Finally, we underline that neglecting components of the mathematical model in choosing the Darcy law can sensibly affect the final solution.

%% file: Sections/Conclusion.tex
In this work, we have presented a PolyDG formulation for the thermo/poro-viscoelasticity problem. The model derivation remarks on which terms are more or less significant depending on the physical parameters of our test cases and, consequently, on the reference application. The stability estimate in the semi-discrete framework highlights the mild requirements on the physical parameters of the model problem. Moreover, it is general with respect to the choice of a fully-inertial or quasi-static model. For the semi-discrete framework, an a-priori error estimate is provided. Numerical simulations are performed to test the convergence and robustness properties of the proposed method. The results also show that the method provides a good approximation when considering limit cases in the ranges of physical parameters and presents some partial-superconvergence results. Finally, some benchmark verifications and physically sound test cases are presented, showing that the PolyDG discretization scheme can be appealing for real problem simulations with appropriate parameters. In particular, we used the Spe10 benchmark to highlight the impact of modeling choice on the fluid flow in heterogeneous media in geophysics.
\par
Further developments of the present work are possible. In particular, we mention the extension to other non-linear models for visco-elasticity. Secondly, to reduce the large computational cost required to deal with the fully coupled problem, two possible approaches are designing effective preconditioning techniques for the resulting system and developing effective splitting schemes. Given the improvements in the linear system resolution strategy, extending the implementation to the three-dimensional case is a further possible improvement. For this, generating and managing the computational mesh is also a point of development and interest. Finally, it would be interesting to analyze the complete thermo-poro-viscoelastic model, which has several applications. From a mathematical point of view, the coupling between the fluid and the temperature fields would provide interesting additional challenges, such as the presence of the nonlinear temperature convective term.

%% file: Sections/Appendix.tex
In this appendix, we show the computations to derive the Newmark-$\beta$-$\theta$ method for discretizing the equations with $\tau_1=0$. We use the Newmark-$\beta$ method for the first equation of \eqref{eq:semidiscr_alg_sys}. To this aim, we introduce the velocity vector $\mathbf{Z} = \dot{\mathbf{U}}$ and the acceleration vector $\mathbf{A} = \ddot{\mathbf{U}}$. Then, at each time-step $t^n$, we solve the following set of equations:
\begin{equation}
    \label{eq:newmark_displacement}
    \left\{
    \begin{aligned}
        & \begin{aligned}
            \bigg( \frac{1}{\beta_N \Delta t^2} \mathbf{M}_u & + \mathbf{A}_e + \mathbf{A}_{\text{div}} \bigg) \mathbf{U}^{n+1} + \left( \mathbf{A}_{e,\delta_1} + \mathbf{A}_{\text{div},\delta_2} \right) \mathbf{Z}^{n+1} - \mathbf{C}^T \boldsymbol{\Phi}^{n+1} = \mathbf{F}^{n+1} + \frac{1}{\beta_N \Delta t^2} \mathbf{M}_u \mathbf{U}^n \\
            & + \frac{1}{\beta_N \Delta t} \mathbf{M}_u \mathbf{Z}^n + \frac{1 - 2\beta_N}{2\beta_N} \mathbf{M}_u \mathbf{A}^n, 
        \end{aligned}\\
        & \mathbf{Z}^{n+1} = \mathbf{Z}^{n} + \Delta t \left( \gamma_N \mathbf{A}^{n+1} + (1-\gamma_N) \mathbf{A}^n \right), \\
        & \mathbf{A}^{n+1} = \frac{1}{\beta_N \Delta t^2}\left(\mathbf{U}^{n+1} - \mathbf{U}^{n}\right) - \frac{1}{\beta_N \Delta t} \mathbf{Z}^{n} + \frac{2\beta_N - 1}{2\beta_N} \mathbf{A}^n.
    \end{aligned}
    \right.
\end{equation}
In the first equation of \eqref{eq:newmark_displacement} we can plug the expressions for $\mathbf{Z}^{n+1}$, $\mathbf{A}^{n+1}$ to get:
\begin{equation}
    \label{eq:newmark_displacement2}
    \begin{aligned}
        & \bigg( \frac{1}{\beta_N \Delta t^2} \mathbf{M}_u + \mathbf{A}_u + \frac{\gamma_N}{\beta_N \Delta t} \mathbf{A}_{u, \delta} \bigg) \mathbf{U}^{n+1}  - \mathbf{C}^T \boldsymbol{\Phi}^{n+1} = \mathbf{F}^{n+1} + \left( \frac{1}{\beta_N \Delta t^2} \mathbf{M}_u + \frac{\gamma_N}{\beta_N \Delta t} \mathbf{A}_{u, \delta} \right) \mathbf{U}^{n} \\
        & + \bigg( \frac{1}{\beta_N \Delta t} \mathbf{M}_u - \frac{\beta_N - \gamma_N}{\beta_N} \mathbf{A}_{u, \delta} \bigg) \mathbf{Z}^n + \left( \frac{1 - 2\beta_N}{2\beta_N} \mathbf{M}_u - \frac{\Delta t (2\beta_N-\gamma_N)}{2 \beta_N} \mathbf{A}_{u, \delta} \right) \mathbf{A}^n, \\
    \end{aligned}
\end{equation}
where, for the sake of clarity, we have introduced the notation $\mathbf{A}_{u} = \mathbf{A}_{e} + \mathbf{A}_{\text{div}}$, $\mathbf{A}_{u, \delta} = \mathbf{A}_{e,\delta_1} + \mathbf{A}_{\text{div},\delta_2}$. Then, we couple \eqref{eq:newmark_displacement} with a $\theta$-method for the pressure equation. To obtain the formulation, we plug into the second equation in \eqref{eq:semidiscr_alg_sys} the definition of the velocity $\mathbf{Z}$ and of the acceleration $\mathbf{A}$ at time-continuous level. Thus, we get the following expression:
\begin{equation}
    \mathbf{M}_{\varphi} \, \dot{\boldsymbol{\Phi}}(t) + \mathbf{C}_{\tau_2} \, \mathbf{A}(t) + \mathbf{C} \, \mathbf{Z}(t) + \mathbf{A}_{\varphi} \boldsymbol{\Phi}(t) = \mathbf{G}(t), \quad t \in (0,T_f].
\end{equation}
Using this form of the pressure equation, the discretized equation reads:
\begin{equation}
    \begin{aligned}
        & \left( \frac{1}{\Delta t}\mathbf{M}_{\varphi} + \theta \mathbf{A}_{\varphi} \right) \boldsymbol{\Phi}^{n+1} + \theta \left(\mathbf{C}_{\tau_2} \, \mathbf{A}^{n+1} + \mathbf{C} \, \mathbf{Z}^{n+1} \right) = \left( \frac{1}{\Delta t}\mathbf{M}_{\varphi} - (1-\theta) \mathbf{A}_{\varphi} \right) \boldsymbol{\Phi}^{n} + \theta \mathbf{G}^{n+1} \\
        & + (1-\theta) \left( \mathbf{G}^{n} - \mathbf{C}_{\tau_2} \, \mathbf{A}^{n} - \mathbf{C} \, \mathbf{Z}^{n} \right) 
    \end{aligned}
\end{equation}
and by using the expressions for $\mathbf{Z}^{n+1}, \mathbf{A}^{n+1}$ in \eqref{eq:newmark_displacement} we obtain:
\begin{equation}
    \label{eq:theta_pressure}
    \begin{aligned}
        & \left( \frac{1}{\Delta t}\mathbf{M}_{\varphi} + \theta \mathbf{A}_{\varphi} \right) \boldsymbol{\Phi}^{n+1} + \mathbf{C}_u \, \mathbf{U}^{n+1} = \theta \, \mathbf{G}^{n+1} + \Tilde{\theta} \, \mathbf{G}^n + \left( \frac{1}{\Delta t}\mathbf{M}_{\varphi} - \Tilde{\theta} \mathbf{A}_{\varphi} \right) \boldsymbol{\Phi}^{n} \\
        & + \mathbf{C}_u \, \mathbf{U}^{n} + \mathbf{C}_z \, \mathbf{Z}^{n} + \mathbf{C}_a \, \mathbf{A}^{n}, 
    \end{aligned}
\end{equation}
where:
\begin{equation}
\begin{aligned}
    & \mathbf{C}_u = \frac{\theta}{\beta_N \Delta t^2}\left( \mathbf{C}_{\tau_2} + \Delta t \gamma \mathbf{C} \right), \\
    & \mathbf{C}_z = \left(\frac{\theta}{\beta_N \Delta t} \mathbf{C}_{\tau_2} + \frac{\theta \gamma_N - \beta_N}{\beta_N}\mathbf{C} \right), \\
    & \mathbf{C}_a = \left(\frac{\theta - 2\beta_N}{2\beta_N} \mathbf{C}_{\tau_2} + \frac{\theta\Delta t(\gamma - 2\beta_N)}{2\beta_N} \mathbf{C} \right),
\end{aligned}
\end{equation}
and $\Tilde{\theta} = 1 - \theta$. Thus, the final algebraic discretized formulation reads:
\begin{equation}
    \label{eq:full_time_discretization}
    \left\{
    \begin{aligned}
        & \begin{aligned}
            \bigg( \frac{1}{\beta_N \Delta t^2} \mathbf{M}_u + \mathbf{A}_u + \frac{\gamma_N}{\beta_N \Delta t} \mathbf{A}_{u, \delta} \bigg) \mathbf{U}^{n+1}  - \mathbf{C}^T \boldsymbol{\Phi}^{n+1} = \mathbf{F}^{n+1} + \left( \frac{1}{\beta_N \Delta t^2} \mathbf{M}_u + \frac{\gamma_N}{\beta_N \Delta t} \mathbf{A}_{u, \delta} \right) \mathbf{U}^{n} \\
            + \bigg( \frac{1}{\beta_N \Delta t} \mathbf{M}_u - \frac{\beta_N - \gamma_N}{\beta_N} \mathbf{A}_{u, \delta} \bigg) \mathbf{Z}^n + \left( \frac{1 - 2\beta_N}{2\beta_N} \mathbf{M}_u - \frac{\Delta t (2\beta_N-\gamma_N)}{2 \beta_N} \mathbf{A}_{u, \delta} \right) \mathbf{A}^n, \\
        \end{aligned} \\ 
        & \begin{aligned}
          \left( \frac{1}{\Delta t}\mathbf{M}_{\varphi} + \theta \mathbf{A}_{\varphi} \right) \boldsymbol{\Phi}^{n+1} + \mathbf{C}_u \, \mathbf{U}^{n+1} = \theta \, \mathbf{G}^{n+1} + \Tilde{\theta} \, \mathbf{G}^n + \left( \frac{1}{\Delta t}\mathbf{M}_{\varphi} - \Tilde{\theta} \mathbf{A}_{\varphi} \right) \boldsymbol{\Phi}^{n} \\
          + \mathbf{C}_u \, \mathbf{U}^{n}
          + \mathbf{C}_z \, \mathbf{Z}^n + \mathbf{C}_a \, \mathbf{A}^n,
        \end{aligned} \\ 
        & \mathbf{Z}^{n+1} = \mathbf{Z}^{n} + \Delta t \left( \gamma_N \mathbf{A}^{n+1} + (1-\gamma_N) \mathbf{A}^n \right),\\
        & \mathbf{A}^{n+1} = \frac{1}{\beta_N \Delta t^2}\left(\mathbf{U}^{n+1} - \mathbf{U}^{n}\right) - \frac{1}{\beta_N \Delta t} \mathbf{Z}^{n} + \frac{2\beta_N - 1}{2\beta_N} \mathbf{A}^n,
    \end{aligned}
    \right.
\end{equation}

%% file: Sections/AppendixConvergence.tex
In this appendix, we report the proof of Theorem \ref{thm:conv}. To start the proof of convergence, we define the following quantities 
\begin{equation}
    \begin{aligned}
        & \mathbf{e}_\mathrm{I}^{\mathbf{u}} = \mathbf{u}-\mathbf{u}_\mathrm{I}, \quad e_\mathrm{I}^\varphi = \varphi-\varphi_\mathrm{I}, \quad \mathbf{e}_\mathrm{I}^{\dot{\mathbf{u}}} = \dot{\mathbf{u}}-\dot{\mathbf{u}}_\mathrm{I}, \quad \mathbf{e}_\mathrm{I}^{\mathbf{W}} = \mathbf{W} - \mathbf{W}_\mathrm{I}, \quad \mathbf{e}_\mathrm{I}^{\Psi} = \Psi-\Psi_\mathrm{I}, \\
        & \mathbf{e}_h^u = \mathbf{u}_\mathrm{I}-\mathbf{u}_{h}, \quad e_h^\varphi = \varphi_\mathrm{I}-\varphi_\mathrm{h}, \quad \mathbf{e}_h^{\dot{\mathbf{u}}} = \dot{\mathbf{u}}_\mathrm{I} - \dot{\mathbf{u}}_h, \quad \mathbf{e}_h^{\mathbf{W}} = \mathbf{W}_\mathrm{I} - \mathbf{W}_h, \quad \mathbf{e}_h^{\Psi} = \Psi_\mathrm{I} - \Psi_h.
    \end{aligned}
\end{equation}
A main assumption of this analysis is the correspondence between the interpolants at time $t=0$ and the discrete initial conditions, such that: $$
\mathbf{e}^\mathbf{u}_h(0)=\mathbf{e}^{\dot{\mathbf{u}}}_h(0)=\mathbf{e}^\mathbf{W}_h(0)=\mathbf{0}, \qquad e^\varphi_h(0)=e^\Psi_h(0)=0. 
$$
After introducing the auxiliary variables and the additional problems as in the proof of Theorems \ref{thm:stab} and \ref{thm:stabdisc}, we can subtract them to obtain:
\begin{equation}
\label{eq:sys1_conv}
    \begin{aligned}
    \mathcal{M}_{u} & (\ddot{\mathbf{u}}-\ddot{\mathbf{u}}_h, \mathbf{v}_h) + \mathcal{M}_{\varphi}(\varphi- {\varphi}_h, \psi_h) + \mathcal{M}_{\varphi, \tau_1}(\dot{\varphi}-\dot{\varphi}_h, \psi_h) + \mathcal{A}_{e,h}(\mathbf{u}-\mathbf{u}_h, \mathbf{v}_h) \\
     & + \mathcal{A}_{e,\delta_1,h}(\dot{\mathbf{u}}-\dot{\mathbf{u}}_h, \mathbf{v}_h) + \mathcal{A}_{\text{div},h}(\mathbf{u}-\mathbf{u}_h, \mathbf{v}_h) + \mathcal{A}_{\text{div},\delta_2,h}(\dot{\mathbf{u}}-\dot{\mathbf{u}}_h, \mathbf{v}_h) \\
    & + \mathcal{A}_{\varphi,h}(\Psi-\Psi_h, \psi_h) - \mathcal{C}_h(\mathbf{v}_h, \varphi - \varphi_h) +\mathcal{C}_h(\mathbf{u}-\mathbf{u}_h, \psi_h) +\mathcal{C}_{\tau_2,h}(\dot{\mathbf{u}}-\dot{\mathbf{u}}_h, \psi_h) = 0,
    \end{aligned}
\end{equation}
and
\begin{equation}
\label{eq:sys2_conv}
    \begin{aligned}
    \mathcal{M}_{u} & (\dot{\mathbf{u}}-\dot{\mathbf{u}}_h, \mathbf{v}_h) + \mathcal{M}_{\varphi}(\varphi- {\varphi}_h, \psi_h) + \mathcal{M}_{\varphi, \tau_1}(\dot{\varphi}-\dot{\varphi}_h, \psi_h) + \mathcal{A}_{e,h}(\mathbf{W}-\mathbf{W}_h, \mathbf{v}_h) \\
     & + \mathcal{A}_{e,\delta_1,h}({\mathbf{u}}-{\mathbf{u}}_h, \mathbf{v}_h) + \mathcal{A}_{\text{div},h}(\mathbf{W}-\mathbf{W}_h, \mathbf{v}_h) + \mathcal{A}_{\text{div},\delta_2,h}({\mathbf{u}}-{\mathbf{u}}_h, \mathbf{v}_h) \\
    & + \mathcal{A}_{\varphi,h}(\Psi-\Psi_h, \psi_h) - \mathcal{C}_h(\mathbf{v}_h, \Psi - \Psi_h) +\mathcal{C}_h(\mathbf{u}-\mathbf{u}_h, \psi_h) +\mathcal{C}_{\tau_2,h}(\dot{\mathbf{u}}-\dot{\mathbf{u}}_h, \psi_h) = 0.
    \end{aligned}
\end{equation}
Now, in \eqref{eq:sys1_conv} we take $(\mathbf{v}_h, \psi_h) = (\tau_2 \mathbf{e}^{\dot{\mathbf{u}}}_h, e^\varphi_h)$ and multiply everything by $\tau_2$, then we exploit the coercivity and continuity properties of the bilinear forms:
\begin{equation}
    \begin{aligned}
    \dfrac{1}{2}\dfrac{\mathrm{d}}{\mathrm{d}t} & \Big(\|\tau_2\sqrt{\rho}\mathbf{e}^{\dot{\mathbf{u}}}_h\|^2 + \|\tau_2 \mathbf{e}^\mathbf{u}_h\|_\mathrm{dG,e}^2 + \|\sqrt{\tau_2 \tau_1 d_0} e^\varphi_h\|^2 + \|\sqrt{\tau_2} {e}^\Psi_h\|_\mathrm{dG,\varphi}^2\Big)  + \|\sqrt{\tau_2 d_0} e^\varphi_h\|^2 \\
    & + \|\tau_2 \mathbf{e}^{\dot{\mathbf{u}}}_h\|_\mathrm{dG,\delta}^2 
    +\mathcal{C}_h(\tau_2 \mathbf{e}^{{\mathbf{u}}}_h, e^\varphi_h) \leq \|\tau_2 \sqrt{\rho} \mathbf{e}^{\ddot{\mathbf{u}}}_\mathrm{I}\|\|\tau_2 \sqrt{\rho}\mathbf{e}^{\dot{\mathbf{u}}}_h\| + \|\sqrt{\tau_2 d_0} e^\varphi_\mathrm{I}\|\,\|\sqrt{\tau_2 d_0}e^\varphi_h\| \\
    &  + \|\sqrt{\tau_2 \tau_1 d_0} e^{\dot{\varphi}}_\mathrm{I}\|\, \|\sqrt{\tau_2 \tau_1 d_0} e^\varphi_h\| + \dfrac{\mathrm{d}}{\mathrm{d}t} \Big(\mathcal{A}_{e,h}(\tau_2\mathbf{e}^{{\mathbf{u}}}_\mathrm{I}, \tau_2 \mathbf{e}^{{\mathbf{u}}}_h) + \mathcal{A}_{\text{div},h}(\tau_2\mathbf{e}^{{\mathbf{u}}}_\mathrm{I}, \tau_2 \mathbf{e}^{{\mathbf{u}}}_h)\Big) \\ & +\tn \tau_2\mathbf{e}^{\dot{\mathbf{u}}}_\mathrm{I}\tn_\mathrm{dG,e} \| \tau_2 \mathbf{e}^{{\mathbf{u}}}_h\|_\mathrm{dG,e} +\tn \tau_2\mathbf{e}^{\dot{\mathbf{u}}}_\mathrm{I}\tn_\mathrm{dG,\delta} \| \tau_2 \mathbf{e}^{\dot{\mathbf{u}}}_h\|_\mathrm{dG,\delta} +  \dfrac{\mathrm{d}}{\mathrm{d}t} \mathcal{A}_{\varphi,h}(\tau_2e^\Psi_\mathrm{I}, e^\Psi_h) \\ & + \|\sqrt{\tau_2}e^\varphi_\mathrm{I}\tn_\mathrm{dG,\varphi} \|\sqrt{\tau_2}e^\Psi_h\|_\mathrm{dG,\varphi} - \mathcal{C}_h(\tau_2 \mathbf{e}^{\dot{\mathbf{u}}}_h, \tau_2e^\varphi_\mathrm{I}) + \mathcal{C}_h(\tau_2 \mathbf{e}^{{\mathbf{u}}}_\mathrm{I}, e^\varphi_h) 
+\mathcal{C}_{\tau_2,h}(\tau_2 \mathbf{e}^{\dot{\mathbf{u}}}_\mathrm{I}, e^\varphi_h), 
    \end{aligned}
\end{equation}
and in \eqref{eq:sys2_conv} we take $(\mathbf{v}_h, \psi_h) =(\mathbf{e}^\mathbf{u}_h,e^\Psi_h)$:
\begin{equation}
    \begin{aligned}
    \dfrac{1}{2}\dfrac{\mathrm{d}}{\mathrm{d}t} & \Big(\|\sqrt{\rho}\mathbf{e}^\mathbf{u}_h\|^2 + \| \sqrt{d_0} e^\Psi_h \|^2 + \|\mathbf{e}^\mathbf{W}_h\|_\mathrm{dG,e}^2 \Big) + \|\mathbf{e}^\mathbf{u}_h\|_\mathrm{dG,\delta}^2 + \|e^\Psi_h\|_\mathrm{dG,\varphi}^2 - \mathcal{C}_{\tau_2,h}(\mathbf{e}^{\mathbf{u}}_h, e^\varphi_h) \\
    &  + \dfrac{\mathrm{d}}{\mathrm{d}t} \mathcal{C}_{\tau_2,h}(\mathbf{e}^{{\mathbf{u}}}_h, e^\Psi_h) \leq \|\sqrt{\rho}\mathbf{e}^{\dot{\mathbf{u}}}_\mathrm{I}\|\,\|\sqrt{\rho} \mathbf{e}^\mathbf{u}_h\| + \left(\|\sqrt{d_0}{e}^\varphi_\mathrm{I}\| + \|\sqrt{d_0}\tau_1{e}^{\dot{\varphi}}_\mathrm{I}\|\right)\|\sqrt{d_0}e^\Psi_h\| \\ & - \dfrac{\mathrm{d}}{\mathrm{d}t}\left(\mathcal{A}_{e,h}(\mathbf{e}^\mathbf{W}_h, \mathbf{e}^\mathbf{W}_h)+\mathcal{A}_{\mathrm{div},h}(\mathbf{e}^\mathbf{W}_h, \mathbf{e}^\mathbf{W}_h) \right) + \tn\mathbf{e}^\mathbf{u}_\mathrm{I}\tn_\mathrm{dG,e} \|\mathbf{e}^\mathbf{W}_h\|_\mathrm{dG,e} \\ & +  \tn\mathbf{e}^\mathbf{u}_\mathrm{I}\tn_\mathrm{dG,\delta} \|\mathbf{e}^\mathbf{u}_h\|_\mathrm{dG,\delta} + \tn e^\Psi_\mathrm{I} \tn_\mathrm{dG,\varphi} \| e^\Psi_h\|_\mathrm{dG,\varphi} - \mathcal{C}_h(\mathbf{e}^\mathbf{u}_h, e^\Psi_\mathrm{I}) \\ & + \mathcal{C}_h(\mathbf{e}^\mathbf{u}_\mathrm{I}, e^\Psi_h) - \mathcal{C}_{\tau_2,h}(\mathbf{e}^{\mathbf{u}}_\mathrm{I}, e^\varphi_h) + \dfrac{\mathrm{d}}{\mathrm{d}t} \mathcal{C}_{\tau_2,h}(\mathbf{e}^{{\mathbf{u}}}_\mathrm{I}, e^\Psi_h) - \mathcal{M}_{\varphi, \tau_1}({e}^{\dot{\varphi}}_h, e^\Psi_h).
    \end{aligned}
\end{equation}
Adapting the proof of the stability (see Theorem \ref{thm:stab}), we will adopt Poincarè and Korn's inequalities to bound some terms in their discrete version (cf. \cite{dipietro2020,BottiDiPietro2020}). To handle the last term on the left-hand side, we derive a new estimate from \eqref{eq:sys2_conv} by $(\tau_2\mathbf{e}^\mathbf{u}_h,0)$ and we obtain:
\begin{equation}
\label{eq:additional_est_c}
    \begin{aligned}
    \dfrac{1}{2}\dfrac{\mathrm{d}}{\mathrm{d}t} & \left(\| \sqrt{\tau_2\rho} \mathbf{e}^\mathbf{u}_h\|^2 + \|\sqrt{\tau_2}\mathbf{e}^\mathbf{W}_h\|_\mathrm{dG,e}^2 \right) + \|\sqrt{\tau_2}\mathbf{e}^\mathbf{u}_h\|_\mathrm{dG,\delta}^2 \leq \mathcal{C}_h(\tau_2\mathbf{e}^\mathbf{u}_h, e^\Psi_h) - \mathcal{C}_h(\tau_2\mathbf{e}^\mathbf{u}_h, e^\Psi_\mathrm{I}) \\
    + & \|\sqrt{\tau_2\rho}\mathbf{e}^{\dot{\mathbf{u}}}_\mathrm{I}\|\;\|\sqrt{\tau_2\rho}\mathbf{e}^\mathbf{u}_h\| + \tn\mathbf{e}^\mathbf{W}_\mathrm{I}\tn_\mathrm{dG,e} \|\tau_2\mathbf{e}^\mathbf{u}_h\|_\mathrm{dG,e}
    + \tn\sqrt{\tau_2}\mathbf{e}^\mathbf{u}_\mathrm{I}\tn_\mathrm{dG,\delta} \|\sqrt{\tau_2}\mathbf{e}^\mathbf{u}_h\|_\mathrm{dG,\delta}.
    \end{aligned}
\end{equation}
Summing up the two equations and integrating in time, we obtain the following estimate, where we neglect the dependence on the physical parameters (using also the equivalence between the two norms $\|\cdot\|_{\mathrm{dG,e}}$ and $\|\cdot\|_{\mathrm{dG,\delta}}$):
\begin{equation}
    \begin{aligned} \dfrac{\mathrm{d}}{\mathrm{d}t} & \left(\|\mathbf{e}^\mathbf{u}_h\|^2 + \|e^\Psi_h\|^2 + \|\mathbf{e}^\mathbf{W}_h\|_\mathrm{dG,e}^2 \right) + \|\mathbf{e}^{\dot{\mathbf{u}}}_h\|^2 + \|\mathbf{e}^\mathbf{u}_h\|^2 + \| e^\varphi_h\|^2 +  \|e^\Psi_h \|^2 + \|\mathbf{e}^\mathbf{u}_h\|_\mathrm{dG,e}^2 \\
    + & \|\mathbf{e}^\mathbf{W}_h\|_\mathrm{dG,e}^2
    + \|e^\Psi_h\|_\mathrm{dG,\varphi}^2 + \int_0^t \left(\|e^\varphi_h\|^2 + \|\mathbf{e}^\mathbf{u}_h\|_\mathrm{dG,e}^2 + \|\mathbf{e}^{\dot{\mathbf{u}}}_h\|_\mathrm{dG,e}^2  + \|e^\Psi_h\|_\mathrm{dG,\varphi}^2 \right) \\ 
    \lesssim & \; \|\mathbf{e}^{\dot{\mathbf{u}}}_\mathrm{I}\|^2 +\|e^\Psi_\mathrm{I}\|^2 + \tn \mathbf{e}^{{\mathbf{u}}}_\mathrm{I}\tn_\mathrm{dG,e}^2 +\tn {e}^\Psi_\mathrm{I}\tn_\mathrm{dG,\varphi}^2 + \tn\mathbf{e}^\mathbf{W}_\mathrm{I}\tn_\mathrm{dG,e}^2 +  \|\mathbf{e}^\mathbf{u}_h\|^2 \\
    + & \int_0^t\left(\|\mathbf{e}^{\ddot{\mathbf{u}}}_\mathrm{I}\|^2 + \|\mathbf{e}^{\dot{\mathbf{u}}}_\mathrm{I}\|^2 + \tn\mathbf{e}^{\dot{\mathbf{u}}}_\mathrm{I}\tn_\mathrm{dG,e}^2 + \tn\mathbf{e}^\mathbf{u}_\mathrm{I}\tn_\mathrm{dG,e}^2 + \|\mathbf{e}^{\dot{\mathbf{u}}}_h\|^2 + \| \mathbf{e}^{{\mathbf{u}}}_h\|_\mathrm{dG,e}^2 + \|\mathbf{e}^\mathbf{u}_h\|^2 + \|\mathbf{e}^\mathbf{W}_h\|_\mathrm{dG,e}^2\right) \\
    + & \int_0^t\left( \|e^\Psi_\mathrm{I}\|^2 + \|e^\varphi_\mathrm{I}\|^2 + \|e^{\dot{\varphi}}_\mathrm{I}\|^2 + \|e^\varphi_\mathrm{I}\tn_\mathrm{dG,\varphi}^2 +\tn e^\Psi_\mathrm{I} \tn_\mathrm{dG,\varphi}^2 +\|e^\varphi_h\|^2 + \|e^\Psi_h\|_\mathrm{dG,\varphi}^2 + \|e^\Psi_h\|^2 \right).
    \end{aligned}
\end{equation}
Finally, by using Gr\"{o}nwall lemma \cite{Quarteroni2017}, integrating in time and neglecting the errors connected with the additional variables at the left-hand side, we obtain that:
\begin{equation}
\label{eq:conv_final_step}
    \begin{aligned} \|\mathbf{e}^\mathbf{u}_h\|^2 + & \int_0^t \left(\|\mathbf{e}^{\dot{\mathbf{u}}}_h\|^2 + \|\mathbf{e}^\mathbf{u}_h\|^2 +  \|\mathbf{e}^\mathbf{u}_h\|_\mathrm{dG,e}^2 \right) 
\lesssim (1 + t e^t) \int_0^t\Big(\|e_h^\varphi\|^2 + \|\mathbf{e}^{\ddot{\mathbf{u}}}_\mathrm{I}\|^2 
    + \|\mathbf{e}^{\dot{\mathbf{u}}}_\mathrm{I}\|^2 
    + \tn\mathbf{e}^{\dot{\mathbf{u}}}_\mathrm{I}\tn_\mathrm{dG,e}^2 \\ 
    + & \tn\mathbf{e}^\mathbf{u}_\mathrm{I}\tn_\mathrm{dG,e}^2
    + \tn\mathbf{e}^\mathbf{W}_\mathrm{I}\tn_\mathrm{dG,e}^2 
    + \|e^\Psi_\mathrm{I}\|^2 + \|e^\varphi_\mathrm{I}\|^2 
    + \|e^{\dot{\varphi}}_\mathrm{I}\|^2 
    + \|e^\varphi_\mathrm{I}\tn_\mathrm{dG,\varphi}^2 
    + \tn e^\Psi_\mathrm{I} \tn_\mathrm{dG,\varphi}^2 \Big).
    \end{aligned}
\end{equation}
The thesis follows by bounding the right-hand side of Equation~\eqref{eq:conv_final_step} using the interpolation estimates of Proposition \ref{prop:interpolant}.